\documentclass[11pt]{amsart}
\usepackage[a4paper,hmargin=28mm,vmargin=25mm]{geometry}
\usepackage{amsmath,amssymb,amsthm,mathtools}
\usepackage{microtype}
\usepackage{longtable}
\usepackage{ytableau}
\ytableausetup{centertableaux}
\usepackage{hyperref}
\hypersetup{hidelinks}

\newtheorem{theorem}{Theorem}[section]
\newtheorem{mainresult}{Theorem}

\newtheorem{proposition}[theorem]{Proposition}
\newtheorem{lemma}[theorem]{Lemma}
\newtheorem{corollary}[theorem]{Corollary}
\theoremstyle{remark}

\newtheorem{example}[theorem]{Example}
\newcommand{\Gr}{\operatorname{Gr}}
\newcommand{\SSYT}{\operatorname{SSYT}}
\newcommand{\bbk}{\mathbb K}

\numberwithin{equation}{section}

\title{Quantum Grassmannians and Quantum Grothendieck Rings in Type A}
\author{Jian-Rong Li}
\address{Jian-Rong Li: Faculty of Mathematics, University of Seville, Calle Tarfia s/n, 41012 Seville, Spain}
\email{lijr07@gmail.com}
\date{}
\begin{document}
\begin{abstract}
Hernandez and Leclerc established an isomorphism between
the Grothendieck ring of the category $\mathcal C_\ell$
in type $A$ and a quotient of the homogeneous coordinate
ring of a Grassmannian.
Chang, Duan, Fraser, and the author further studied
this isomorphism using semistandard Young tableaux
and obtained explicit formulas for the basis of this quotient
corresponding to simple modules, where the non-wrapping
consecutive Pl\"ucker coordinates are set to $1$.
The passage to the quantum setting is nontrivial since
the non-wrapping consecutive quantum Pl\"ucker coordinates
are not central in general and cannot simply be specialized
to $1$ as in the classical case.
We construct a cocycle twist of the quantum Grassmannian
under which these coordinates become central and specialize
them to suitable powers of the deformation parameter.
We establish an explicit isomorphism between the resulting
quotient and the quantum Grothendieck ring of $\mathcal C_\ell$
and determine the images of all dual canonical basis elements.
Combining this isomorphism with the basis transition formula
of Brundan, we obtain
an explicit formula for the $(q,t)$-characters of arbitrary
simple modules in $\mathcal C_\ell$ in terms of fundamental
$(q,t)$-characters and Kazhdan--Lusztig polynomials of
symmetric groups.
We also construct a normalized quantum cyclic shift
which acts on the dual canonical basis of the quantum Grassmannian by
Sch\"utzenberger promotion and describe the images of
the shifted basis elements in the central quotient.
After localization at the images of the wrapping consecutive
quantum Pl\"ucker coordinates, a frozen-variable correction
of the normalized quantum cyclic shift yields quantum cluster
quasi-automorphisms of the localized quotient and the
localized quantum Grothendieck ring, which permute
their respective localized dual canonical and canonical bases.
\end{abstract}

\subjclass[2020]{Primary 17B37; Secondary 13F60, 14M15, 20C08.}

\keywords{Quantum Grassmannians, quantum affine algebras,
quantum Grothendieck rings, dual canonical bases,
$(q,t)$-characters, Kazhdan--Lusztig polynomials,
cocycle twists, cluster quasi-automorphisms}

\maketitle

\setcounter{tocdepth}{1}
\tableofcontents

\section{Introduction}

Hernandez and Leclerc~\cite{HL10a,HL2015} developed
cluster algebra structures on the Grothendieck rings of
certain monoidal subcategories of finite-dimensional
representations of quantum affine algebras.
Another important class of cluster algebras arises from
the homogeneous coordinate rings of Grassmannians,
as shown by Scott~\cite{Scott2006}.
In type $A$, these two classes of cluster algebras are
related through an isomorphism between the Grothendieck
ring of a Hernandez--Leclerc category and a quotient
of a Grassmannian coordinate ring~\cite{HL10a}.
Chang, Duan, Fraser, and the author~\cite{CDFL20}
described this isomorphism in terms of semistandard
Young tableaux and investigated its relation to the
dual canonical basis \cite{Lusztig1993}
(equivalently, the upper global basis~\cite{Kashiwara1993}).
In this paper, we establish a quantum analogue of this
correspondence, identify the respective canonical bases,
and construct quantum cluster quasi-automorphisms
compatible with Sch\"utzenberger promotion~\cite{Schutzenberger1972}.

\subsection{The classical correspondence}

Let $k\geq 2$, $\ell\geq 0$, $n=k+\ell+1$, and
$m=\ell+1$. Let $\Gr(k,n)$ be the
Grassmannian of $k$-dimensional subspaces of $\mathbb C^n$,
and let $\mathbb C[\Gr(k,n)]$ be its
homogeneous coordinate ring with respect to the
Pl\"ucker embedding.
We write $[n]=\{1,\ldots,n\}$ and
$[a,b]=\{a,a+1,\ldots,b\}$ for integers $a\leq b$.
We use $\mathbb N=\mathbb Z_{\geq0}$.

Let $\mathsf I=\{1,\ldots,k-1\}$ be the vertex set of
the Dynkin diagram of type $A_{k-1}$.
A height function on this diagram is a function
$\xi:\mathsf I\to\mathbb Z$ such that
$|\xi_i-\xi_j|=1$ whenever $|i-j|=1$.
Following~\cite{CDFL20}, we fix the height function
$\xi_i=i-2$ throughout the paper and set
\[
 \mathcal X
 =
 \{(i,\xi_i-2r):i\in\mathsf I,\ 0\leq r\leq\ell\}
 =
 \{(i,i-2b):i\in\mathsf I,\ 1\leq b\leq m\}.
\]
We assume that $q\in\mathbb C^\times$ is not a root of
unity and consider finite-dimensional type~$1$ representations of the quantum affine algebra
$U_q(\widehat{\mathfrak{sl}_k})$.
For $(i,p)\in\mathsf I\times\mathbb Z$, let $Y_{i,p}$
denote the variable corresponding to the fundamental
$\ell$-weight at node $i$ with spectral parameter $q^p$.
A Laurent monomial in the variables $Y_{i,p}$ is called dominant
if all its exponents are nonnegative.
Let $\mathcal P_\ell^+$ be the monoid of dominant monomials
in the variables $Y_{i,p}$ with $(i,p)\in\mathcal X$.
We write $L(M)$ for the simple module with highest
$\ell$-weight monomial $M$.
The Hernandez--Leclerc category $\mathcal C_\ell$ is
the full subcategory of finite-dimensional type~$1$ $U_q(\widehat{\mathfrak{sl}_k})$-modules
whose simple composition factors have highest
$\ell$-weight monomials in the monoid generated by
the variables $Y_{i,p}$ with $(i,p)\in\mathcal X$.
The ambient category $\mathcal C_{\mathbb Z}$ is defined in
the same way using all $p\in\xi_i+2\mathbb Z$.

For a $k$-element subset $J$ of $[n]$, let $P_J$
denote the corresponding Pl\"ucker coordinate.
Hernandez and Leclerc~\cite[Section~13]{HL10a}
established an isomorphism of $\mathbb C$-algebras
\begin{equation}\label{eq:classical-isomorphism}
 \Phi_{\mathrm{cl}}:
 \mathbb C[\Gr(k,n,\sim)]
 \xrightarrow{\sim}
 \mathbb C\otimes_{\mathbb Z}K_0(\mathcal C_\ell),
\end{equation}
where $K_0(\mathcal C_\ell)$ is the Grothendieck ring
of $\mathcal C_\ell$ (with the relations $[V]=[U]+[W]$
for short exact sequences $0\to U\to V\to W\to0$ and multiplication
$[U][V]=[U\otimes V]$), and
\[
 \mathbb C[\Gr(k,n,\sim)]
 :=
 \mathbb C[\Gr(k,n)]/
 \left\langle
 P_{[a,a+k-1]}-1:1\leq a\leq n-k+1
 \right\rangle.
\]

Building on this isomorphism, Chang, Duan, Fraser,
and the author~\cite[Theorem~2.1, Lemma~3.9 and
Theorem~3.17]{CDFL20} described the simple modules
in $\mathcal C_\ell$ in terms of semistandard Young
tableaux.
Let $\SSYT(k,[n])$ denote the set of
rectangular semistandard Young tableaux with $k$ rows
and entries in $[n]$, including the empty tableau.
The rows are weakly increasing from left to right and the columns
are strictly increasing from top to bottom.
For two such tableaux $T$ and $S$, their union
$T\cup S$ is obtained by taking the multiset union
of the entries in each row and arranging the entries
in weakly increasing order.
A column of the form
$C_a=[a,a+k-1]$, $1\leq a\leq n-k+1$,
is called a trivial column.
Let $\sim$ be the equivalence relation on
$\SSYT(k,[n])$ generated by
$T\sim T\cup C_a$ for all trivial columns $C_a$.
The quotient $\SSYT(k,[n],\sim)
 :=
 \SSYT(k,[n])/\sim$
parametrizes the isomorphism classes of simple modules
in $\mathcal C_\ell$ \cite{CDFL20}. More explicitly, for a column
$J=\{a_1<\cdots<a_k\}$, set
\[
 M_J=
 \prod_{r=1}^{k-1}
 \prod_{a=a_r+1}^{a_{r+1}-1}
 Y_{k-r,k+r-2a}.
\]
For a tableau $T$ with columns $T_1,\ldots,T_d$, put $M_T=\prod_{j=1}^d M_{T_j}$.
Since $M_{C_a}=1$ for every trivial column $C_a$,
the monomial $M_T$ depends only on the equivalence
class of $T$.
The correspondence $T\mapsto M_T$ induces a bijection
between $\SSYT(k,[n],\sim)$ and the monoid
of dominant monomials generated by the variables
$Y_{i,p}$ with $(i,p)\in\mathcal X$.
In~\cite[Proposition~5.17 and the subsequent discussion]{CDFL20},
the authors gave evidence for the identification of their
Grassmannian basis with the dual canonical basis of the Grassmannian coordinate ring.
The identification of the class of $L(M_T)$ with the image
of the latter basis element indexed by $T$ will also follow
from Theorem~\ref{thm:intro-A} by specialization.
In particular, the class of the fundamental module
$L(Y_{i,i-2b})$ corresponds to the image of the
Pl\"ucker coordinate indexed by the one-column
tableau
\begin{equation}\label{eq:fundamental-columns}
 I_{b,i}=[b,b+k]\setminus\{b+k-i\},
 \qquad 1\leq b\leq m,\quad 1\leq i<k.
\end{equation}
Feigin and the author~\cite[Theorem~B]{FeiginLi2026}
proved that reduced standard monomials and reduced dual
canonical basis elements give bases throughout the family
of classical quotients obtained by specializing the
non-wrapping consecutive Pl\"ucker coordinates to arbitrary
scalars.

\subsection{Quantum cluster algebras and monoidal categorification}

We now turn to the quantum setting.
Let $v$ be an indeterminate and put
$\bbk=\mathbb Q(v^{1/2})$.
Let $\mathcal O_v(M_{k,n})$ be the algebra of quantum
$k\times n$ matrices over $\bbk$, with generators
$X_{r,a}$ for $1\leq r\leq k$ and $1\leq a\leq n$ and the
relations in Section~\ref{sec:twist}.
The quantum Grassmannian
$A_v=\mathcal O_v(\Gr(k,n))$
is the subalgebra of $\mathcal O_v(M_{k,n})$ generated
by all maximal quantum minors; see
\cite{Lauve2006,TaftTowber1991}.
For a $k$-element subset
$J=\{a_1<\cdots<a_k\}$ of $[n]$, the corresponding
quantum Pl\"ucker coordinate is
\begin{equation}\label{eq:quantum-plucker-coordinate}
 P_J
 =
 \sum_{\sigma\in\mathfrak S_k}
 (-v)^{\ell(\sigma)}
 X_{1,a_{\sigma(1)}}\cdots X_{k,a_{\sigma(k)}},
\end{equation}
where $\mathfrak S_k$ is the symmetric group and
$\ell(\sigma)=|\{(a,b):a<b,\ \sigma(a)>\sigma(b)\}|$.
We use the same notation $P_J$ for classical and
quantum Pl\"ucker coordinates whenever the meaning
is clear from the context.

Quantum cluster algebras were introduced by
Berenstein and Zelevinsky~\cite{BZ05}.
Geiss, Leclerc, and Schr\"oer~\cite{GLS13} established
quantum cluster algebra structures on quantized
coordinate rings of unipotent subgroups of symmetric
Kac--Moody groups associated with Weyl group elements.
Grabowski and Launois~\cite{GL11,GL14} established
quantum cluster algebra structures on quantum
Grassmannians $A_v=\mathcal O_v(\Gr(k,n))$ quantizing the classical cluster
algebra structure of Scott~\cite{Scott2006}.

On the representation-theoretic side,
Frenkel and Reshetikhin~\cite{FR99} introduced
$q$-characters $\chi_q$ of finite-dimensional representations
of quantum affine algebras.
Nakajima~\cite{Nakajima04} introduced $t$-analogues
of $q$-characters using the geometry of quiver varieties.
Hernandez~\cite{Hernandez04} developed an algebraic
approach to $(q,t)$-characters and quantum Grothendieck
rings.
Hernandez and Leclerc~\cite[Theorems~1.1 and~7.3]{HL2015}
gave a presentation of the quantum Grothendieck ring
$K_t(\mathcal C_{\mathbb Z})$ over $\mathbb C(t^{1/2})$
in simply-laced types and identified a specialization
of this ring with a twisted derived Hall algebra
of a Dynkin quiver.
For the subcategory $\mathcal C_Q$ associated with an
orientation $Q$ of the Dynkin diagram, they established
an isomorphism between its quantum Grothendieck ring and
the quantum coordinate ring of the maximal unipotent
subgroup, identifying the $(q,t)$-characters of simple
modules with the dual canonical basis up to scalar
normalization; see~\cite[Theorems~1.2 and~6.1]{HL2015}.

Hernandez and Leclerc~\cite{HL10a} introduced monoidal
categorification of cluster algebras and conjectured
in~\cite[Conjecture~13.2]{HL10a} that the cluster monomials
in the Grothendieck ring of $\mathcal C_\ell$ are precisely
the classes of real simple modules, where a simple
module $S$ is called real if $S\otimes S$ is simple.
Kang, Kashiwara, Kim, and Oh~\cite{KKKO2018} established
monoidal categorifications of quantum unipotent coordinate rings
in symmetric types.
Using common triangular bases,
Qin~\cite[Theorem~9.4.2]{Qin2017} established the
monoidal categorification of Hernandez--Leclerc
cluster algebras associated with $\mathcal C_\ell$
in simply-laced types, showing that every quantum
cluster monomial in $K_t(\mathcal C_\ell)$ is the
$(q,t)$-character of a simple module.
Kashiwara, Kim, Oh, and Park
\cite{KKOP2020,KKOP2021,KKOP2022,KKOP2024}
established monoidal categorifications for a broad
class of quantum affine representation categories,
including non-simply-laced and twisted types.
In particular, cluster monomials in the corresponding
Grothendieck rings are classes of real simple modules.
Oh and Park~\cite{OP2025} developed the PBW theory for
bosonic extensions of quantum groups of arbitrary finite type
and established Levendorskii--Soibelman relations for PBW root vectors.
Kashiwara, Kim, Oh,
and Park established the global basis theory and braid
symmetries for bosonic extensions of quantum unipotent
coordinate rings~\cite{KKOPGlobal2025,KKOPBraid2024}.
In simply-laced finite type, after extension of scalars,
the bosonic extension is isomorphic to the quantum
Grothendieck ring of $\mathcal C_{\mathbb Z}$, with the
normalized global basis corresponding to its canonical
basis~\cite[Theorem~8.4]{KKOPGlobal2025}.
Qin~\cite[Theorems~1.6 and~6.6]{Qin2024Infinite}
established quantum cluster algebra structures on subalgebras
of bosonic extensions associated with positive braid monoid
elements in finite type, together with monoidal categorifications
in type $ADE$.
In simply-laced finite types, Bi~\cite[Theorems~3.17 and~4.8]{Bi2025}
and Kashiwara, Kim, Oh, and Park~\cite[Theorem~9.7]{KKOPIII2025}
established such structures with quantum cluster monomials
belonging to the global basis, with the appropriate normalization.
The latter work treats the standard integral form and gives
monoidal categorifications using representations of quantum
affine algebras, at the classical level for complete duality data
and at the quantum level in the untwisted affine setting
\cite[Theorems~9.4 and~9.12]{KKOPIII2025}.
Bi~\cite[Theorems~1.2 and~1.3]{Bi2026Bosonic}
constructed quantum cluster structures on the standard integral
forms for arbitrary symmetrizable generalized Cartan matrices,
and monoidal categorifications for symmetric generalized Cartan
matrices using localized categories of graded modules over
quiver Hecke algebras.

Bittmann~\cite[Theorem~5.2.4]{Bittmann21} established
a quantum cluster algebra structure on the quantum
Grothendieck ring of the subcategory
$\mathcal C_{\mathbb Z}^-$ of $\mathcal C_{\mathbb Z}$,
whose simple objects have highest $\ell$-weight
monomials involving only the variables $Y_{i,p}$
with $p\leq 0$, in simply-laced types.
Fujita, Hernandez, Oh, and Oya
\cite[Theorem~8.6 and Corollary~8.7]{FHOOPos2022}
extended the quantum unipotent coordinate ring
realization associated with $\mathcal C_Q$ to
subcategories associated with $Q$-data
introduced in~\cite[Definition~2.5]{FO2021}, and established
canonical basis preserving isomorphisms between
quantum Grothendieck rings associated with different
types having the same simply-laced unfolding
\cite[Theorem~10.6]{FHOOPos2022}.
They further established quantum cluster algebra
structures on $K_t(\mathcal C^{\leq\xi})$, where
$\mathcal C^{\leq\xi}$ is the subcategory with highest monomials
supported on spectral indices at or below the height function
of a $Q$-datum (see~\cite[Section~4.3]{FHOOCluster2023}),
and showed
that quantum cluster monomials belong to the canonical
basis of these quantum Grothendieck rings
\cite[Theorems~5.16 and~6.6]{FHOOCluster2023}.

Qin~\cite{Qin2024Lie} established the existence of
common triangular bases for a broad class of quantum
cluster algebras arising from Lie theory, including
those associated with double Bott--Samelson cells,
braid varieties, and Grassmannians.
For symmetric Cartan matrices, he also established
quasi-categorifications, yielding positivity of the
structure constants of these bases.
In particular, he related the distinguished seed
of a Grassmannian to a seed associated with an
adapted word in type $A$
\cite[Sections~7.5--7.6]{Qin2024Lie}.
Cluster algebra structures on braid varieties were
constructed in~\cite{CGGLSS2025}.
Huh, Jung, Kim, and Park
\cite[Theorem~7.7 and Corollary~9.3]{HJKP2026}
established an explicit isomorphism between the localized
classical specializations of bosonic extensions associated
with positive braid monoid elements and the coordinate
rings of the corresponding braid varieties, relating
the cluster structures arising from admissible chains
of $i$-boxes and Demazure weaves.
They also obtained monoidal categorifications of the
associated partially compactified cluster algebras.
Bi~\cite[Theorems~8.15 and~8.16]{Bi2026Braid}
established basis-preserving isomorphisms between quantum
cluster algebras associated with braid varieties,
localized quantum Grothendieck rings, and localized
subalgebras of bosonic extensions.

Let $K_t(\mathcal C_\ell)$ be the quantum Grothendieck
ring of $\mathcal C_\ell$, realized inside the quantum torus
(see Section~\ref{subsec:general-spectral-words}).
The ordinary $q$-character of a finite-dimensional module
encodes its decomposition into common generalized
eigenspaces for the Drinfeld Cartan currents, with the
dimensions of these eigenspaces as coefficients of the
corresponding Laurent monomials~\cite{FR99}.
For a dominant monomial $M$ corresponding to a
simple module in $\mathcal C_\ell$, we denote by
$\chi_{q,t}(L(M))$ the $(q,t)$-character of $L(M)$,
viewed as an element of $K_t(\mathcal C_\ell)$.
By~\cite[Section~3.3 and Proposition~5.7]{FHOOPos2022},
$K_t(\mathcal C_\ell)$ is generated over $\mathbb Z[t^{\pm1/2}]$
by the fundamental $(q,t)$-characters $\chi_{q,t}(L(Y_{i,p}))$,
$(i,p)\in\mathcal X$.
The elements $\chi_{q,t}(L(M))$, where $M$ ranges over the
dominant monomials in the variables $Y_{i,p}$ with
$(i,p)\in\mathcal X$, form the canonical basis of
$K_t(\mathcal C_\ell)$ over $\mathbb Z[t^{\pm1/2}]$.
Following~\cite[Section~3.2]{HL2015}, we use commutative
Laurent monomials in the variables $Y_{i,p}$ as the
basis of the quantum torus $\mathcal Y_t$.
We write $*$ for its quantum multiplication and
use juxtaposition of Laurent monomials for their
commutative product; see Section~\ref{subsec:general-spectral-words}.
The elements $\chi_{q,t}(L(M))$ are invariant under the bar of $\mathcal Y_t$
and have leading monomial $M$ with coefficient one.
Their normalization is characterized by the bar and the standard
basis in Section~\ref{subsec:general-spectral-words}.

Throughout this paper, $q$ denotes the parameter of
the quantum affine algebra
$U_q(\widehat{\mathfrak{sl}_k})$, $v$ denotes the
deformation parameter of the quantum Grassmannian,
and $t$ denotes the deformation parameter of the
quantum Grothendieck ring.
Our isomorphism identifies $v$ with $t$.
The parameter $q$ is independent of these deformation
parameters. Under the identification $v=t$, we use $t$
also for scalar factors in the central quotients.

\subsection{The quantum correspondence}

Our starting point is the classical Grassmannian quotient
in~\eqref{eq:classical-isomorphism} and the correspondence
between semistandard Young tableaux and simple modules
established in~\cite{CDFL20}.
We study the quantization of this quotient together
with its dual canonical basis.
More precisely, we construct a cocycle twist of the
homogeneous quantum Grassmannian
(Proposition~\ref{prop:center}) and establish an
isomorphism between a central quotient of the twisted
algebra and the quantum Grothendieck ring of
$\mathcal C_\ell$
(Theorem~\ref{thm:general-based-isomorphism}).
We determine the images of all quantum Pl\"ucker
coordinates and, more generally, of all dual canonical basis
elements of the quantum Grassmannian indexed by semistandard
Young tableaux; see
\eqref{eq:scaled-plucker-images}
and~\eqref{eq:scaled-tableau-images}.

A difficulty is that the non-wrapping consecutive
quantum Pl\"ucker coordinates are not central in
general. Consequently, the desired quotient cannot
be obtained simply by setting these coordinates equal
to $1$ in the quantum Grassmannian.
For example, in the quantum Grassmannian associated
with $\Gr(3,6)$, we have $P_{123}P_{124}=vP_{124}P_{123}$.
Setting $P_{123}=1$ would imply $(1-v)P_{124}=0$,
and hence $P_{124}=0$ over $\bbk$.
This is incompatible with the desired correspondence,
since $P_{124}$ corresponds in the classical case
to the class of the fundamental module
$L(Y_{1,-1})$ in $\mathcal C_2$.

To overcome this difficulty, we construct a cocycle twist
of the quantum Grassmannian under which all non-wrapping
consecutive quantum Pl\"ucker coordinates become central,
allowing us to specialize them to nonzero scalars.
Our twist differs from the one used by Launois and
Lenagan~\cite[Theorem~5.9]{LaunoisLenagan2011}
to obtain a quantum analogue of the cyclic permutation
of Pl\"ucker coordinates.
Let $e_1,\ldots,e_n$ be the standard basis of
$\mathbb Z^n$. The algebra $\mathcal O_v(M_{k,n})$
is graded by $\mathbb Z^n$ by setting $\deg X_{r,a}=e_a$.
The quantum Grassmannian $A_v$ inherits the content grading. In particular,
each quantum Pl\"ucker coordinate $P_J$ is homogeneous
of degree $\deg P_J=\sum_{a\in J}e_a$.
Define a skew-symmetric bilinear form
$\Omega:\mathbb Z^n\times\mathbb Z^n\to\mathbb Z$
by
\[
 \Omega(e_a,e_b)=
 \begin{cases}
 -1,&a<b,\ k\mid(b-a),\\
 1,&a>b,\ k\mid(a-b),\\
 0,&\text{otherwise}.
 \end{cases}
\]
Let $A_v^\Omega$ be the $\bbk$-vector space $A_v$
equipped with the multiplication
$f\star g = v^{\Omega(\deg f,\deg g)/2}fg$
for homogeneous elements $f,g\in A_v$, extended
bilinearly to arbitrary elements.
The associativity of this multiplication follows
from the bilinearity of $\Omega$.

Put $s=n-k+1=m+1$. For $1\leq a\leq s$, write
$F_a=P_{C_a}=P_{[a,a+k-1]}$.
By Proposition~\ref{prop:center}, all these elements
are central in $A_v^\Omega$. 
The wrapping consecutive columns are the length-$k$ cyclic
intervals starting at $s+1,\ldots,n$ (indices taken modulo $n$).
Their Pl\"ucker coordinates are retained as frozen variables
in the central quotient $B_\chi$ defined below.
We localize at their images in
Section~\ref{subsec:cyclic-quasi-automorphism}
to construct quantum cluster quasi-automorphisms.

Let $b_T$ denote the dual canonical basis element
of $A_v$ indexed by a tableau $T$, normalized with
respect to the bar involution that sends $v^{1/2}$
to $v^{-1/2}$, fixes the quantum matrix generators,
and reverses the order of multiplication. Its standard monomial
characterization is given in Section~\ref{subsec:homogeneous-standard-monomials}.
In particular, for a one-column tableau $J$,
we have $b_{J}=P_J$.
The elements $b_T$ are homogeneous for the content
grading.
We also regard them as elements of the underlying
vector space of $A_v^\Omega$.

Following~\cite[Section~3.2]{CDFL20}, a tableau $T$
is called reduced if it cannot be written as
$T=S\cup C_a$ for a tableau $S$ and a trivial column
$C_a$.
Every equivalence class in
$\SSYT(k,[n],\sim)$ has a unique reduced
representative, denoted by $T_{\mathrm{red}}$.

We extend scalars in the quantum Grothendieck ring
by setting
$K= \bbk\otimes_{\mathbb Z[t^{\pm1/2}]} K_t(\mathcal C_\ell)$,
where the homomorphism
$\mathbb Z[t^{\pm1/2}]\to\bbk$ sends
$t^{1/2}$ to $v^{1/2}$. To specify the scalar factors in our isomorphism,
we fix the additive map
$\chi:\mathbb Z^n\to\mathbb Z$ given by
$\chi(e_a)=\delta_{a,n}-\delta_{a,1}$.
For a column $J$ and a tableau $T$, write
$\chi(J)=\chi(\deg P_J)$, $\chi(T)=\chi(\deg b_T)$.
Thus $\chi(J)$ is the difference between the
number of occurrences of $n$ and that of $1$
in the column $J$, and similarly for $\chi(T)$.
The prescribed scalar
values of the consecutive minors are
$F_1\mapsto t^{-1/2}$, $F_a\mapsto 1 \ (2\leq a\leq m)$, $F_{m+1}\mapsto t^{1/2}$.
This choice of scalar factors fixes the
normalization of our isomorphism.

\begin{mainresult}[Theorem~\ref{thm:general-based-isomorphism}]
\label{thm:intro-A}
Let $k\geq 2$ and $\ell\geq 0$.
The elements $F_a$, $1\leq a\leq\ell+2$,
are central in $A_v^\Omega$.
There is an isomorphism of $\bbk$-algebras
\[
 \Phi_\chi:B_\chi\xrightarrow{\sim}K,
\]
where
\[
 B_\chi
 =
 A_v^\Omega/
 \left\langle
 F_a-t^{\chi(C_a)/2}:
 1\leq a\leq\ell+2
 \right\rangle.
\]
Let $\pi_\chi:A_v^\Omega\to B_\chi$ be the
quotient map.

For every $k$-element subset
$J=\{a_1<\cdots<a_k\}$ of $[n]$, we have
\begin{equation}\label{eq:scaled-plucker-images}
 \Phi_\chi(\pi_\chi(P_J))
 =
 t^{\chi(J)/2}\chi_{q,t}(L(M_J)),
 \qquad
 M_J=
 \prod_{r=1}^{k-1}
 \prod_{a=a_r+1}^{a_{r+1}-1}
 Y_{k-r,k+r-2a}.
\end{equation}
In particular,
\[
 \Phi_\chi(\pi_\chi(P_{I_{b,i}}))
 =
 t^{\chi(I_{b,i})/2}\chi_{q,t}(L(Y_{i,i-2b})),
 \qquad
 1\leq b\leq m,\quad 1\leq i<k.
\]
The elements $\pi_\chi(P_{I_{b,i}})$ generate
$B_\chi$ as a $\bbk$-algebra, and their images
determine the isomorphism uniquely.

The elements $\pi_\chi(b_T)$, where $T$ ranges
over all reduced tableaux, form a $\bbk$-basis
of $B_\chi$.
Moreover, for every tableau $T$, we have
\begin{equation}\label{eq:scaled-tableau-images}
 \Phi_\chi(\pi_\chi(b_T))
 =
 t^{\chi(T)/2}\chi_{q,t}(L(M_T)).
\end{equation}

The map on homogeneous elements of $A_v^\Omega$
given by
$f\mapsto
 t^{-\chi(\deg f)}\overline f$
induces a bar involution on $B_\chi$.
The isomorphism $\Phi_\chi$ intertwines this
involution with the bar involution on $K$.
In particular, it identifies the dual canonical basis
$\{t^{-\chi(T)/2}\pi_\chi(b_T): T\text{ is reduced}\}$
of $B_\chi$ with the canonical basis
$\{\chi_{q,t}(L(M)):M\in\mathcal P_\ell^+\}$ of $K$.
Both bases are bar invariant.
\end{mainresult}

At $v^{1/2}=t^{1/2}=1$, the cocycle twist and the scalar
factors disappear. The isomorphism on the integral forms
spanned by the respective normalized bases specializes to
the classical isomorphism in~\cite{HL10a,CDFL20}, and identifies
the images of the dual canonical basis of the classical Grassmannian
coordinate ring with the classes of simple modules.

The proof proceeds by a direct comparison of the
algebraic structures and bases on the two sides.
We first establish a multiplication formula for
the dual canonical basis elements of $A_v$ and
the consecutive quantum Pl\"ucker coordinates.
After the cocycle twist and central specialization,
this formula gives a basis of the quotient indexed
by reduced tableaux.
We then show that suitably normalized ordered
products of the fundamental Pl\"ucker coordinates
$P_{I_{b,i}}$ form another basis of the quotient
and establish the triangular transition between
these two bases (Proposition~\ref{thm:general-source-pbw}).
Using the defining relations of quantum minors
and the ordered-product basis, we establish a
presentation of the quotient by the fundamental
Pl\"ucker coordinates and an explicit family of
quadratic relations (Proposition~\ref{thm:general-defining}).

Hernandez and Leclerc~\cite[Theorem~7.3]{HL2015}
gave a presentation of the quantum Grothendieck ring of
$\mathcal C_{\mathbb Z}$ over $\mathbb C(t^{1/2})$ in
simply-laced types, using fundamental $(q,t)$-characters
indexed by simple roots and integers.
Fujita, Hernandez, Oh, and Oya~\cite[Theorem~10.4]{FHOOPos2022}
extended this presentation to non-simply-laced types
using $Q$-data.
A corresponding presentation for quantum virtual
Grothendieck rings was established by Jang, Lee,
and Oh~\cite[Theorem~7.2]{JLO2023Braid}.
For our purposes, we derive a presentation of
$K_t(\mathcal C_\ell)$ by all fundamental $(q,t)$-characters
in $\mathcal C_\ell$, with explicit quadratic relations,
from the segment relations of Kang, Kashiwara and Kim
and the normalization used in the realization of the
quantum Grothendieck ring in~\cite{Fujita2024}
(Proposition~\ref{thm:general-target-relations}).
The two presentations of the central quotient and the
quantum Grothendieck ring coincide, yielding the desired
algebra isomorphism.
The compatibility with the bar involutions and the
triangular expansions of the two bases then determine
the images of all dual canonical basis elements of $A_v$.

As a consequence of the isomorphism, the
appropriately scaled $(q,t)$-characters of simple modules
corresponding to quantum Pl\"ucker coordinates
satisfy all the twisted quantum Young symmetry
relations~\eqref{eq:young-twisted};
see~\eqref{eq:young-character}.

\subsection{Kazhdan--Lusztig formulas for
\texorpdfstring{$(q,t)$}{(q,t)}-characters}

The computation of $(q,t)$-char\-acters has been studied
by geometric methods using quiver varieties
\cite[Sections~5 and~8]{Nakajima04}
and by algebraic methods
\cite[Sections~5.7 and~6.4]{Hernandez04}.
Hernandez and Leclerc~\cite[Theorem~3.1]{HL2016} gave a
cluster algebra algorithm for computing the $q$-characters
of Kirillov--Reshetikhin modules.
Bittmann~\cite[Proposition~6.3.1]{Bittmann21} extended
this approach to the computation of their $(q,t)$-characters.
See also Hernandez and Oya~\cite{HO2019} for an approach
to the computation of $(q,t)$-characters using isomorphisms
between quantum Grothendieck rings and the Kazhdan--Lusztig
algorithm.

Du~\cite{Du1992,Du1995} obtained Kazhdan--Lusztig formulas
for the transition between semistandard and canonical bases
of irreducible polynomial representations of quantum
$\mathrm{GL}_n$.
Brundan~\cite{Brundan2005} recovered these formulas using
tensor products and Schur--Weyl duality, and obtained
corresponding formulas for the dual canonical basis.
Combining the latter formulas with
Theorem~\ref{thm:general-based-isomorphism},
we obtain an explicit finite formula for the
$(q,t)$-characters of arbitrary simple modules in
$\mathcal C_\ell$ in terms of fundamental
$(q,t)$-characters and Kazhdan--Lusztig polynomials
of symmetric groups.

More precisely, we have the following formula.
\begin{mainresult}[Theorem~\ref{thm:kl-character-formula}]
\label{thm:intro-B}
For every $M\in\mathcal P_\ell^+$, we have
\begin{equation}\label{eq:intro-kl-character-formula}
 \chi_{q,t}(L(M))
 =
 \sum_{s\geq 0}(-1)^s
 \sum_{\substack{M=M_0>M_1>\cdots>M_s\\
                 M_a\in\mathcal P_\ell^+}}
 \left(
 \prod_{a=0}^{s-1}c(M_a,M_{a+1};t)
 \right)
 \mathcal E(M_s),
\end{equation}
where the inner sum runs over strictly decreasing
chains of dominant monomials with respect to the
partial order $M>M'$ if $M'/M$ is a nontrivial product of inverse simple root monomials
$A_{i,p}^{-1}$, with $A_{i,p}$ in~\eqref{eq:simple-root-monomials}.
The sum is finite, and the term corresponding to
$s=0$ is $\mathcal E(M)$.

For $M'=y_1\cdots y_w$, where $y_a=Y_{i_a,p_a}$ and
$p_1\geq\cdots\geq p_w$, the normalized ordered product is
\[
 \mathcal E(M')=t^{-\frac12\sum_{a<b}J(y_a,y_b)}
 \chi_{q,t}(L(y_1))*\cdots*\chi_{q,t}(L(y_w)).
\]
Here $J$ is the skew-symmetric bilinear form defined in
Section~\ref{subsec:general-spectral-words}, and $\mathcal E(1)=1$.
The coefficients $c(M,M';t)$ are given explicitly in terms of
Kazhdan--Lusztig polynomials of symmetric groups
by~\eqref{eq:kl-tableau-coefficient}
and~\eqref{eq:kl-monomial-coefficient}.
\end{mainresult}

The formulas in~\cite[Theorems~5.4 and~5.8]{CDFL20}
use values of Kazhdan--Lusztig polynomials at $1$ to express
$q$-characters and Grassmannian basis elements.
In~\eqref{eq:intro-kl-character-formula}, the coefficients
retain the deformation parameter and are expressed as finite
sums of products of explicit expressions involving
Kazhdan--Lusztig polynomials. The resulting formula expresses the $(q,t)$-character of
a simple module in terms of fundamental $(q,t)$-characters
and these coefficients.

\subsection{Quantum cyclic shifts and promotion}

We finally study a quantum analogue of the cyclic shift on
the homogeneous Grassmannian and the images of the shifted
dual canonical basis elements in the central quotient.
After inverting the remaining frozen variables, we construct
a quantum cluster quasi-automorphism of the localized central quotient
and of the localized quantum Grothendieck ring, in the sense
of~\cite{CHLQuasi2024}.
Write $\sigma$ for the cyclic permutation sending $n$ to $1$
and $a$ to $a+1$ for $a<n$. For a subset $I\subseteq[n]$,
set $\sigma I=\{\sigma(i):i\in I\}$. On content vectors, the action is
$(\sigma\alpha)_{\sigma(a)}=\alpha_a$. The indices in
$P_{\sigma^j I}$ are always taken in the natural increasing order
on $[n]$. Write $\operatorname{pr}$ for Sch\"utzenberger promotion
on rectangular semistandard tableaux, with the convention that
it removes the entries $n$ and transforms the content by $\sigma$.
For $\alpha\in\mathcal L=\{\alpha\in\mathbb Z^n:
|\alpha|\in k\mathbb Z\}$, put $d_\alpha=|\alpha|/k$,
where $|\alpha|=\sum_a\alpha_a$.
For the normalized specialization, write
$\mathcal D(f)=v^{-\chi(\deg f)/2}\pi_\chi(f)$ on homogeneous
$f$, extended linearly. For $1\leq a\leq n$, use $C_a$ for the
cyclic consecutive column starting at $a$ and put
$z_a=\mathcal D(P_{C_a})$; thus $z_1=\cdots=z_s=1$.
For bar-invariant elements with $uw=v^awu$, write $[uw]=v^{-a/2}uw$.
For several pairwise quasi-commuting factors, the brackets denote the
normalized product described in Section~\ref{subsec:cyclic-quasi-automorphism}.
The same notation is used with $v=t$ on the Grothendieck ring side.

\begin{mainresult}[Theorems~\ref{thm:cyclic-canonical-basis} and~\ref{thm:cyclic-quasi-automorphism}]
\label{thm:intro-C}
\leavevmode

\noindent\textup{(i)} There is a unique $\bbk$-linear quantum cyclic shift
$\widehat\rho_v$ on $A_v$ which maps $P_I$ to $P_{\sigma I}$ and satisfies
$\widehat\rho_v(fg)=v^{d_\alpha\beta_n-\alpha_n d_\beta}
\widehat\rho_v(f)\widehat\rho_v(g)$ for homogeneous $f,g$
of contents $\alpha,\beta$. It satisfies $\widehat\rho_v^n=1$
and $\widehat\rho_v(b_T)=b_{\operatorname{pr}(T)}$.
Over $\mathbb Q(v^{1/(2k)})$, the cyclic shift $\rho_v^{\rm LL}$
of Launois and Lenagan satisfies
\[
 \rho_v^{\rm LL}(b_T)=
 v^{-r(d+1)+r(r-1)/k}b_{\operatorname{pr}(T)},
\]
where $T$ has $d$ columns and $r$ entries equal to $n$.
For every tableau $T$ and $0\leq j<n$, specialization to the
central quotient $B_\chi$ gives
\[
 \Phi_\chi\bigl(v^{-\chi(\operatorname{pr}^j(T))/2}
 \pi_\chi(\widehat\rho_v^j(b_T))\bigr)
 =\chi_{q,t}(L(M_{\operatorname{pr}^j(T)})).
\]
At $v^{1/2}=1$, $\widehat\rho_v$ specializes to the classical signed cyclic
shift $\rho$ on $\mathbb C[\Gr(k,n)]$.
Neither map descends to its central quotient.

\noindent\textup{(ii)} The localization
$B_\chi^{\mathrm{loc}}=B_\chi[z_{s+1}^{-1},\ldots,z_n^{-1}]$
has the quantum cluster structure obtained by twisting and
specializing the Grassmannian quantum seeds. There is a quantum
cluster quasi-automorphism $\widetilde\rho_v$ of this algebra
satisfying
\[
 \widetilde\rho_v(\mathcal D(b_T))
 =[z_{s+1}^{-r_n(T)}\mathcal D(b_{\operatorname{pr}(T)})],
\]
where $r_n(T)$ is the number of entries equal to $n$ in $T$.
In particular, for a single column $I$ the formula is
$\widetilde\rho_v(\mathcal D(P_I))=
 [z_{s+1}^{-\mathbf1_{n\in I}}\mathcal D(P_{\sigma I})]$.
Let $Z_i=\chi_{q,t}(L(\prod_{b=1}^{m}Y_{i,i-2b}))$ for $1\leq i<k$.
Under Theorem~\ref{thm:intro-A} this gives a quantum cluster
quasi-automorphism $\mathcal R_t$ of
$K^{\mathrm{loc}}=K[Z_1^{-1},\ldots,Z_{k-1}^{-1}]$, with
\[
 \mathcal R_t(\chi_{q,t}(L(M_T)))=
 [Z_{k-1}^{-r_n(T)}\chi_{q,t}(L(M_{\operatorname{pr}(T)}))].
\]
The maps $\widetilde\rho_v$ and $\mathcal R_t$ permute the localized
dual canonical and canonical bases, respectively, obtained by
normalized multiplication of basis elements by frozen Laurent monomials.
Put $w_b=z_{s+b}$ for $1\leq b<k$, and let $c_b(T)$ count the
entries of $T$ congruent to $b$ modulo $k$, for $1\leq b\leq k$.
Then
\[
 \widetilde\rho_v^n(\mathcal D(b_T))=
 \left[\left(\prod_{b=1}^{k-1}w_b^{c_{b+1}(T)-c_b(T)}\right)
 \mathcal D(b_T)\right],
\]
where products involving several frozen factors are normalized as
quantum Laurent monomials. The order of $\widetilde\rho_v$ (and of
$\mathcal R_t$) as an algebra automorphism is
\[
 \begin{cases}
 k,&n=k+1,\\
 kn,&n\geq k+2,\ \gcd(k,n)=1,\\
 \infty,&\gcd(k,n)>1.
 \end{cases}
\]
The corresponding formulas at $v^{1/2}=1$ give classical cluster
quasi-automorphisms of the localized quotients, with the same orders.
\end{mainresult}

The formula for $\widehat\rho_v(b_T)$ is a quantum version of the
classical cyclic symmetry of the dual canonical basis
\cite[Proposition~5.5]{Rhoades2010}, \cite[Theorem~2]{Lam2019}.
In Section~\ref{sec:cyclic-shifts},
we compare the normalization with that of Launois and Lenagan
and distinguish the uncorrected specialization from the quantum cluster
quasi-automorphism of the localized quotient.

\subsection{Organization of the paper} 
In Section~\ref{sec:twist}, we construct the cocycle
twist and describe the central quotients.
Section~\ref{sec:canonical-specialization} establishes the standard
and dual canonical bases of the central quotient and their triangular
relation. Section~\ref{subsec:general-spectral-words} recalls the
standard and canonical bases of the quantum Grothendieck ring and
determines the dominant monomials in products of fundamental characters.
In Section~\ref{sec:defining-relations}, we compare normalization
factors, establish the two presentations, and compare our presentation
of the quantum Grothendieck ring with existing presentations.
We then prove Theorem~\ref{thm:intro-A} and give an explicit example.
In Section~\ref{sec:kl-characters}, we derive the
Kazhdan--Lusztig formula for $(q,t)$-characters of simple modules
and illustrate it with an example.
Section~\ref{sec:cyclic-shifts} proves Theorem~\ref{thm:intro-C}.
It studies the specialization of shifted basis elements, constructs
the quantum cluster quasi-automorphisms of the localized quotient
and the localized quantum Grothendieck ring, and determines their orders.

\subsection*{Acknowledgments}
The author would like to thank Deniz Kus for helpful discussions.
The author was supported by the Beatriz Galindo Senior grant BG24/00114 from the Spanish Ministry of Science, Innovation and Universities. 

\subsection*{AI Use Disclosure}
ChatGPT 6 was used solely for language editing, readability improvement, and writing computer programs. The author takes full responsibility for the content of this paper.

\section{A cocycle twist and central quotients of the quantum Grassmannian}
\label{sec:twist}

In this section, we construct a cocycle twist of the
quantum Grassmannian under which the non-wrapping
consecutive quantum Pl\"ucker coordinates become central.
We describe the defining relations of the twisted algebra
and of the quotients obtained by specializing these
central elements to scalars. We also study the
compatibility of the bar involution with the cocycle
twist and central specialization.

Let $s=n-k+1=\ell+2$ and $\bbk=\mathbb Q(v^{1/2})$.
As in the introduction, we write
$A_v=\mathcal O_v(\Gr(k,n))$, $C_a=[a,a+k-1]$, $F_a=P_{C_a},  1\leq a\leq s$.

\subsection{Quantum matrices and Young symmetry relations}

We realize $A_v$ as the subalgebra of the quantum matrix
algebra $\mathcal O_v(M_{k,n})$ generated by its maximal
quantum minors.
We use the Faddeev--Reshetikhin--Takhtajan (FRT)
convention~\cite[Equations~(2) and~(6)]{FRT90};
see also~\cite[Section~2]{LenaganRussell2008}
for rectangular quantum matrices.
The matrix generators $X_{r,a}$ satisfy
\[
\begin{aligned}
X_{r,a}X_{r,b}&=vX_{r,b}X_{r,a},
&\qquad X_{r,a}X_{r',a}&=vX_{r',a}X_{r,a},\\
X_{r,b}X_{r',a}&=X_{r',a}X_{r,b},
&\qquad X_{r,a}X_{r',b}-X_{r',b}X_{r,a}
&=(v-v^{-1})X_{r,b}X_{r',a}.
\end{aligned}
\]
where the first relation holds for $1\leq a<b\leq n$,
the second for $1\leq r<r'\leq k$, and the last two
for $1\leq r<r'\leq k$ and $1\leq a<b\leq n$.

For $I\in\binom{[n]}{k}$, let $P_I$ be the quantum
Pl\"ucker coordinate defined in
\eqref{eq:quantum-plucker-coordinate}.
Recall from the introduction that
$\deg X_{r,a}=e_a$ and $\deg P_I=\sum_{a\in I}e_a$
in the content grading.

For subsets $A,B\subseteq[n]$, set
$\operatorname{inv}(A,B)=|\{(a,b)\in A\times B:a>b\}|$.

By the presentation of Taft and Towber
\cite{TaftTowber1991}, as recalled in
\cite[Definition~12 and Equation~(9)]{Lauve2006},
the algebra $A_v$ is generated by the quantum
Pl\"ucker coordinates $P_I$,
$I\in\binom{[n]}{k}$, subject to the quantum Young
symmetry relations
\begin{equation}\label{eq:young-ordinary}
 \sum_{\substack{\Lambda\subset U\setminus V\\
                 |\Lambda|=r}}
 (-v)^{a(U,V,\Lambda)}
 P_{U\setminus\Lambda}P_{\Lambda\cup V}=0,
\end{equation}
where
$a(U,V,\Lambda) = \operatorname{inv}(U\setminus\Lambda,\Lambda) + \operatorname{inv}(\Lambda,V)$,
$1\leq r\leq\min(k,n-k)$,
$U\in\binom{[n]}{k+r}$, and
$V\in\binom{[n]}{k-r}$.
All minor indices in~\eqref{eq:young-ordinary}
are taken in increasing order.

The restriction $\Lambda\cap V=\varnothing$ removes
the terms for which the second minor has repeated
indices. Our convention is obtained from that
of Lauve by replacing his quantum parameter by
$v^{-1}$; see~\cite[Definitions~5 and~6,
Theorem~8]{Lauve2006}.

The relations~\eqref{eq:young-ordinary} generate
the defining ideal in the free associative algebra
on the increasing $k$-subsets of $[n]$.
By~\cite[Proposition~13 and Corollary~14]{Lauve2006},
it suffices to take the relations with $r=1$.
We retain all values of $r$ in order to state the
quantum Young symmetry relations uniformly.

\subsection{The cocycle twist}

Every consecutive quantum Pl\"ucker coordinate
quasi-commutes with every quantum Pl\"ucker coordinate
\cite[Lemma~2.2]{LenaganRussell2008}.
For the non-wrapping consecutive minors $F_a$,
the quasi-commutation formula of Leclerc and Zelevinsky
\cite{LZ98}, see also~\cite[Proposition~3]{Scott2005}, gives
\begin{equation}\label{eq:frozen-commutation}
 F_aP_I=v^{d_a(I)}P_IF_a,\qquad
 d_a(I)=|I\cap[a+k,n]|-|I\cap[1,a-1]|.
\end{equation}
Here an interval $[u,w]$ is understood to be empty
if $u>w$.

Define a skew-symmetric bilinear form
$\Omega:\mathbb Z^n\times\mathbb Z^n\to\mathbb Z$ by
\begin{equation}\label{eq:Omega}
 \Omega(e_i,e_j)=
 \begin{cases}
  -1,&i<j\ \text{and }k\mid(j-i),\\
  1,&i>j\ \text{and }k\mid(i-j),\\
  0,&\text{otherwise}.
 \end{cases}
\end{equation}
For homogeneous elements $f,g\in A_v$, define
\begin{equation}\label{eq:star}
 f\star g=v^{\Omega(\deg f,\deg g)/2}fg,
\end{equation}
and extend this multiplication bilinearly to $A_v$.
By the bilinearity of $\Omega$, the multiplication
\eqref{eq:star} is associative. We denote the resulting
algebra by $A_v^\Omega$.

The form $\Omega$ is related to the bilinear form
$B$ used by Kang, Kashiwara, and Kim
\cite[Lemma~4.5.5 and Definition~4.5.6]{KKK2018}
to modify the monoidal structure on a quotient
category of modules over quiver Hecke algebras.
Let $P=\bigoplus_{a\in\mathbb Z}\mathbb Z\epsilon_a$
with the standard inner product, and let
$S_k(\epsilon_a)=\epsilon_{a+k}$.
Their bilinear form is given by
$B(x,y)=-\sum_{r>0}(S_k^r(x),y)$.
Identifying $\mathbb Z^n$ with the sublattice of $P$
spanned by $\epsilon_1,\ldots,\epsilon_n$, we have
$\Omega(x,y)=B(x,y)-B(y,x)$ for $x,y\in\mathbb Z^n$.
While their construction uses $B$ directly, our
cocycle twist of the quantum Grassmannian uses
the skew-symmetric form $\Omega$.

\begin{proposition}\label{prop:center}
For every $1\leq a\leq s$, the element $F_a$
is central in $A_v^\Omega$.
\end{proposition}

\begin{proof}
The interval $C_a$ contains exactly one integer
in each residue class modulo $k$.
It follows from~\eqref{eq:Omega} that
\[
 \Omega(\deg F_a,e_j)
 =
 \begin{cases}
  1,&j<a,\\
  0,&a\leq j\leq a+k-1,\\
  -1,&j>a+k-1.
 \end{cases}
\]
Consequently, for every $I\in\binom{[n]}{k}$,
$\Omega(\deg F_a,\deg P_I) = |I\cap[1,a-1]|-|I\cap[a+k,n]| =-d_a(I)$.
Combining~\eqref{eq:frozen-commutation} with
\eqref{eq:star}, we obtain
\[
 F_a\star P_I = v^{-d_a(I)/2}F_aP_I = v^{d_a(I)/2}P_IF_a = P_I\star F_a.
\]
Since the quantum Pl\"ucker coordinates generate
$A_v^\Omega$, the element $F_a$ is central.
\end{proof}

We will specialize only the non-wrapping consecutive
minors $F_1,\ldots,F_s$.
The remaining $k-1$ cyclic consecutive minors are
neither specialized to scalars nor inverted when forming the
central quotient.

For subsets $A,B\subseteq[n]$, we use the notation
$\Omega(A,B) = \Omega\left(\sum_{a\in A}e_a,\sum_{b\in B}e_b\right)$.
By the homogeneous presentation~\eqref{eq:young-ordinary}, the twisted Young symmetry relations are obtained
by rescaling the terms of~\eqref{eq:young-ordinary}.
The algebra $A_v^\Omega$ is generated by the quantum
Pl\"ucker coordinates $P_I$,
$I\in\binom{[n]}{k}$, subject to the relations
\begin{equation}\label{eq:young-twisted}
 \sum_{\substack{\Lambda\subset U\setminus V\\
                 |\Lambda|=r}}
 (-v)^{a(U,V,\Lambda)}
 v^{-\Omega(U\setminus\Lambda,\Lambda\cup V)/2}
 P_{U\setminus\Lambda}\star P_{\Lambda\cup V}
 =0,
\end{equation}
where $r,U,V$ range over the same indices as
in~\eqref{eq:young-ordinary}.

For any $\boldsymbol\zeta=(\zeta_1,\ldots,\zeta_s)\in\bbk^s$,
adjoining the relations $F_a=\zeta_a$, $1\leq a\leq s$,
gives the defining relations of the quotient
$A_v^\Omega/\langle F_a-\zeta_a:1\leq a\leq s\rangle$.

Indeed, twisting the content-graded free algebra on the Pl\"ucker
symbols preserves its freeness and the homogeneous ideal of relations
as a vector space. Each product of homogeneous elements is rescaled by
\eqref{eq:star}; thus the homogeneous generators of the defining ideal
in~\eqref{eq:young-ordinary} become exactly~\eqref{eq:young-twisted}.
Adjoining the central scalar relations gives the quotient presentation.

\subsection{Central specialization and the bar involution}
\label{subsec:central-specialization}

By Proposition~\ref{prop:center}, for any
$\boldsymbol\zeta=(\zeta_1,\ldots,\zeta_s)\in\bbk^s$,
we may consider the central quotient
$A_v^\Omega(\boldsymbol\zeta) = A_v^\Omega/ \langle F_a-\zeta_a:1\leq a\leq s\rangle$.

Let $\chi:\mathbb Z^n\to\mathbb Z$ be an additive map.
Consider the central quotient
$B_\chi = A_v^\Omega/ \langle F_a-v^{\chi(C_a)/2}: 1\leq a\leq s\rangle$,
and let $\pi_\chi:A_v^\Omega\to B_\chi$ be the quotient map.
Here and below, for a subset $I\subseteq[n]$,
we write
$\chi(I)=\chi\left(\sum_{a\in I}e_a\right)$.
The specialization in Theorem~\ref{thm:intro-A} is
obtained by taking the map $\chi$ defined
in the introduction.

More generally, every specialization of
the form
$F_a=v^{d_a/2}$, $d_a\in\mathbb Z$,
can be obtained in this way.
Indeed, the equations
$\sum_{j=a}^{a+k-1}\chi(e_j)=d_a$, $1\leq a\leq s$,
can be solved recursively over $\mathbb Z$
after choosing the first $k-1$ values of
$\chi(e_j)$.

For a homogeneous element $f\in A_v^\Omega$,
put
$\mathcal D(f) = v^{-\chi(\deg f)/2}\pi_\chi(f)$,
and extend $\mathcal D$ linearly.
Since $\chi$ is additive, $\mathcal D$ is an algebra homomorphism.
Moreover, $\mathcal D(F_a)=1$ for $1\leq a\leq s$.
Let $B_0=A_v^\Omega/\langle F_a-1:1\leq a\leq s\rangle$,
with quotient map $\pi_0$. Multiplication of each homogeneous
$f\in A_v^\Omega$ by $v^{-\chi(\deg f)/2}$ defines an algebra
automorphism of $A_v^\Omega$, with inverse given by multiplication
by $v^{\chi(\deg f)/2}$. It sends $F_a-1$ to
$v^{-\chi(C_a)/2}(F_a-v^{\chi(C_a)/2})$ and hence induces an
isomorphism $D_\chi:B_0\to B_\chi$ satisfying
$D_\chi(\pi_0(f))=\mathcal D(f)$.

The standard matrix bar $\overline{\phantom{x}}$ is the semilinear
anti-involution of the quantum matrix algebra fixing the matrix
generators and sending $v^{1/2}$ to $v^{-1/2}$. Its restriction to
the quantum Grassmannian is denoted by the same symbol.
See~\cite{Zhang2004} and~\cite[Remark~17]{Brundan2005}.
In particular,
$\overline{P_I}=P_I$ for all $I\in\binom{[n]}{k}$.
Since the bar involution preserves the content
grading and $\Omega$ is skew-symmetric, for
homogeneous elements $f,g\in A_v$ we have
\begin{equation}\label{eq:twisted-bar}
 \overline{f\star g} = v^{-\Omega(\deg f,\deg g)/2} \overline g\,\overline f = \overline g\star\overline f.
\end{equation}
Thus the bar involution of $A_v$ is also a
semilinear anti-involution of $A_v^\Omega$.
It descends to the quotient
$A_v^\Omega(\boldsymbol\zeta)$ whenever
$\overline{\zeta_a}=\zeta_a$ for all $a$.

For the specialization defining $B_\chi$, we introduce the modified bar by setting
$\overline f^{\,\chi} = v^{-\chi(\deg f)}\overline f$
for homogeneous elements $f\in A_v^\Omega$,
and extend this map semilinearly.
For homogeneous elements $f,g\in A_v^\Omega$,
the additivity of $\chi$ gives
\[
 \overline{f\star g}^{\,\chi}
 =v^{-\chi(\deg f)-\chi(\deg g)}\overline g\star\overline f
 =\overline g^{\,\chi}\star\overline f^{\,\chi}.
\]
Moreover, since $\overline v=v^{-1}$
and the bar involution preserves the content grading,
we have
$\overline{(\overline f^{\,\chi})}^{\,\chi}=v^{\chi(\deg f)}v^{-\chi(\deg f)}f=f$.
Hence $\overline{\phantom{x}}^{\,\chi}$ is a
semilinear anti-involution of $A_v^\Omega$.

For every $1\leq a\leq s$, we have
\[
 \overline{ F_a-v^{\chi(C_a)/2} }^{\,\chi} = v^{-\chi(C_a)}F_a - v^{-\chi(C_a)/2} = v^{-\chi(C_a)} \left(F_a-v^{\chi(C_a)/2}\right).
\]
Therefore the ideal defining $B_\chi$ is invariant
under $\overline{\phantom{x}}^{\,\chi}$, which
induces a bar involution on $B_\chi$.
We denote this induced involution by the same symbol.
For homogeneous $f\in A_v^\Omega$, we have
\[
 \overline{\mathcal D(f)}^{\,\chi}
 =v^{\chi(\deg f)/2}v^{-\chi(\deg f)}\pi_\chi(\overline f)
 =\mathcal D(\overline f).
\]
By semilinearity, the equality holds for all $f\in A_v^\Omega$.
In particular, $D_\chi$ intertwines the bar involutions on
$B_0$ and $B_\chi$.

\section{Dual canonical and standard bases of the central quotient}
\label{sec:canonical-specialization}

We first recall the normalized standard monomials and the dual
canonical basis of the quantum Grassmannian. We then establish
the dual canonical basis of the central quotient by a free module
decomposition over the central subalgebra generated by the non-wrapping
consecutive minors. Finally, we compare the dominance order on tableaux
with the partial order on dominant monomials and use the resulting
triangularity to construct a standard basis of the quotient.

\subsection{Standard and dual canonical bases of the quantum Grassmannian}
\label{subsec:homogeneous-standard-monomials}

We use the standard monomial and dual canonical bases of the
quantum Grassmannian with the following normalizations.
For earlier work on IC bases of quantum linear groups and their
coordinate algebras, see~\cite{Du1994,Du1996}.

Let $T\in\SSYT(k,[n])$ be a tableau of width $w$,
with columns $T_1,\ldots,T_w$ from left to right.
Let $\gamma_a$ be the number of occurrences of $a$
in $T$. The corresponding row and column degrees
in the quantum matrix algebra are
$\mu=(w,\ldots,w)$ and $\gamma=(\gamma_1,\ldots,\gamma_n)$;
the pair $(\mu,\gamma)$ is the bidegree.
The normalization exponent in~\cite[Equation~(2.9) and Lemma~4.13]{LOT2026} is given by
\begin{equation}\label{eq:main-beta}
 \beta(T)=\frac{(\mu,\mu)-(\gamma,\gamma)}{2}
 =k\binom{w}{2}-\sum_{a=1}^n\binom{\gamma_a}{2}
 =\sum_{1\leq a<b\leq w}(k-|T_a\cap T_b|),
\end{equation}
where $(\cdot,\cdot)$ denotes the standard inner
product. The second equality follows from
$\sum_{a=1}^n\gamma_a=kw$, and the third follows
by counting pairs of columns containing each
entry of $[n]$. In particular, $\beta(T)$ depends
only on the bidegree $(\mu,\gamma)$.

In the convention of~\cite[Section~4.1 and Equation~(4.53)]{LOT2026},
the standard monomial and its normalization are given by
$P_T=P_{T_1}\cdots P_{T_w}$ and $\widetilde P_T=v^{-\beta(T)/2}P_T$.
For the empty tableau, both expressions are $1$.
The product defining $P_T$ is taken in $A_v$.

Let $\psi_B$ be the bar involution in the convention of
\cite[Theorem~16]{Brundan2005}, with the parameter $q$ there
identified with $v$. It is semilinear under $v^{1/2}\mapsto v^{-1/2}$,
fixes the matrix generators,
and satisfies, for $f,g$ of bidegrees $(\mu,\gamma)$ and
$(\mu',\gamma')$, respectively,
\[
 \psi_B(fg)=v^{(\mu,\mu')-(\gamma,\gamma')}
 \psi_B(g)\psi_B(f).
\]
In particular, for $f$ of the same bidegree as $T$,
\[
 \psi_B(f)=v^{\beta(T)}\overline f,
\]
where $\overline f$ is the matrix bar of
Section~\ref{sec:twist}; see~\cite[Remark~17]{Brundan2005}.
By~\cite[Theorem~26 and Remark~27]{Brundan2005}, the standard
monomials $P_S$ of this bidegree form a basis, and the dual
canonical basis element $L_T^B$ in this convention is uniquely
characterized by
\[
 \psi_B(L_T^B)=L_T^B,\qquad
 L_T^B\in P_T+\sum_{S\ne T} v^{-1}\mathbb Z[v^{-1}]P_S,
\]
where the sum is over semistandard tableaux $S$ of the same
shape and content as $T$. In the notation of
\cite[Theorem~26]{Brundan2005}, $P_T$ is the product of flag
minors $V_T$ (all of size $k$ here), and $L_T^B$ is indexed
by the row-sorted tableau $R(T)$, which is $T$ in our
rectangular semistandard convention.
Since $\beta(S)=\beta(T)$ for all such $S$, the rescaled element
$v^{-\beta(T)/2}L_T^B$ is invariant under the matrix bar and lies in
$\widetilde P_T+\sum_{S\ne T} v^{-1}\mathbb Z[v^{-1}]\widetilde P_S$.
Thus, in the normalization of~\cite[Proposition~B.1, Equation~(B.3)]{LOT2026},
\begin{equation}\label{eq:main-Brundan-rescaling}
 b_T=v^{-\beta(T)/2}L_T^B.
\end{equation}

The elements $b_T$ belong to the quantum matrix algebra over
$\mathbb Z[v^{1/2},v^{-1/2}]$ and specialize at $v^{1/2}=1$
to the dual canonical basis elements of the classical Grassmannian
coordinate ring indexed by the same tableaux; see~\cite[Equations~(6.1)--(6.2)]{Brundan2005}
and~\cite[Definition~2.4]{LOT2026}, together with
\eqref{eq:main-Brundan-rescaling}.
Thus the normalized standard basis and the dual canonical basis
of $A_v$ are, respectively,
\[
 \{\widetilde P_T:T\in\SSYT(k,[n])\},\qquad
 \{b_T:T\in\SSYT(k,[n])\}.
\]
The normalized standard monomials need not be bar invariant.
Their triangular relation with the dual canonical basis is
described in Section~\ref{sec:main-support}.

\subsection{The dual canonical basis of the central quotient}

Since $A_v$ and $A_v^\Omega$ have the same underlying vector space,
the elements $b_T$ also form a basis of $A_v^\Omega$.
We first establish the multiplication formula for the non-wrapping
consecutive minors and deduce from it a basis of each central
quotient indexed by reduced tableaux.

By the quantum multiplication property used in the proof of
\cite[Lemma~3.11]{FeiginLi2026}, which follows from
\cite{Caldero2003}, the product $F_ab_T$ is a power of $v^{1/2}$
times a dual canonical basis element. Its classical specialization
identifies that element as $b_{T\cup C_a}$ by the classical
multiplication formula in the same lemma. The twisted product $F_a\star b_T$ is also a power of $v^{1/2}$
times $b_{T\cup C_a}$. Since $F_a$ is central in the twisted algebra
and both factors are bar invariant, the product is bar invariant
by~\eqref{eq:twisted-bar}. As $v$ is an indeterminate and
$b_{T\cup C_a}$ is bar invariant, this scalar power must be $1$.
Thus, for every $T\in\SSYT(k,[n])$ and $1\leq a\leq s$,
\begin{equation}\label{eq:exact-frozen}
 F_a\star b_T=b_{T\cup C_a}=b_T\star F_a.
\end{equation}

Equivalently, in the quantum Grassmannian $A_v$,
we have
\begin{equation}\label{eq:untwisted-frozen}
 F_ab_T
 =
 v^{d_a(T)/2}b_{T\cup C_a},
 \qquad
 b_TF_a
 =
 v^{-d_a(T)/2}b_{T\cup C_a},
\end{equation}
where
$d_a(T)=\sum_{j=1}^w d_a(T_j)$
and $d_a(T_j)$ is defined
in~\eqref{eq:frozen-commutation}.

By~\cite[Lemma~4.2]{FeiginLi2026}, every tableau
$T\in\SSYT(k,[n])$ admits a unique decomposition
\begin{equation}\label{eq:reduction}
 T
 =
 T_{\mathrm{red}}
 \cup C_1^{m_1(T)}
 \cup\cdots\cup C_s^{m_s(T)},
\end{equation}
where $T_{\mathrm{red}}$ is reduced and
\begin{equation}\label{eq:ma}
 m_a(T)
 =
 \min_{1\leq r\leq k}h_{r,a+r-1}(T),
 \qquad 1\leq a\leq s.
\end{equation}
Here $h_{r,x}(T)$ denotes the multiplicity of $x$
in row $r$ of $T$, and $C_a^{m_a(T)}$ denotes
the rowwise union of $m_a(T)$ copies of $C_a$.

Combining this decomposition with~\eqref{eq:exact-frozen}, we obtain a quantum
analogue of the free module decomposition
in~\cite[Theorem~4.3]{FeiginLi2026}.

\begin{proposition}\label{thm:free}
The elements $F_1,\ldots,F_s$ are algebraically
independent in $A_v^\Omega$.
Let
$Z=\bbk[F_1,\ldots,F_s]\subseteq A_v^\Omega$.
Then $A_v^\Omega$ is a free $Z$-module with basis
$\{b_R:R\text{ is reduced}\}$.
More precisely,
\begin{equation}\label{eq:free}
 A_v^\Omega
 =
 \bigoplus_{R\text{ reduced}} Zb_R.
\end{equation}

For any
$\boldsymbol\zeta=(\zeta_1,\ldots,\zeta_s)\in\bbk^s$,
let $\pi_{\boldsymbol\zeta}$ be the quotient map
\[
 \pi_{\boldsymbol\zeta}:
 A_v^\Omega\longrightarrow
 A_v^\Omega(\boldsymbol\zeta)
 =
 A_v^\Omega/
 \langle F_a-\zeta_a:1\leq a\leq s\rangle.
\]
Then
$\{\pi_{\boldsymbol\zeta}(b_R): R\text{ is reduced}\}$
is a $\bbk$-basis of $A_v^\Omega(\boldsymbol\zeta)$.
Moreover, for every $T\in\SSYT(k,[n])$, we have
\begin{equation}\label{eq:quotient-action}
 \pi_{\boldsymbol\zeta}(b_T)
 =
 \left(\prod_{a=1}^s\zeta_a^{m_a(T)}\right)
 \pi_{\boldsymbol\zeta}(b_{T_{\mathrm{red}}}).
\end{equation}
\end{proposition}

\begin{proof}
By~\eqref{eq:exact-frozen}
and~\eqref{eq:reduction}, we have
$b_T = F_1^{\star m_1(T)}\star\cdots\star F_s^{\star m_s(T)} \star b_{T_{\mathrm{red}}}$.
By the uniqueness of the decomposition~\eqref{eq:reduction},
rowwise union with $m_a$ copies of $C_a$ gives a bijection
from pairs consisting of a reduced tableau and an $s$-tuple
of nonnegative integers to $\SSYT(k,[n])$.
Consequently, the elements
\[
 F_1^{\star m_1}\star\cdots\star
 F_s^{\star m_s}\star b_R,
 \qquad
 m_1,\ldots,m_s\geq 0,\quad R\text{ reduced},
\]
form a $\bbk$-basis of $A_v^\Omega$.

Taking $R$ to be the empty tableau proves that
$F_1,\ldots,F_s$ are algebraically independent.
Since these elements are central by
Proposition~\ref{prop:center}, allowing arbitrary
reduced tableaux gives the free $Z$-module
decomposition~\eqref{eq:free}.

The quotient $A_v^\Omega(\boldsymbol\zeta)$ is
obtained by base change along the algebra
homomorphism
$Z\to\bbk$, $F_a\mapsto\zeta_a$.
The freeness of $A_v^\Omega$ over $Z$ gives the
asserted basis of the quotient.
Applying $\pi_{\boldsymbol\zeta}$ to the
multiplication formula for $b_T$ gives
\eqref{eq:quotient-action}.
\end{proof}

For every reduced tableau $T$, put
$\widehat b_T=\mathcal D(b_T)=v^{-\chi(T)/2}\pi_\chi(b_T)$.
By Proposition~\ref{thm:free}, the set
\[
 \mathcal B_\chi=\{\widehat b_T:T\in\SSYT(k,[n])\text{ is reduced}\}
\]
is a $\bbk$-basis of $B_\chi$. We call it the dual canonical
basis of the central quotient $B_\chi$.

The homomorphism $\mathcal D$ sends each $F_a$ to $1$.
Applying it to the factorization of $b_V$ in the proof of
Proposition~\ref{thm:free} gives, for every $V\in\SSYT(k,[n])$,
\begin{equation}\label{eq:normalized-reduction}
 \mathcal D(b_V)
 =
 \widehat b_{V_{\mathrm{red}}}.
\end{equation}

By the compatibility of $\mathcal D$ with the bar involutions
proved in Section~\ref{subsec:central-specialization}, the basis
elements $\widehat b_T$ are invariant under the modified bar involution.

\subsection{Quantum straightening and the dominance order}
\label{sec:main-support}

We compare the dominance order on tableaux with the partial order
on monomials in $q$-character theory. This comparison allows us
to express quantum straightening and the homogeneous basis transition
in terms of dominant monomials.

For $i\in\mathsf I$ and $p\in\mathbb Z$, let
\begin{equation}\label{eq:simple-root-monomials}
 A_{i,p}=Y_{i,p-1}Y_{i,p+1}
 \prod_{\substack{j\in\mathsf I\\ |i-j|=1}}Y_{j,p}^{-1};
\end{equation}
see~\cite{FR99}.
The usual partial order on monomials is given by
\[
 M'\leq M
 \quad\Longleftrightarrow\quad
 M(M')^{-1}=\prod_{i,p}A_{i,p}^{u_{i,p}},
 \qquad u_{i,p}\in\mathbb Z_{\geq0};
\]
see~\cite[Theorem~4.1]{FM01}
and~\cite[Definition~2.9]{Nakajima04}.
Here only finitely many of the exponents $u_{i,p}$
are nonzero.
We write $M'<M$ if $M'\leq M$ and $M'\ne M$.

For a tableau $T$, write $h(T)=(h_{r,a}(T))_{r,a}$, where
$h_{r,a}(T)$ is the multiplicity of $a$ in row $r$ of $T$.
For tableaux $S,T$ of the same content, the dominance order
of~\cite[Definition~3.2]{CDFL20} and~\cite[Equation~(4.2)]{LOT2026}
is expressed in terms of row multiplicities by
\[
 h(S)\unlhd h(T)
 \quad\Longleftrightarrow\quad
 \sum_{i=1}^r\sum_{x=1}^a h_{i,x}(S)
 \leq \sum_{i=1}^r\sum_{x=1}^a h_{i,x}(T)
 \quad\text{for all }r,a.
\]
Indeed, the length of row $i$ of the restriction of $T$ to entries
in $[a]$ is $\sum_{x=1}^a h_{i,x}(T)$.
By~\cite[Proposition~3.28]{CDFL20}, we have
\begin{equation}\label{eq:tableau-monomial-order}
 h(S)\unlhd h(T)\quad\Longleftrightarrow\quad M_S\leq M_T.
\end{equation}

For a sequence of increasing columns $W=(I_1,\ldots,I_d)$,
let $M_W=\prod_{j=1}^d M_{I_j}$ and let $T^{\max}$ be the
semistandard tableau obtained by sorting each row of $W$.
By the column formula in the introduction, the exponent of
$Y_{k-r,k+r-2a}$ in $M_W$ is the number of entries strictly
less than $a$ in row $r$ minus the number of entries at most
$a$ in row $r+1$. Thus $M_W$ depends only on the row
multiplicities, and $M_{T^{\max}}=M_W$.
By~\cite[Lemma~4.7]{LOT2026}, every tableau $T$ occurring in the
standard monomial expansion of $P_{I_1}\cdots P_{I_d}$ satisfies
$h(T)\unlhd h(T^{\max})$, and hence $M_T\leq M_W$ by
\eqref{eq:tableau-monomial-order}. Equality implies
$h(T)=h(T^{\max})$ by the same equivalence and antisymmetry
of the dominance order. Since a semistandard tableau is determined
by its row multiplicities, equality holds only for $T=T^{\max}$.
The cocycle twist rescales the product by a nonzero scalar and
therefore preserves this support.

In particular, $M_T$ distinguishes tableaux of the same content.
We will also use the fact that, for $M\in\mathcal P_\ell^+$, the set
$\{M'\in\mathcal P_\ell^+:M'\leq M\}$ is finite.
Indeed, define an additive function on Laurent monomials by
$\operatorname{wt}(\prod_{i,p}Y_{i,p}^{u_{i,p}})
=\sum_{i,p}i(k-i)u_{i,p}$.
Since $\operatorname{wt}(A_{i,p})=2$, we have
$\operatorname{wt}(M')\leq\operatorname{wt}(M)$ whenever $M'\leq M$.
The function $\operatorname{wt}$ is positive on every variable
$Y_{i,p}$, and only finitely many such variables occur in the
monomials of $\mathcal P_\ell^+$.
Thus each exponent of $M'$ is bounded by the fixed number
$\operatorname{wt}(M)$, proving the finiteness assertion.

By the unitriangular basis transitions in
\cite[Theorem~26]{Brundan2005} and~\cite[Proposition~B.1]{LOT2026},
together with~\eqref{eq:tableau-monomial-order}, for every
$T\in\SSYT(k,[n])$ we have
\begin{equation}\label{eq:homogeneous-root-inverse}
 \widetilde P_T
 =b_T+
 \sum_{\substack{V\text{ of the same content as }T\\
                 M_V<M_T}}
 c_{T,V}(v)b_V,
 \qquad c_{T,V}(v)\in v^{-1}\mathbb Z[v^{-1}].
\end{equation}
Here $c_{T,V}(v)=k^*_{T,V}(v^{-1})$ in the notation of~\cite{Brundan2005}.
On each fixed bidegree, the strictly triangular part of the
transition matrix is nilpotent. Its inverse is therefore also
unitriangular, with off-diagonal coefficients in
$v^{-1}\mathbb Z[v^{-1}]$.

\subsection{The standard basis of the central quotient}
\label{subsec:source-standard-basis}

Using the preceding triangularity, we construct a standard basis
of $B_\chi$ from ordered products of fundamental Pl\"ucker coordinates.

The fundamental columns $I_{b,i}$ in~\eqref{eq:fundamental-columns}
form a chain in the componentwise order (see~\cite[Lemma~3.12]{CDFL20}).
Indeed, writing $r=k-i$, the entry of $I_{b,i}$ in row $j$ is $b+j-1$
if $j\leq r$ and $b+j$ if $j>r$. For fixed $b$, increasing $i$ increases the column componentwise.
If $b<c$, every entry of $I_{b,i}$ is at most the corresponding
entry of $I_{c,j}$ for any $i,j$. Thus the lexicographic order
on the pairs $(b,i)$ gives the componentwise chain order.
Throughout, the lexicographic order on fundamental columns refers to
this order on the pairs $(b,i)$.

Put $y_{b,i}=Y_{i,i-2b}$ and $X_{b,i}=\mathcal D(P_{I_{b,i}})$.
For $\mathbf u=(u_{b,i})\in\mathbb N^{m(k-1)}$, put
$M_{\mathbf u}=\prod_{b,i}y_{b,i}^{u_{b,i}}$ and let $U(\mathbf u)$
be the tableau obtained by arranging $u_{b,i}$ copies of $I_{b,i}$
in lexicographic order. In particular, $M_{U(\mathbf u)}=M_{\mathbf u}$.
Define
\begin{equation}\label{eq:source-standard-definition}
 S_{\mathbf u}=\mathcal D(\widetilde P_{U(\mathbf u)}).
\end{equation}
Let $I_1,\ldots,I_h$ be the columns of $U(\mathbf u)$ in this order,
and write $X_a=\mathcal D(P_{I_a})$. By~\eqref{eq:main-beta} and
\eqref{eq:star},
\begin{equation}\label{eq:source-chain-standard}
 S_{\mathbf u}=v^{-\frac12\sum_{a<b}\Lambda(I_a,I_b)}X_1\cdots X_h,
 \qquad \Lambda(I,I')=k-|I\cap I'|+\Omega(I,I').
\end{equation}
For the empty tableau this product is $1$. This normalization is determined by the standard monomials
of $A_v$ and the cocycle twist. Since each $P_{I_{b,i}}$
is bar invariant, the elements $X_{b,i}$ are invariant under the
modified bar involution by Section~\ref{subsec:central-specialization}.

Recall that the correspondence between
tableaux and dominant monomials induces
a bijection
$\{T\in\SSYT(k,[n]):T\text{ is reduced}\}\to\mathcal P_\ell^+$, given by $T\mapsto M_T$;
see~\cite[Lemma~3.9, Lemma~3.13
and Equation~(3.5)]{CDFL20}.

For each reduced tableau $T$, let
$\mathbf u(T)$ be determined by
$M_T= \prod_{b=1}^m\prod_{i=1}^{k-1} y_{b,i}^{u_{b,i}(T)}$.
Since $M_{U(\mathbf u(T))}=M_T$ and the monomial associated
with a trivial column is one, the bijection above gives
$U(\mathbf u(T))_{\mathrm{red}}=T$.

Put
\[
 \mathcal S_\chi=\{S_{\mathbf u}:\mathbf u\in\mathbb N^{m(k-1)}\}.
\]
We show that $\mathcal S_\chi$ is a basis of $B_\chi$ and
compare it with the dual canonical basis $\mathcal B_\chi$.
Both sets are indexed by the dominant monomials in $\mathcal P_\ell^+$,
via $\mathbf u\mapsto M_{\mathbf u}$ and $T\mapsto M_T$, respectively.

\begin{proposition}\label{thm:general-source-pbw}
For every reduced tableau $T$, we have
\begin{equation}\label{eq:source-standard-triangularity}
 \widehat b_T-S_{\mathbf u(T)}
 \in
 \sum_{M_{\mathbf u'}<M_T}
 v^{-1}\mathbb Z[v^{-1}]S_{\mathbf u'}.
\end{equation}
In particular, $\mathcal S_\chi$ is a $\bbk$-basis of $B_\chi$,
which we call the standard basis of the central quotient.
The standard basis $\mathcal S_\chi$ and the dual canonical
basis $\mathcal B_\chi$ span the same
$\mathbb Z[v^{1/2},v^{-1/2}]$-module.
Moreover, $\widehat b_T$ is the unique element invariant under
the modified bar involution satisfying~\eqref{eq:source-standard-triangularity}.
Consequently, $B_\chi$ is generated by
the elements
$X_{b,i}$, $1\leq b\leq m, 1\leq i<k$.

\end{proposition}

\begin{proof}
Let $\mathbf u\in\mathbb N^{m(k-1)}$ and
put $U=U(\mathbf u)$.
Applying $\mathcal D$ to~\eqref{eq:homogeneous-root-inverse} for $T=U$
and using
\eqref{eq:source-standard-definition} and~\eqref{eq:normalized-reduction}, we
obtain
\begin{equation}\label{eq:source-standard-expansion}
 S_{\mathbf u}
 =
 \widehat b_{U_{\mathrm{red}}}
 +
 \sum_{M_V<M_{\mathbf u}}
 c_{U,V}(v)\widehat b_{V_{\mathrm{red}}}.
\end{equation}

Removing trivial columns does not change
the associated monomial.
Thus
$M_{U_{\mathrm{red}}}=M_{\mathbf u},  M_{V_{\mathrm{red}}}=M_V$.
In particular, every term in the sum
in~\eqref{eq:source-standard-expansion}
has an associated monomial strictly
smaller than $M_{\mathbf u}$.
After collecting terms indexed by the
same reduced tableau, we obtain a
unitriangular transition between the
elements $S_{\mathbf u}$ and the basis
elements $\widehat b_T$, with
off-diagonal coefficients in
$v^{-1}\mathbb Z[v^{-1}]$.

As observed in Section~\ref{sec:main-support},
for every $M\in\mathcal P_\ell^+$, the set
$\{M'\in\mathcal P_\ell^+:M'\leq M\}$
is finite.
We may therefore invert the transition
on each such finite set of monomials.
Its inverse is again unitriangular,
with off-diagonal coefficients in
$v^{-1}\mathbb Z[v^{-1}]$.
Since $\mathcal B_\chi$ is a basis, this proves that
$\mathcal S_\chi$ is a basis and establishes
\eqref{eq:source-standard-triangularity}. The two transitions
also show that the bases span the same integral module.

For uniqueness, the difference between two bar-invariant elements
satisfying~\eqref{eq:source-standard-triangularity} is a linear
combination of lower standard basis elements with coefficients
in $v^{-1}\mathbb Z[v^{-1}]$. By~\eqref{eq:source-standard-expansion},
its coefficients in the dual canonical basis also belong to
$v^{-1}\mathbb Z[v^{-1}]$. Since that basis is bar invariant,
these coefficients must be invariant under $v\mapsto v^{-1}$,
and hence vanish.

Since each $S_{\mathbf u}$ is a power
of $v^{1/2}$ times a product of the
generators $X_{b,i}$, these
generators generate $B_\chi$.

\end{proof}

\section{Standard and canonical bases of the quantum Grothendieck ring}
\label{subsec:general-spectral-words}

We recall the quantum torus realization and the two bases of
$K_t(\mathcal C_\ell)$, using the decreasing spectral order for
normalized standard products. We also determine the dominant monomials
in products of two fundamental characters for use in the presentation
of the quantum Grothendieck ring.

\subsection{The quantum torus and the standard basis}

Let $\mathsf I=[1,k-1]$.
Recall that the quantum Cartan matrix $C(z)$ of
type $A_{k-1}$ has entries~\cite[Section~2.1]{FR99}
\[
 C_{ij}(z)=
 \begin{cases}
  z+z^{-1},&i=j,\\
  -1,&|i-j|=1,\\
  0,&\text{otherwise}.
 \end{cases}
\]
Following~\cite[Section~2.4]{HL2015}, write
$C(z)^{-1}_{ij}=\sum_{h\geq1}c^+_{ij}(h)z^h$
for the expansion at $z=0$, and put
$c^+_{ij}(h)=0$ for $h\leq0$.

Let $\mathcal M$ denote the multiplicative group of Laurent
monomials in the variables $Y_{i,p}$, $i\in\mathsf I$, $p\in\mathbb Z$,
regarded as a lattice
via their exponent vectors. Define a bilinear form
$N:\mathcal M\times\mathcal M\to\mathbb Z$ by
$N(Y_{i,p},Y_{j,p'}) =c^+_{ij}(p'-p-1)-c^+_{ij}(p'-p+1)$.
Put $J(M,M')=N(M',M)-N(M,M')$.
Then $J$ is the skew-symmetric form
$\mathcal N$ in~\cite[Equations~(3.1)
and~(3.5)]{FHOOPos2022} for type $A_{k-1}$.

In the notation of~\cite[Section~3.2]{HL2015},
the quantum torus $\mathcal Y_t$ is the free
$\mathbb Z[t^{\pm1/2}]$-module with basis $\mathcal M$.
Its multiplication is given by
\begin{equation}\label{eq:quantum-torus-product}
 M*M'=t^{J(M,M')/2}MM',
 \qquad M,M'\in\mathcal M,
\end{equation}
where $MM'$ denotes the commutative product of Laurent
monomials. The bar involution is the semilinear
anti-involution determined by
$\overline{t^{1/2}}=t^{-1/2}$ and $\overline M=M$
for $M\in\mathcal M$. In particular,
$M*M'=t^{J(M,M')}M'*M$.
All quantum products of characters are taken with
respect to $*$. The bar preserves $K_t(\mathcal C_\ell)$, and its
restrictions to this ring and its scalar extension $K$ are denoted by
the same symbol (with $v=t$ on $K$).

Extend the notation $y_{b,i}=Y_{i,i-2b}$
to every $b\in\mathbb Z$ and $1\leq i<k$,
and put $y_{b,0}=y_{b,k}=1$.
Let $G_{b,i}=\chi_{q,t}(L(y_{b,i}))$ and
$G_{b,0}=G_{b,k}=1$.

The set
$\mathcal X=\{(i,i-2b):i\in\mathsf I,\
1\leq b\leq m\}$
is convex in the repetition quiver of type
$A_{k-1}$.
Indeed, along an arrow
$(i,p)\to(j,p+1)$ with $|i-j|=1$,
the value of $p-i$ changes by
$1+i-j\in\{0,2\}$.
Hence every directed path whose endpoints
belong to $\mathcal X$ is contained in
$\mathcal X$.
This is the simply-laced case of the convexity
condition in~\cite[Definition~5.5]{FHOOPos2022}.

By~\cite[Lemma~5.6 and Proposition~5.7]{FHOOPos2022},
the quantum Grothendieck ring $K_t(\mathcal C_\ell)$
is a subalgebra of the quantum Grothendieck ring
of $\mathcal C_{\mathbb Z}$ and has bases
$\{\chi_{q,t}(L(M)):M\in\mathcal P_\ell^+\}$ and
$\{E_t(M):M\in\mathcal P_\ell^+\}$,
where $E_t(M)$ denotes the $(q,t)$-character of
the standard module associated with $M$; its normalization
as an ordered product of fundamental characters is given below.
The former is the canonical basis and the latter is the standard
basis in increasing spectral order.

Following~\cite[Section~3.3 and
Theorem~3.7]{FHOOPos2022}, for
$M=y_1\cdots y_h$, with
$y_a=Y_{i_a,p_a}$ and $p_1\leq\cdots\leq p_h$,
we have
\[
 E_t(M)=t^{-\frac12\sum_{a<b}J(y_a,y_b)}
 \chi_{q,t}(L(y_1))*\cdots*\chi_{q,t}(L(y_h)).
\]
The order among factors with the same spectral
parameter is immaterial.
Applying the bar involution gives
\[
 \overline{E_t(M)}
 =t^{\frac12\sum_{a<b}J(y_a,y_b)}
 \chi_{q,t}(L(y_h))*\cdots*\chi_{q,t}(L(y_1)).
\]
Thus $\overline{E_t(M)}$ is the normalized product
of fundamental $(q,t)$-characters in decreasing
order of their spectral parameters.
This is the standard basis convention used in
\cite[Equation~(3.1)]{Fujita2024}.

For the standard basis in decreasing spectral order, we use the following notation:
for $M=y_1\cdots y_h$ with factors ordered by decreasing spectral
parameter (ties may be resolved by the lexicographic order), put
\begin{equation}\label{eq:kl-standard-product}
 \mathcal E(M)=\overline{E_t(M)}
 =t^{-\frac12\sum_{a<b}J(y_a,y_b)}
 \chi_{q,t}(L(y_1))*\cdots*\chi_{q,t}(L(y_h)).
\end{equation}
For $M=1$, both standard products are $1$. The standard basis
used in the rest of the paper and the canonical basis are, respectively,
\[
 \mathcal S_K=\{\mathcal E(M):M\in\mathcal P_\ell^+\},\qquad
 \mathcal B_K=\{\chi_{q,t}(L(M)):M\in\mathcal P_\ell^+\}.
\]
These are bases of $K_t(\mathcal C_\ell)$ over $\mathbb Z[t^{\pm1/2}]$
and of $K$ over $\bbk$.

\subsection{The canonical basis and triangularity}

By~\cite[Theorem~8.1]{Nakajima04},
\cite[Theorem~6.9]{Hernandez04}, and
\cite[Theorem~3.7]{FHOOPos2022},
the element $\chi_{q,t}(L(M))$ is uniquely determined by
bar invariance and
$\chi_{q,t}(L(M))-E_t(M)\in
\sum_{M'<M}t\mathbb Z[t]E_t(M')$.
Equivalently, we have
\begin{equation}\label{eq:main-opposite-standard}
\begin{aligned}
 &\overline{\chi_{q,t}(L(M))}=\chi_{q,t}(L(M)),\\
 &\chi_{q,t}(L(M))-\mathcal E(M)
 \in\sum_{M'<M}
 t^{-1}\mathbb Z[t^{-1}]\mathcal E(M').
\end{aligned}
\end{equation}
See also~\cite[Theorem~3.3]{Fujita2024}.
For $M\in\mathcal P_\ell^+$, the monomials $M'$
in these expansions also belong to
$\mathcal P_\ell^+$.
Under the identification $v=t$, the latter
characterization uses the same coefficient
ring as the triangular expansion in
\eqref{eq:homogeneous-root-inverse}.

\subsection{Dominant monomials in products of fundamental characters}

By~\cite[Proposition~4.6 and Equation~(4.7)]
{NakajimaTableaux}, the fundamental $(q,t)$-characters
of type $A_{k-1}$ are given by the single-column
tableau formula, with all coefficients equal to one
in the basis of commutative monomials;
see also~\cite[Section~3.2]{Fujita2024}.
In our notation, this gives
\begin{equation}\label{eq:fundamental-path-formula}
 G_{b,i}
 =
 \sum_{1\leq a_1<\cdots<a_i\leq k}
 \prod_{l=1}^i
 \frac{y_{b-i+l,a_l}}
      {y_{b-i+l-1,a_l-1}}.
\end{equation}
The products and quotients in this formula are commutative.
In terms of the variables $Y_{a,p}$, the factor
associated with the entry $a_l$ in position $l$
is
\[
 \frac{Y_{a_l,p+i+a_l-2l}}
 {Y_{a_l-1,p+i+1+a_l-2l}},
 \qquad p=i-2b,
\]
where $Y_{0,p}=Y_{k,p}=1$.
In particular, the column $(1,2,\ldots,i)$
gives the highest monomial $y_{b,i}$.

The following lemma determines all dominant
monomials in a product of two fundamental
characters.

\begin{lemma}\label{lem:fundamental-dominants}
Let $(b,i)<(c,j)$ in the lexicographic order, and put
$d=c-b$.
The product of the ordinary $q$-characters
of the fundamental modules with highest monomials
$y_{b,i}$ and $y_{c,j}$ has highest monomial
$y_{b,i}y_{c,j}$ with multiplicity one.

This product has another dominant monomial
if and only if
$1\leq d\leq\min(k-i,j)$ and $i>j-d$.
In this case, there is exactly one other
dominant monomial, namely
$y_{b,j-d}y_{c,i+d}$,
also with multiplicity one.

More precisely, putting $h_0=j-d$, the latter
monomial is obtained from the column
$(1,\ldots,i)$ in the first factor and the
column
$(1,\ldots,h_0,i+1,i+2,\ldots,i+d)$
in the second factor.
\end{lemma}

\begin{proof}
Consider a monomial
\[
 m=\prod_{l=1}^i
 \frac{y_{b-i+l,a_l}}
      {y_{b-i+l-1,a_l-1}}
\]
associated with a column
$1\leq a_1<\cdots<a_i\leq k$.

For $1\leq l<i$, the factors with first index
$h=b-i+l$ give the quotient
$y_{h,a_l}/y_{h,a_{l+1}-1}$.
These factors cancel if $a_{l+1}=a_l+1$.
Otherwise, the denominator has strictly larger
second index than the numerator.

At the first index $h=b-i$, the contribution is
$y_{b-i,a_1-1}^{-1}$, which is nontrivial precisely
when $a_1>1$.
At the final index $h=b$, the contribution is
$y_{b,a_i}$, which is one if $a_i=k$.

We call a step with $a_{l+1}>a_l+1$, or an initial
entry with $a_1>1$, a jump.
Every jump gives a nontrivial negative factor.
For an internal jump, the second index of the
negative factor is strictly larger than that of
the positive factor with the same first index.

Consider a product of two monomials associated
with two columns.
If both columns have a jump at the same first
index, their negative factors cannot both be
canceled by the positive factors of the other
column: mutual cancellation would give
contradictory strict inequalities between their
second indices.
Consequently, in a dominant product, a jump
in one column can occur only at the final first
index of the other column, and its negative
factor must be canceled by the final positive
factor of that column.

If $b=c$, neither column can have a jump.
Both columns must therefore be the highest
columns, giving $y_{b,i}y_{c,j}$.

Suppose that $b<c$.
The first column cannot have a jump, since
its jumps occur at first indices strictly
less than $b$, whereas cancellation would
require the first index $c$.
Hence the first column is $(1,\ldots,i)$.

The second column can have at most one jump,
which must occur at the first index $b$.
If there is no jump, the second column is
$(1,\ldots,j)$ and the product gives the
highest monomial.

Suppose that there is a jump at the first
index $b$.
The number of entries before this jump is
$h_0=j-d$, and these entries must be
$1,\ldots,h_0$.
In particular, $0\leq h_0<j$.
The negative factor at the jump must cancel
$y_{b,i}$, so the next entry must be $i+1$.
The jump is nontrivial precisely when $i>h_0$.
Since there are no subsequent jumps, the
remaining $d$ entries must be
$i+1,\ldots,i+d$.

Such a column exists if and only if
$1\leq d\leq j$, $i>j-d$, and $i+d\leq k$.
For this column, the product of its monomial
with $y_{b,i}$ is $y_{b,j-d}y_{c,i+d}$.
The columns giving each of the two dominant
monomials are uniquely determined.
This proves the multiplicity assertions,
including the cases $j-d=0$ and $i+d=k$.
\end{proof}

\section{Presentations and the main isomorphism}
\label{sec:defining-relations}
\label{sec:general-fundamental-presentation}

We first compare the normalization factors in the central
quotient of the quantum Grassmannian and the quantum
Grothendieck ring. We then establish presentations of
these algebras by fundamental generators and compare
our presentation of the quantum Grothendieck ring with
existing presentations. We prove Theorem~\ref{thm:intro-A}
and finally illustrate our results in the case of $\Gr(3,6)$.

\subsection{Comparison of normalization factors}
\label{subsec:normalization-comparison}

For the fundamental columns defined
in~\eqref{eq:fundamental-columns}, we use the notation
$J_{b,r}=I_{b,k-r}=[b,b+k]\setminus\{b+r\}$,
$1\leq r<k$.
For the calculation below put $p(b,r)=k-2b-r$ and use $z_{b,r}=M_{J_{b,r}}=y_{b,k-r}$.

As shown in Section~\ref{subsec:source-standard-basis},
the columns $J_{b,r}$, ordered by increasing
$b$ and then decreasing $r$, form a chain
in the componentwise order. This is the lexicographic order on
$(b,i)$ with $i=k-r$.
Ordering these columns by decreasing spectral
parameter corresponds to ordering them by
increasing $2b+r$.

The following proposition compares the
normalization factors for products of fundamental
quantum Pl\"ucker coordinates with those for
products of the corresponding monomials in the
quantum torus $\mathcal Y_t$.

\begin{proposition}\label{prop:general-chain-pair}
Let $I=J_{b,r}$ and $I'=J_{c,s}$ be two
fundamental columns such that $I$ precedes
or equals $I'$ in the lexicographic order.
Then
\begin{equation}\label{eq:general-chain-pair}
 k-|I\cap I'|+\Omega(I,I')
 =J(z_{b,r},z_{c,s}).
\end{equation}

Put $d=c-b$ and define
\[
 \delta_k(x)=
 \begin{cases}
  1,&k\mid x,\\
  0,&k\nmid x.
 \end{cases}
\]
If $d=0$, the common value
in~\eqref{eq:general-chain-pair} is
$1-\delta_{r,s}$.
If $d>0$, it is
\begin{equation}\label{eq:general-chain-deltas}
 -\delta_k(d)
 +\delta_k(d+s)
 +\delta_k(d-r)
 -\delta_k(d+s-r).
\end{equation}

Moreover, if $I\ne I'$ and
$p(b,r)\leq p(c,s)$, then
$P_IP_{I'}=v^{d+1}P_{I'}P_I$ and
$P_I\star P_{I'}
=v^{J(z_{b,r},z_{c,s})}P_{I'}\star P_I$.
In the latter equality, the exponent is $1$
if $d=0$ and $0$ if $d>0$.
In particular, distinct fundamental quantum
Pl\"ucker coordinates associated with the same
spectral parameter commute in $A_v^\Omega$.
\end{proposition}

\begin{proof}
We first compute the left-hand side
of~\eqref{eq:general-chain-pair}.
Define a bilinear form on $\mathbb Z^n$ by
$H(\gamma,\gamma')=\gamma\cdot\gamma'
-\Omega(\gamma,\gamma')$.
By~\eqref{eq:Omega}, for standard basis vectors
$e_x,e_y$, we have
$H(e_x,e_y)=h_k(y-x)$, where
\[
 h_k(z)=
 \begin{cases}
  \delta_k(z),&z\geq0,\\
  -\delta_k(z),&z<0.
 \end{cases}
\]

Put $B_b=\deg P_{C_b}$ and
$D_{b,r}=e_{b+k}-e_{b+r}$.
Then $\deg P_{J_{b,r}}=B_b+D_{b,r}$.
Since each consecutive column $C_b$ contains
exactly one element in each residue class
modulo $k$, for $c\geq b$ we have
$H(B_b,B_c)=k$ and $H(B_b,D_{c,s})=0$.
Moreover,
\[
 H(D_{b,r},B_c)=
 \begin{cases}
  -2,&c=b,\\
  0,&c>b.
 \end{cases}
\]
Indeed, for a standard basis vector $e_x$,
the value of $H(e_x,B_c)$ is $1$ if the
unique element of $C_c$ congruent to $x$
modulo $k$ is greater than or equal to $x$,
and is $-1$ otherwise.

By the definition of $H$, we have
$k-|I\cap I'|+\Omega(I,I')
=k-H(\deg P_I,\deg P_{I'})$.
If $d>0$, the preceding identities give
\[
 k-H(\deg P_I,\deg P_{I'})
 =-H(D_{b,r},D_{c,s})
 =-h_k(d)+h_k(d+s-k)+h_k(d+k-r)-h_k(d+s-r).
\]
Each of the four arguments of $h_k$ is
strictly greater than $-k$.
Consequently, each nonzero value is $1$,
and the above expression is precisely
\eqref{eq:general-chain-deltas}.

If $d=0$, the lexicographic order implies $r\geq s$.
In this case,
$H(D_{b,r},D_{b,s})=1+\delta_{r,s}$.
Therefore
$k-H(\deg P_I,\deg P_{I'})
=2-H(D_{b,r},D_{b,s})=1-\delta_{r,s}$.

We now compute $J(z_{b,r},z_{c,s})$.
Put $L=p(b,r)-p(c,s)=2d+s-r$.
By~\cite[Example~3.7]{Fujita2024},
for $L>0$, the type $A_{k-1}$ commutation
form is given by
\[
 J(z_{b,r},z_{c,s})
 =-\delta_{2k}(L-r+s)-\delta_{2k}(L+r-s)
 +\delta_{2k}(L-r-s)+\delta_{2k}(L+r+s).
\]
Substituting $L=2d+s-r$ gives
\eqref{eq:general-chain-deltas}.

Suppose that $d>0$ and $L\leq0$.
Then $r-s\geq2d$.
The four arguments $d$, $d+s$, $d-r$ and
$d+s-r$ in~\eqref{eq:general-chain-deltas}
are nonzero and have absolute value strictly
less than $k$.
Hence~\eqref{eq:general-chain-deltas} is zero.

If $L=0$, we have $J(z_{b,r},z_{c,s})=0$.
If $L<0$, we apply the above formula after
interchanging the two arguments of $J$.
Since the displayed expression in terms of
$\delta_{2k}$ is symmetric under $L\mapsto-L$,
it follows that $J(z_{b,r},z_{c,s})=0$
whenever $d>0$ and $L\leq0$.

Finally, if $d=0$ and $r>s$, then
$L=s-r<0$.
Applying the formula with $-L=r-s$ and using
skew-symmetry gives $J(z_{b,r},z_{b,s})=1$.
If $r=s$, the value is zero.
This proves~\eqref{eq:general-chain-pair}
in all cases.

It remains to prove the quasi-commutation
assertions.
Suppose that $I\ne I'$ and
$p(b,r)\leq p(c,s)$.

If $d=0$, then $r>s$, and
$I\setminus I'=\{b+s\}$,
$I'\setminus I=\{b+r\}$.
If $d>0$, then $r-s\geq2d$, and
$I\setminus I'=[b,c-1]\cup\{c+s\}$ and $I'\setminus I=\{b+r\}\cup[b+k+1,c+k]$.
In both cases, every element of
$I\setminus I'$ is smaller than every
element of $I'\setminus I$.
By the quasi-commutation formula for
quantum minors~\cite{LZ98}, see also
\cite[Theorem~6.2 and Corollary~6.5]{Goodearl2006},
we have
$P_IP_{I'}=v^{|I\setminus I'|}P_{I'}P_I
=v^{d+1}P_{I'}P_I$.

By the definition of the twisted multiplication,
$P_I\star P_{I'}
=v^{d+1+\Omega(I,I')}P_{I'}\star P_I$.
Since $|I\setminus I'|=k-|I\cap I'|=d+1$,
Equation~\eqref{eq:general-chain-pair} gives
$P_I\star P_{I'}
=v^{J(z_{b,r},z_{c,s})}P_{I'}\star P_I$.
The preceding calculation of
$J(z_{b,r},z_{c,s})$ proves the remaining
assertions.
\end{proof}

For a multiplicity vector $\mathbf u$, let $y_1,\ldots,y_h$ be the
factors of $M_{\mathbf u}$ in decreasing spectral order, with ties
resolved by the lexicographic order. Denote the corresponding normalized
fundamental generators by $X_1,\ldots,X_h$.

By~\eqref{eq:source-chain-standard} and Proposition~\ref{prop:general-chain-pair},
we obtain
\begin{equation}\label{eq:general-spectral-standard}
 S_{\mathbf u}=v^{-\frac12\sum_{a<b}J(y_a,y_b)}X_1\cdots X_h.
\end{equation}
Indeed, to pass from the lexicographic order to the decreasing spectral
order, interchange adjacent factors with spectral parameters
$p_a<p_b$. By Proposition~\ref{prop:general-chain-pair} and the
skew-symmetry of $J$, each such interchange preserves the normalized
product by the quasi-commutation formula for the corresponding
fundamental generators.
All other pair contributions to the exponent are unchanged.
Factors with equal spectral parameters commute.

\subsection{Presentation of the central quotient}

We work over $\bbk=\mathbb Q(v^{1/2})$ and identify
the parameters by $v=t$.
Let $\mathcal D:A_v^\Omega\to B_\chi$ be the
algebra homomorphism defined in
Section~\ref{subsec:central-specialization}.
We use the generators $X_{b,i}$ and monomials $y_{b,i}$ defined in
Section~\ref{subsec:source-standard-basis}.
We extend this notation to the boundary indices
by setting
$I_{b,0}=C_b, I_{b,k}=C_{b+1}$,
and $X_{b,0}=X_{b,k}=1$. Recall that $y_{b,0}=y_{b,k}=1$.

Using the skew-symmetric form $J$ defined in
Section~\ref{subsec:general-spectral-words}, put
$C_{b,i;c,j} = v^{-J(y_{b,i},y_{c,j})/2}X_{b,i}X_{c,j}$,
and
$C^{\rm op}_{b,i;c,j} = v^{J(y_{b,i},y_{c,j})/2}X_{c,j}X_{b,i}$.
These definitions also apply when one or both
indices belong to $\{0,k\}$, with
$J(1,M)=J(M,1)=0$.

The following proposition gives a presentation of
$B_\chi$ using only the fundamental quantum
Pl\"ucker coordinates.

\begin{proposition}\label{thm:general-defining}
The algebra $B_\chi$ is generated by the elements
$X_{b,i}$,
$1\leq b\leq m$, $1\leq i<k$,
subject to the following relations.

For every pair $(b,i)<(c,j)$ in the lexicographic order,
put $d=c-b$. Then
\begin{equation}\label{eq:uniform-pair-relation}
 C^{\rm op}_{b,i;c,j}-C_{b,i;c,j}
 =
 \begin{cases}
 (v-v^{-1})C_{b,j-d;c,i+d},
 &
 1\leq d\leq\min(k-i,j),\quad i>j-d,\\[2mm]
 0,&\text{otherwise}.
 \end{cases}
\end{equation}
Here the boundary generators are interpreted as
$X_{b,0}=X_{b,k}=1$.

Whenever the correction term is nonzero, the
element $C_{b,j-d;c,i+d}$ is invariant under
the modified bar involution of $B_\chi$.

Moreover, the ordered monomials in the generators
$X_{b,i}$, with the factors arranged in the chain
order, form a $\bbk$-basis of $B_\chi$.
\end{proposition}

\begin{proof}
We first establish the corresponding relations
between the quantum minors in $A_v$.

Put
$r=k-i$, $s=k-j$, $I=I_{b,i}$, and $H=I_{c,j}$,
and let
$u=b+r, w=c+s$
be the omitted entries of the two fundamental
columns, respectively.
Thus
$I=[b,b+k]\setminus\{u\}$, $H=[c,c+k]\setminus\{w\}$.
Put $h=k-|I\cap H|$.

We claim that every element of $I\setminus H$
is smaller than every element of $H\setminus I$,
except when
\begin{equation}\label{eq:crossing-holes}
 b<c,\qquad c\leq u<w\leq b+k.
\end{equation}
Indeed, if $b=c$, the lexicographic order implies $i<j$,
and hence $w<u$.
The two set differences are then
$I\setminus H=\{w\},  H\setminus I=\{u\}$.
If $b<c$, the elements of $I\setminus H$
outside the common interval lie to the left
of that interval, whereas the elements of
$H\setminus I$ outside the common interval
lie to its right.
Inside the common interval, the only possible
element of $I\setminus H$ is $w$, and the only
possible element of $H\setminus I$ is $u$.
Thus the two set differences fail to be separated
precisely when~\eqref{eq:crossing-holes} holds.

Substituting
$u=b+k-i, w=c+k-j$
into~\eqref{eq:crossing-holes}, we obtain
$1\leq d\leq\min(k-i,j)$, $i>j-d$.
These are precisely the conditions in the first
case of~\eqref{eq:uniform-pair-relation}.

We now apply the commutation formula for quantum
minors in~\cite[Corollary~5.3,
Equations~(5.8) and~(5.9)]{Goodearl2006}.
Since the two maximal minors have the same row
index set $[1,k]$, the formula takes the form
\begin{equation}\label{eq:main-Goodearl-specialized}
 v^hP_HP_I-P_IP_H
 =
 \sum_K\lambda_KP_KP_{K'},
\end{equation}
where the sum is over the increasing $k$-subsets
$K$ satisfying
$I\cap H\subseteq K\subseteq I\cup H,  K<I$
in the componentwise order, and
$K'=(I\cap H)\cup((I\cup H)\setminus K)$.
The coefficient is given by
\[
 \lambda_K
 =
 (v-v^{-1})^{|I\setminus K|}
 (-v)^{\kappa(K,I,H)}
 \xi_v(I\setminus K;K\setminus I),
\]
where
$\kappa(K,I,H) = \operatorname{inv}(D,I\setminus K) - \operatorname{inv}(D,K\setminus I)$,
and
$D=(K\setminus K')\cup(H\setminus I)$.
The function $\xi_v$ is defined in
\cite[Section~4.4, Equation~(4.12)]{Goodearl2006}.
We will only need the equality
$\xi_v(\{w\};\{u\})=1$ if $u<w$.

If every element of $I\setminus H$ is smaller
than every element of $H\setminus I$, there is
no subset $K$ satisfying the above conditions.
Indeed, replacing any element of $I\setminus H$
by an element of $H\setminus I$ cannot produce
a componentwise smaller subset.
Consequently,
$v^hP_HP_I=P_IP_H$.

Suppose now that~\eqref{eq:crossing-holes} holds.
Then
$I\setminus H=[b,c-1]\cup\{w\}$ and $H\setminus I=\{u\}\cup[b+k+1,c+k]$.
Every subset $K$ occurring in
\eqref{eq:main-Goodearl-specialized} must contain
all the elements of $[b,c-1]$, since they are
the smallest available entries.
The common entries $I\cap H$ must also belong
to $K$.
There is therefore precisely one possible
replacement that gives a componentwise smaller
subset, namely replacing $w$ by $u$.

It follows that the sum in
\eqref{eq:main-Goodearl-specialized} has exactly
one term, with
$K=[b,b+k]\setminus\{w\}$, $K'=[c,c+k]\setminus\{u\}$.
For these two subsets, we have
$I\setminus K=\{w\},  K\setminus I=\{u\}$,
and
$D=[b,c-1]\cup\{u\}\cup[b+k+1,c+k]$.
The numbers of elements of $D$ strictly greater
than $u$ and strictly greater than $w$ are
both equal to $d$.
Hence
$\kappa(K,I,H)=0$.
Since $\xi_v(\{w\};\{u\})=1$, the coefficient
$\lambda_K$ is $v-v^{-1}$.
We obtain
\begin{equation}\label{eq:uniform-FRT-pair}
 v^hP_HP_I-P_IP_H
 =
 (v-v^{-1})P_KP_{K'}.
\end{equation}

In terms of the fundamental columns, the two
new subsets are
$K=I_{b,j-d},  K'=I_{c,i+d}$.
The inequalities in~\eqref{eq:crossing-holes}
imply that
$0\leq j-d<i<i+d\leq k$.
If $j-d=0$, then $K=C_b$.
If $i+d=k$, then $K'=C_{c+1}$.
Moreover,
$K\cap K'=I\cap H$,
and every element of $K\setminus K'$ is smaller
than every element of $K'\setminus K$.
In particular,
$P_KP_{K'}=v^hP_{K'}P_K$.

We next pass to the twisted algebra and its
central quotient.
Put
$J_0=J(y_{b,i},y_{c,j})$.
By Proposition~\ref{prop:general-chain-pair},
$J_0=h+\Omega(I,H)$.
Using the definition of the twisted multiplication
and the multiplicativity of $\mathcal D$ with
respect to this multiplication, we obtain
\[
\begin{aligned}
 C_{b,i;c,j}&=v^{-J_0/2}\mathcal D(P_I\star P_H)=v^{-h/2}\mathcal D(P_IP_H),\\
 C^{\rm op}_{b,i;c,j}&=v^{J_0/2}\mathcal D(P_H\star P_I)=v^{h/2}\mathcal D(P_HP_I).
\end{aligned}
\]
Here $\mathcal D$ is applied to ordinary products
by viewing them as elements of the underlying
vector space of $A_v^\Omega$.

If both correction indices $j-d$ and $i+d$
belong to $[1,k-1]$, then
Proposition~\ref{prop:general-chain-pair} gives
$J(y_{b,j-d},y_{c,i+d}) = h+\Omega(K,K')$.
If at least one of the correction indices is
a boundary index, the corresponding consecutive
minor is central in $A_v^\Omega$.
Since
$P_KP_{K'}=v^hP_{K'}P_K$,
centrality implies
$h+\Omega(K,K')=0$.
This agrees with our convention that $J$ is zero
when one of its arguments is the unit monomial.

In either case, we obtain
$v^{-h/2}\mathcal D(P_KP_{K'}) = C_{b,j-d;c,i+d}$.
Multiplying~\eqref{eq:uniform-FRT-pair}
by $v^{-h/2}$ and applying $\mathcal D$ gives
the first case of~\eqref{eq:uniform-pair-relation}.
The second case follows in the same way from
$v^hP_HP_I=P_IP_H$.

Since $P_KP_{K'}=v^hP_{K'}P_K$, we have
$\overline{v^{-h/2}P_KP_{K'}} = v^{h/2}P_{K'}P_K = v^{-h/2}P_KP_{K'}$.
The compatibility of $\mathcal D$ with the bar
involutions therefore implies that
$C_{b,j-d;c,i+d}$ is bar invariant.

It remains to show that the relations
\eqref{eq:uniform-pair-relation} generate the
defining ideal.

Let $\mathcal F$ be the free associative
$\bbk$-algebra on generators $z_{b,i}$,
$1\leq b\leq m$, $1\leq i<k$.
Let $\mathcal I$ be the two-sided ideal generated
by the relations~\eqref{eq:uniform-pair-relation}
written in these generators, with
$z_{b,0}=z_{b,k}=1$.
Put
$\mathcal P=\mathcal F/\mathcal I$.
The preceding argument shows that the assignment
$z_{b,i}\mapsto X_{b,i}$ defines an algebra
homomorphism
$\rho:\mathcal P\to B_\chi$.
By Proposition~\ref{thm:general-source-pbw},
the elements $X_{b,i}$ generate $B_\chi$.
Thus $\rho$ is surjective.

Assign the weight
$f(i)=i(k-i)$
to the generator $z_{b,i}$ and extend this
weight additively to words.
For the boundary indices, put $f(0)=f(k)=0$.

Each relation in~\eqref{eq:uniform-pair-relation}
expresses a product of two generators in reversed
lexicographic order as a unit multiple of the
product in lexicographic order, possibly together with
a correction term.
In the latter case,
$f(i)+f(j)-f(j-d)-f(i+d) =2d(i-j+d)>0$.
Thus the correction term has strictly smaller
weight than the original product.
This remains true when one or both of the
correction generators are units.

The ordered product has the same weight as
the original product, but replacing an adjacent
reversed pair by an ordered pair decreases
the number of inversions of the word by one.
It follows by induction on the pair
$(\text{weight},\text{number of inversions})$,
ordered lexicographically in $\mathbb N^2$,
that every word in the generators of
$\mathcal P$ is a linear combination of words
in lexicographic order. The same induction works over
$\mathbb Z[v^{\pm1/2}]$, since the leading coefficients
used in reordering are powers of $v^{1/2}$.

For every multiplicity vector $\mathbf u$,
let $z^{\mathbf u}$ be the corresponding
chain-ordered word.
By~\eqref{eq:source-chain-standard},
its image under $\rho$ is a nonzero power of
$v^{1/2}$ times $S_{\mathbf u}$.
The elements $S_{\mathbf u}$ form the standard basis
$\mathcal S_\chi$ of $B_\chi$ by Proposition~\ref{thm:general-source-pbw}.
Consequently, the images of the ordered words
$z^{\mathbf u}$ are linearly independent.

Since the ordered words span $\mathcal P$,
this proves that $\rho$ is injective.
Hence $\rho$ is an isomorphism, and the ordered
monomials form a basis of the presented algebra.
\end{proof}

We may also express the defining relations
directly in terms of the images of the original
quantum Pl\"ucker coordinates, without the
rescaling by $\mathcal D$.

For $0\leq i\leq k$, put
$Z_{b,i}=\pi_\chi(P_{I_{b,i}})$.
In particular,
$Z_{b,0}=t^{\chi(C_b)/2}$, $Z_{b,k}=t^{\chi(C_{b+1})/2}$.
The relation between these generators and the
normalized generators is
$X_{b,i} = t^{-\chi(I_{b,i})/2}Z_{b,i}$, $0\leq i\leq k$.

The same relations also give a presentation in the generators
$Z_{b,i}$, $1\leq b\leq m$, $1\leq i<k$, with boundary values
$Z_{b,0}=t^{\chi(C_b)/2}$ and $Z_{b,k}=t^{\chi(C_{b+1})/2}$,
and with $C_{b,i;c,j}$ and $C^{\rm op}_{b,i;c,j}$ defined using
$Z$ in place of $X$. Indeed, whenever the correction term is nonzero,
$\deg P_{I_{b,i}}+\deg P_{I_{c,j}}=\deg P_{I_{b,j-d}}+\deg P_{I_{c,i+d}}$,
so the invertible change of generators preserves
\eqref{eq:uniform-pair-relation} by the additivity of $\chi$.

\subsection{Presentation of the quantum Grothendieck ring}

We now establish the presentation of $K_t(\mathcal C_\ell)$
independently of the Grassmannian. The relations follow
from the segment exact sequences of Kang, Kashiwara and Kim
\cite[Proposition~4.2.3]{KKK2018} and their realization of
the quantum Grothendieck ring as recalled in
\cite[Section~5 and Theorem~5.1]{Fujita2024}.
We give the change of normalization explicitly.

We use the segment convention of~\cite{CDFL20}. For integers
$a\leq b$, put
\begin{equation}\label{eq:segment-convention}
 \Delta=[a,b],\qquad
 |\Delta|=b-a+1,\qquad
 \beta_\Delta=\epsilon_{a-1}-\epsilon_b,\qquad
 Y_\Delta=Y_{b-a+1,a+b-1},
\end{equation}
where the last expression is defined for $1\leq |\Delta|<k$.
We use $[a,a-1]=\varnothing$ for the empty segment, with
$|\varnothing|=0$ and $\beta_\varnothing=0$.
The root lattice of type $A_\infty$ is realized in the lattice
with orthonormal basis $\epsilon_r$ ($r\in\mathbb Z$), with simple
roots $\epsilon_r-\epsilon_{r+1}$. The segment module $L_\Delta$
is the one-dimensional graded module of degree zero for the
quiver Hecke algebra of weight $\beta_\Delta$: the idempotent
of the word $(a-1,a,\ldots,b-1)$ acts as the identity, the other
idempotents act as zero, and all polynomial and crossing generators
act as zero. This is the module $L(a-1,b-1)$ of~\cite[Section~4.2]{KKK2018}. Thus the shift by one in the
segment labels accounts for the difference between our
convention and that of~\cite{KKK2018,Fujita2024}.
More precisely, the root $\alpha_{u,w}=\epsilon_u-\epsilon_w$
in~\cite[Equation~(5.5)]{Fujita2024} corresponds to our
segment $[u+1,w]$. We have
$Y_\Delta=Y_{w-u,u+w}$ for $u=a-1$ and $w=b$.
In particular, the segment corresponding to $y_{b,i}$ is
\begin{equation}\label{eq:fundamental-segment}
 \Delta_{b,i}=[1-b,i-b],\qquad
 \beta_{\Delta_{b,i}}=\epsilon_{-b}-\epsilon_{i-b},\qquad
 Y_{\Delta_{b,i}}=y_{b,i}=Y_{i,i-2b}.
\end{equation}
We extend the interval formula to the boundary indices by
$\Delta_{b,0}=\varnothing$ and $\Delta_{b,k}=[1-b,k-b]$.
The grading parameter of~\cite{KKK2018} is $t^{-1}$
in our convention, as in~\cite[Section~5]{Fujita2024}.

For the fundamental $(q,t)$-characters, define
\[
 C^G_{b,i;c,j}
 =t^{-J(y_{b,i},y_{c,j})/2}G_{b,i}*G_{c,j},
 \qquad
 C^{G,{\rm op}}_{b,i;c,j}
 =t^{J(y_{b,i},y_{c,j})/2}G_{c,j}*G_{b,i}.
\]

\begin{proposition}\label{thm:general-target-relations}
Over $\mathbb Z[t^{\pm1/2}]$, the algebra $K_t(\mathcal C_\ell)$
is generated by the fundamental $(q,t)$-characters $G_{b,i}$,
$1\leq b\leq m$, $1\leq i<k$, subject to the following relations.
For every pair $(b,i)<(c,j)$ in the lexicographic order,
put $d=c-b$. Then
\begin{equation}\label{eq:target-pair-relation}
 C^{G,{\rm op}}_{b,i;c,j}-C^G_{b,i;c,j}
 =
 \begin{cases}
 (t-t^{-1})C^G_{b,j-d;c,i+d},
 &
 1\leq d\leq\min(k-i,j),\quad i>j-d,\\[2mm]
 0,&\text{otherwise}.
 \end{cases}
\end{equation}
Here $G_{b,0}=G_{b,k}=1$.
The chain-ordered monomials in the $G_{b,i}$ form a
$\mathbb Z[t^{\pm1/2}]$-basis. In particular, after extension
of scalars, this gives the same presentation as
Proposition~\ref{thm:general-defining} under $v=t$.
\end{proposition}

\begin{proof}
For a segment $\Delta=[a,b]$ as in
\eqref{eq:segment-convention}, let $H_\Delta$ denote the
image of the class of $L_\Delta$ under the localization
and the isomorphism of~\cite[Theorem~5.1]{Fujita2024}.
By the localization and equation (5.5) of that paper,
with the root $\alpha_{a-1,b}$, we have
\[
 H_\Delta=
 \begin{cases}
 1,& |\Delta|=0,\\
 t^{-1/2}\chi_{q,t}(L(Y_\Delta)),&1\leq |\Delta|<k,\\
 1,& |\Delta|=k,\\
 0,& |\Delta|>k.
 \end{cases}
\]
In particular, $G_{b,i}=t^{1/2}H_{\Delta_{b,i}}$ for
$1\leq i<k$. For a segment $S$, write $\beta_S$ for its
weight and $h_S$ for its class in the graded Grothendieck ring of
quiver Hecke modules, with $h_\varnothing=1$. Multiplication in that
ring is induced by convolution (induction from the tensor product
of the two quiver Hecke algebras). In the quotient and localization
used above, segments of length greater than $k$ vanish and segments
of length $k$ are identified with the unit.
The bilinear form used in~\cite{Fujita2024} is
$B_k(\beta,\gamma)=\sum_{r>0}(S_k^r\beta,\gamma)$, where $S_k(\epsilon_a)=\epsilon_{a+k}$.
The image of the untwisted product $h_Sh_T$ is therefore
$t^{-B_k(\beta_S,\beta_T)}H_S*H_T$,
by~\cite[Equation~(5.3)]{Fujita2024}.
The skew form of the quantum torus in
\cite[Section~3.1]{Fujita2024} is precisely our $J$.

Put $R=\Delta_{b,i}=[1-b,i-b]$ and
$L=\Delta_{c,j}=[1-c,j-c]$.
Two nonempty segments are linked if their union is a segment and neither
contains the other; adjacent disjoint segments are included.
By~\cite[Proposition~4.2.3(iii)--(vi)]{KKK2018},
if the segments are not linked, their convolution
classes commute up to a power of the grading parameter.
After twisting and localization this gives
quasi-commutation of $G_{b,i}$ and $G_{c,j}$.
The exponent is necessarily $J(y_{b,i},y_{c,j})$,
by comparing their highest monomial, whose coefficient
in the two normalized products is one.
This proves the relation in the unlinked case.

The linked case occurs precisely when $R=[u,f]$,
$L=[a,e]$ with $a<u\leq e+1\leq f$.
Put $U=L\cup R=[a,f]$ and $V=L\cap R=[u,e]$,
where $V=\varnothing$ if $u=e+1$, and set
$\varepsilon=\delta_{u,e+1}$.
The two exact sequences in
\cite[Proposition~4.2.3(v), (vi)]{KKK2018} give
\[
 h_Rh_L-t^{\varepsilon}h_Lh_R
 =(t^{-1}-t)h_Uh_V,
\]
where $h_V=1$ if $u=e+1$.
The two lengths of $R,L$ are less than $k$.
In terms of the original indices,
$U=[1-c,i-b]=\Delta_{c,i+d}$ and
$V=[1-b,j-c]=\Delta_{b,j-d}$,
with lengths $i+d$ and $j-d$, respectively.
In particular, linkedness is equivalent to
$1\leq d\leq j$ and $i>j-d$.

For brevity put $B(S,T)=B_k(\beta_S,\beta_T)$.
In the linked case, evaluating the four endpoint
pairings in the definition of $B_k$ gives
$B(R,L)=0, B(L,R)=-\delta_{f-a+1,k}$.
Indeed, the two weights are $\epsilon_{a-1}-\epsilon_e$
and $\epsilon_{u-1}-\epsilon_f$. The only possible equality
between a positively shifted endpoint of the weight of $L$
and an endpoint of the weight of $R$ is $a-1+k=f$,
contributing with sign $-1$.
The length of $U$ is less than $2k$.
If $f-a+1>k$, the correction term is zero after
localization, and comparison of highest monomials gives
$C^{G,{\rm op}}_{b,i;c,j}=C^G_{b,i;c,j}$.

Suppose now that $f-a+1\leq k$. Set
$\delta=\delta_{f-a+1,k}$ and $\lambda=\varepsilon+\delta$.
The endpoint pairings also give $B(U,V)=B(V,U)=0$.
The number of nonempty correction segments of length
strictly less than $k$ is $2-\varepsilon-\delta$.
Thus, substituting the images $H_S$ into the segment
relation and multiplying by $t$ gives
\[
 G_{b,i}*G_{c,j}-t^\lambda G_{c,j}*G_{b,i}
 =(t^{-1}-t)t^{\lambda/2}G_{c,i+d}*G_{b,j-d}.
\]
The correction segments are strictly nested if both are
internal. They commute, by the same segment relation
for nested segments and the vanishing of the two twist
exponents; the assertion also holds if a segment is a unit.
Comparison of highest monomials gives
$J(y_{b,j-d},y_{c,i+d})=0$.
The correction highest monomial is strictly below
$y_{b,i}y_{c,j}$ in the dominance order (it is the
additional dominant monomial of Lemma~\ref{lem:fundamental-dominants}).
Hence comparison of the highest monomial in the displayed
identity gives $J(y_{b,i},y_{c,j})=\lambda$.
Dividing by $t^{\lambda/2}$ and reversing the difference
now proves~\eqref{eq:target-pair-relation} in our normalization.

It remains to show that these are defining relations.
The straightening argument in the proof of
Proposition~\ref{thm:general-defining} applies over
$\mathbb Z[t^{\pm1/2}]$ and shows that the chain-ordered monomials
span the algebra defined by~\eqref{eq:target-pair-relation}.
If $(b,i)<(c,j)$ and $i-2b<j-2c$, the correction term in
\eqref{eq:target-pair-relation} is zero: either $d=0$ or
$i<j-2d\leq j-d$.
Thus the corresponding adjacent interchanges preserve the
normalized product when passing to decreasing spectral order.
By~\eqref{eq:kl-standard-product}, the images of the chain-ordered
monomials are therefore, up to powers of $t^{1/2}$, precisely the
standard basis elements $\mathcal E(M)$, $M\in\mathcal P_\ell^+$.
Since these form a $\mathbb Z[t^{\pm1/2}]$-basis of
$K_t(\mathcal C_\ell)$, the natural homomorphism from the
presented algebra to $K_t(\mathcal C_\ell)$ is an isomorphism.
\end{proof}

For $k\geq4$, the correction term in
\eqref{eq:target-pair-relation} can be quadratic.
For instance, taking $k=4$, $i=j=2$ and $c=b+1$ gives
\[
 G_{b+1,2}*G_{b,2}-G_{b,2}*G_{b+1,2}
 =(t-t^{-1})G_{b,1}*G_{b+1,3}.
\]

\subsection{Comparison with existing presentations}
Hernandez and Leclerc~\cite[Theorem~7.3]{HL2015}
gave a presentation of the quantum Grothendieck ring of
$\mathcal C_{\mathbb Z}$ over $\mathbb C(t^{1/2})$
in simply-laced types.
Fujita, Hernandez, Oh, and Oya~\cite[Theorem~10.4]{FHOOPos2022}
extended this presentation to non-simply-laced types using
$Q$-data. Jang, Lee, and Oh~\cite[Theorem~7.2]{JLO2023Braid}
obtained a corresponding presentation for quantum virtual
Grothendieck rings. For type $A$, see also
\cite[Definition~3.9 and Theorem~3.10]{Fujita2024}.
We use the corresponding presentation over $\mathbb Q(t^{1/2})$;
see also~\cite[Corollary~8.3]{HL2015}.
Its generators $e_{j,h}$ are fundamental $(q,t)$-characters
corresponding to the pairs
$(\alpha_j,h)$ under the Auslander--Reiten bijection
$\phi_\xi:\{(i,p):i\in\mathsf I,\ p\in\xi_i+2\mathbb Z\}
\to\Delta_+\times\mathbb Z$, where $\Delta_+$ is the set of
positive roots of type $A_{k-1}$ and $\alpha_j$ is its $j$th simple root. Suppressing the
symbol $*$ for multiplication, the relations are
\begin{equation}\label{eq:hl-bosonic-presentation}
\begin{aligned}
 e_{i,h}e_{j,h}&=e_{j,h}e_{i,h} && (|i-j|>1),\\
 e_{i,h}^2e_{j,h}-(t+t^{-1})e_{i,h}e_{j,h}e_{i,h}
       +e_{j,h}e_{i,h}^2&=0 && (|i-j|=1),\\
 e_{i,h}e_{j,h+1}&=t^{-(\alpha_i,\alpha_j)}e_{j,h+1}e_{i,h}
                  +\delta_{ij}(1-t^{-2}),\\
 e_{i,h}e_{j,p}&=t^{(-1)^{p-h}(\alpha_i,\alpha_j)}e_{j,p}e_{i,h}
                  &&(p>h+1).
\end{aligned}
\end{equation}
Here $(\alpha_i,\alpha_j)$ is the symmetric Cartan
pairing of type $A_{k-1}$. These are the quantum Serre relations,
quantum boson relations, and quasi-commutation relations for
integer indices differing by at least two.
Write $\widehat{\mathcal A}_{k,[h_1,h_2]}$ for the subalgebra
of this bosonic extension generated by the indicated
simple-root generators with $h_1\leq h\leq h_2$.
It has the presentation obtained by restricting
\eqref{eq:hl-bosonic-presentation} to $h_1\leq h,p\leq h_2$:
the relations for distinct integer indices allow products to be ordered
by these indices, and for each fixed $h$ the quantum Serre relations
give a Poincar\'e--Birkhoff--Witt (PBW) basis.
The resulting ordered products are independent in the
ambient algebra by the proof of~\cite[Theorem~7.3]{HL2015}.
For each fixed $h$, the fundamental characters corresponding to
$\Delta_+\times\{h\}$ belong to the subalgebra generated by the
$e_{j,h}$; see~\cite[Section~7.1]{HL2015}.

For our increasing height function, the explicit
bijection is given in~\cite[Example~3.8, equations (3.3)--(3.4)]
{Fujita2024}. Write $1-b=kq_0+r$ with $1\leq r\leq k$
and let $r'$ be the representative of $r+i$ modulo $k$
in $[1,k]$. With $\alpha_{a,e}=\epsilon_a-\epsilon_e$
for $1\leq a<e\leq k$, one has
\[
 \phi_\xi(i,i-2b)=
 \begin{cases}
 (\alpha_{r,r'},2q_0+1),&r<r',\\
 (\alpha_{r',r},2q_0+2),&r>r'.
 \end{cases}
\]
In particular, every block $ak+1\leq b\leq(a+1)k$
($a\geq0$) corresponds to $\Delta_+\times\{-2a-1,-2a\}$.
The corresponding simple-root characters are
\begin{equation}\label{eq:hl-window-generators}
 e_{j,-2a-1}\longmapsto G_{ak+k+1-j,1},\qquad
 e_{j,-2a}\longmapsto G_{ak+k-j,k-1}
 \qquad(1\leq j<k).
\end{equation}
Consequently, if $\ell+1=r_0k$ for a positive integer $r_0$,
Proposition~\ref{thm:general-target-relations} and the
generation result for each fixed $h$ in~\cite[Section~7.1]{HL2015}
give an identification
\begin{equation}\label{eq:complete-layer-window}
 \widehat{\mathcal A}_{k,[-2r_0+1,0]}
 \xrightarrow{\sim}
 \mathbb Q(t^{1/2})\otimes_{\mathbb Z[t^{\pm1/2}]}
 K_t(\mathcal C_{r_0k-1}).
\end{equation}
Indeed, the spectral indices defining $\mathcal C_{r_0k-1}$
correspond to all pairs in $\Delta_+\times[-2r_0+1,0]$, so the two
subalgebras of the ambient ring coincide. If the image of the
spectral indices defining $\mathcal C_\ell$ under $\phi_\xi$ contains
only a proper nonempty subset of $\Delta_+\times\{h\}$ for some $h$,
restricting the integer indices in~\eqref{eq:hl-bosonic-presentation}
does not directly give a presentation of the corresponding subalgebra.

More generally, Oh and Park~\cite[Theorem~7.11]{OP2025}
established PBW bases and Levendorskii--Soibelman relations
for subalgebras of bosonic extensions associated with positive
braid monoid elements. Together, these give presentations
by PBW root vectors over the rational function field.
See also~\cite[Section~4]{KKOPBraid2024} for the PBW and global
bases of these subalgebras.
The presentation in Proposition~\ref{thm:general-target-relations}
uses all fundamental $(q,t)$-characters in $\mathcal C_\ell$
and gives explicit quadratic relations over
$\mathbb Z[t^{\pm1/2}]$ for every $\ell\geq0$.

\subsection{The main isomorphism and canonical bases}

We now combine the presentation of the central
quotient with the relations among the fundamental
$(q,t)$-characters and the comparison of bases
in Section~\ref{subsec:normalization-comparison}.

\begin{theorem}\label{thm:general-based-isomorphism}
The elements $F_a$, $1\leq a\leq s$, are central in
$A_v^\Omega$. For $\chi(e_a)=\delta_{a,n}-\delta_{a,1}$ and
$B_\chi$ as defined in Section~\ref{subsec:central-specialization}, there is a unique
$\bbk$-algebra isomorphism $\Phi_\chi:B_\chi\xrightarrow{\sim}K$
with
$\Phi_\chi(\pi_\chi(P_{I_{b,i}}))=
t^{\chi(I_{b,i})/2}\chi_{q,t}(L(Y_{i,i-2b}))$.
The images of the fundamental Pl\"ucker coordinates generate
$B_\chi$. For every tableau $T$,
$\Phi_\chi(\pi_\chi(b_T))=t^{\chi(T)/2}\chi_{q,t}(L(M_T))$;
in particular, the images of all quantum Pl\"ucker coordinates
are given by~\eqref{eq:scaled-plucker-images}.
The elements $\pi_\chi(b_T)$ for reduced $T$ form a basis
of $B_\chi$, and $\Phi_\chi$ is compatible with the modified
bar involution of Section~\ref{sec:twist}.
More precisely, for every $\mathbf u$ and every reduced tableau $T$,
\[
 \Phi_\chi(S_{\mathbf u})=\mathcal E(M_{\mathbf u}),\qquad
 \Phi_\chi(\widehat b_T)=\chi_{q,t}(L(M_T)).
\]
Thus $\Phi_\chi$ identifies the standard bases $\mathcal S_\chi$ and
$\mathcal S_K$, and the dual canonical basis $\mathcal B_\chi$ with
the canonical basis $\mathcal B_K$.
\end{theorem}

\begin{proof}
The centrality, reduced basis, and modified bar involution
are established in Proposition~\ref{prop:center},
Proposition~\ref{thm:free}, and Section~\ref{sec:twist}, respectively.
Propositions~\ref{thm:general-defining}
and~\ref{thm:general-target-relations} give identical
presentations after extension of scalars and $v=t$. The assignment
$X_{b,i}\mapsto \chi_{q,t}(L(Y_{i,i-2b}))$ defines a
$\bbk$-algebra homomorphism $\Phi_\chi:B_\chi\to K$.
The identical presentations make $\Phi_\chi$ an isomorphism.
By \eqref{eq:general-spectral-standard} and
\eqref{eq:kl-standard-product},
$\Phi_\chi(S_{\mathbf u})=\mathcal E(M_{\mathbf u})$.
The modified bar on $B_\chi$ and the quantum torus bar on $K$ reverse
multiplication and fix the corresponding normalized generators. Hence
$\Phi_\chi(\overline x^{\,\chi})=\overline{\Phi_\chi(x)}$ for $x\in B_\chi$.
For each reduced tableau $T$, applying $\Phi_\chi$ to
\eqref{eq:source-standard-triangularity} gives
\[
 \Phi_\chi(\widehat b_T)-\mathcal E(M_T)
 \in\sum_{M'<M_T}t^{-1}\mathbb Z[t^{-1}]\mathcal E(M').
\]
By bar invariance and the uniqueness in~\eqref{eq:main-opposite-standard},
$\Phi_\chi(\widehat b_T)=\chi_{q,t}(L(M_T))$.

For an arbitrary tableau $T\in\SSYT(k,[n])$,
we have $M_T=M_{T_{\mathrm{red}}}$ and,
by~\eqref{eq:normalized-reduction},
$\mathcal D(b_T)=\widehat b_{T_{\mathrm{red}}}$.
Consequently,
\[
 \Phi_\chi(\pi_\chi(b_T))
 =t^{\chi(T)/2}\Phi_\chi(\mathcal D(b_T))
 =t^{\chi(T)/2}\chi_{q,t}(L(M_T)),
\]
which proves~\eqref{eq:scaled-tableau-images}.
Taking $T=(I)$ and using $b_{(I)}=P_I$, we obtain
\eqref{eq:scaled-plucker-images}.
The uniqueness of the isomorphism follows from
the generation assertion of
Proposition~\ref{thm:general-source-pbw}.
\end{proof}

By composing with the isomorphism $D_\chi$ of
Section~\ref{subsec:central-specialization}, we obtain the
isomorphism $\Phi_0=\Phi_\chi\circ D_\chi:B_0\to K$ satisfying
\begin{equation}\label{eq:all-plucker-images}
 \Phi_0(\pi_0(P_I))=\chi_{q,t}(L(M_I)).
\end{equation}

Theorem~\ref{thm:intro-A} also gives an explicit
description of the inverse map.
For $M\in\mathcal P_\ell^+$, write
$M= \prod_{b=1}^m\prod_{i=1}^{k-1} Y_{i,i-2b}^{u_{b,i}}$.
Let $U$ be the tableau obtained by arranging
$u_{b,i}$ copies of $I_{b,i}$ in the chain
order, and put $T=U_{\mathrm{red}}$.
Then $M_T=M$, and
\begin{equation}\label{eq:inverse-based-isomorphism}
 \Phi_\chi^{-1}(\chi_{q,t}(L(M)))
 =
 t^{-\chi(T)/2}\pi_\chi(b_T).
\end{equation}

Finally, applying $\Phi_\chi$ to~\eqref{eq:young-twisted} and using
\eqref{eq:scaled-plucker-images}, we obtain the corresponding
relations for $(q,t)$-characters of simple modules.
For every $I\in\binom{[n]}{k}$, put
$Q_I=t^{\chi(I)/2}\chi_{q,t}(L(M_I))\in K$.
For every
$1\leq r\leq\min(k,n-k)$,
$U\in\binom{[n]}{k+r}$ and
$V\in\binom{[n]}{k-r}$, we have
\begin{equation}\label{eq:young-character}
 \sum_{\substack{\Lambda\subset U\setminus V\\
                 |\Lambda|=r}}
 (-t)^{a(U,V,\Lambda)}
 t^{-\Omega(U\setminus\Lambda,\Lambda\cup V)/2}
 Q_{U\setminus\Lambda}*Q_{\Lambda\cup V}
 =
 0,
\end{equation}
where $a(U,V,\Lambda)$ is defined as in
\eqref{eq:young-ordinary}.
All products $*$ are taken in the quantum Grothendieck ring.

\subsection{Presentations and isomorphism for \texorpdfstring{$\Gr(3,6)$}{Gr(3,6)}}
\label{subsec:isomorphism-examples}

Let $k=3$, $\ell=2$, and $n=6$. Write $Z_I=\pi_\chi(P_I)$,
abbreviating $\{a,b,c\}$ by $abc$. The central specialization
is given by $Z_{123}=v^{-1/2}$, $Z_{234}=Z_{345}=1$, and
$Z_{456}=v^{1/2}$. The six normalized generators of $B_\chi$
and their corresponding fundamental highest monomials and
segments in the convention~\eqref{eq:segment-convention} are
\[
\begin{array}{c|c|c|c|c}
 a&(b,i)&x_a=X_{b,i}&y_a=y_{b,i}&\Delta_{b,i}\\ \hline
 1&(1,1)&v^{1/2}Z_{124}&Y_{1,-1}&[0,0]\\
 2&(1,2)&v^{1/2}Z_{134}&Y_{2,0}&[0,1]\\
 3&(2,1)&Z_{235}&Y_{1,-3}&[-1,-1]\\
 4&(2,2)&Z_{245}&Y_{2,-2}&[-1,0]\\
 5&(3,1)&v^{-1/2}Z_{346}&Y_{1,-5}&[-2,-2]\\
 6&(3,2)&v^{-1/2}Z_{356}&Y_{2,-4}&[-2,-1]
\end{array}
\]
The order $x_1,\ldots,x_6$ is the lexicographic order.
For type $A_2$, the inverse quantum Cartan matrix is
\[
 C(z)^{-1}=\frac{1}{1+z^2+z^4}
 \begin{pmatrix}z+z^3&z^2\\z^2&z+z^3\end{pmatrix}.
\]
Using this matrix and~\eqref{eq:uniform-pair-relation},
we obtain the following fifteen defining relations.
The eight quasi-commutation relations are
\begin{equation}\label{eq:gr36-quasicommutation}
\begin{alignedat}{4}
 x_1x_2&=vx_2x_1, &\quad x_1x_4&=v^{-1}x_4x_1,
 &\quad x_1x_5&=vx_5x_1, &\quad x_2x_5&=v^{-1}x_5x_2,\\
 x_2x_6&=vx_6x_2, &\quad x_3x_4&=vx_4x_3,
 &\quad x_3x_6&=v^{-1}x_6x_3, &\quad x_5x_6&=vx_6x_5.
\end{alignedat}
\end{equation}
The four relations with a linear correction term are
\begin{equation}\label{eq:gr36-linear-relations}
\begin{alignedat}{2}
 v^{1/2}x_3x_1-v^{-1/2}x_1x_3&=(v-v^{-1})x_4,
 &\quad v^{1/2}x_4x_2-v^{-1/2}x_2x_4&=(v-v^{-1})x_1,\\
 v^{1/2}x_5x_3-v^{-1/2}x_3x_5&=(v-v^{-1})x_6,
 &\quad v^{1/2}x_6x_4-v^{-1/2}x_4x_6&=(v-v^{-1})x_3.
\end{alignedat}
\end{equation}
The three relations with a constant correction term are
\begin{equation}\label{eq:gr36-constant-relations}
\begin{gathered}
 vx_3x_2-v^{-1}x_2x_3=v-v^{-1},\qquad
 vx_5x_4-v^{-1}x_4x_5=v-v^{-1},\\
 vx_6x_1-v^{-1}x_1x_6=v-v^{-1}.
\end{gathered}
\end{equation}
These relations account for all $\binom{6}{2}=15$ pairs of
distinct generators. By Proposition~\ref{thm:general-defining},
they give a presentation of $B_\chi$ as a quotient of the
free associative algebra on $x_1,\ldots,x_6$.
The ordered monomials
$x_1^{u_1}\cdots x_6^{u_6}$, with $u_a\in\mathbb Z_{\geq0}$,
form a PBW basis of $B_\chi$.

On the representation-theoretic side, let
$G_a=\chi_{q,t}(L(y_a))$, $1\leq a\leq6$, be the six
fundamental $(q,t)$-characters of $\mathcal C_2$ for
$U_q(\widehat{\mathfrak{sl}_3})$, with the $y_a$ as in the
above table. The relations in
Proposition~\ref{thm:general-target-relations} are precisely
\eqref{eq:gr36-quasicommutation}--\eqref{eq:gr36-constant-relations}
with $v$ replaced by $t$, each $x_a$ replaced by $G_a$, and
multiplication replaced by $*$. For instance, the three types
of relations give, respectively,
\[
\begin{aligned}
 G_1*G_2&=tG_2*G_1,\\
 t^{1/2}G_3*G_1-t^{-1/2}G_1*G_3&=(t-t^{-1})G_4,\\
 tG_6*G_1-t^{-1}G_1*G_6&=t-t^{-1}.
\end{aligned}
\]
Proposition~\ref{thm:general-target-relations} gives the
presentation of $K_t(\mathcal C_2)$ by these fifteen
relations over $\mathbb Z[t^{\pm1/2}]$ and the basis
$G_1^{*u_1}*\cdots*G_6^{*u_6}$, $u_a\geq0$.
Here $G_a^{*u_a}$ is the $u_a$-fold quantum product.
After extension of scalars and identification $v=t$,
$x_a\mapsto G_a$ identifies the two presentations and
ordered monomial bases, illustrating the algebraic
part of Theorem~\ref{thm:general-based-isomorphism}.

For comparison, the simple-root generators of the
Hernandez--Leclerc presentation are, by
\eqref{eq:hl-window-generators},
$e_{1,-1}=G_5$, $e_{2,-1}=G_3$, $e_{1,0}=G_4$, and $e_{2,0}=G_2$.
The remaining two fundamental characters are obtained
from the linear correction relations:
\[
 G_6=\frac{t^{1/2}G_5*G_3-t^{-1/2}G_3*G_5}{t-t^{-1}},\qquad
 G_1=\frac{t^{1/2}G_4*G_2-t^{-1/2}G_2*G_4}{t-t^{-1}}.
\]
For instance, writing $A=G_5$, $B=G_3$, $C=G_6$
and omitting $*$, the relations $AC=tCA$ and
$BC=t^{-1}CB$ yield, upon substituting the expression
for $C$,
\[
 A^2B-(t+t^{-1})ABA+BA^2=0,\qquad
 B^2A-(t+t^{-1})BAB+AB^2=0.
\]
The same calculation gives the two quantum Serre relations
for $h=0$. The four relations between the generators with
integer indices $-1$ and $0$ are
\[
\begin{aligned}
 G_5*G_4&=t^{-2}G_4*G_5+1-t^{-2},&
 G_3*G_2&=t^{-2}G_2*G_3+1-t^{-2},\\
 G_5*G_2&=tG_2*G_5,&G_3*G_4&=tG_4*G_3.
\end{aligned}
\]
There are no pairs of integer indices differing by at least two
in this example.
By~\eqref{eq:complete-layer-window}, the four Serre and
four boson relations give an alternative presentation
with four generators over $\mathbb Q(t^{1/2})$.
In particular, $B_\chi$ is isomorphic to
$\widehat{\mathcal A}_{3,[-1,0]}$, the quantum boson algebra
with integer indices $-1,0$ (cf.~\cite[Remark~7.5(b)]{HL2015}).
The six-generator presentation, unlike the formulas with
$t-t^{-1}$ in the denominator, specializes to the
classical polynomial presentation.

We now illustrate the identification of the canonical
bases in Theorem~\ref{thm:general-based-isomorphism}.

By Theorem~\ref{thm:general-based-isomorphism},
we have
\[
 \Phi_\chi(Z_I)
 =t^{(\mathbf 1_{6\in I}-\mathbf 1_{1\in I})/2}\chi_{q,t}(L(M_I)),
 \qquad
 M_I=\prod_{a=a_1+1}^{a_2-1}Y_{2,4-2a}
     \prod_{a=a_2+1}^{a_3-1}Y_{1,5-2a}.
\]
The following table gives the images of all twenty
quantum Pl\"ucker coordinates under $\Phi_\chi\circ\pi_\chi$
(we abbreviate $\{a,b,c\}$ by $abc$).
\[
\begin{array}{c|l@{\quad}c|l}
 I&\Phi_\chi(Z_I)&I&\Phi_\chi(Z_I)\\\hline
123 & t^{-1/2} & 234 & 1 \\
124 & t^{-1/2}\chi_{q,t}(L(Y_{1,-1})) & 235 & \chi_{q,t}(L(Y_{1,-3})) \\
125 & t^{-1/2}\chi_{q,t}(L(Y_{1,-1}Y_{1,-3})) & 236 & t^{1/2}\chi_{q,t}(L(Y_{1,-3}Y_{1,-5})) \\
126 & \chi_{q,t}(L(Y_{1,-1}Y_{1,-3}Y_{1,-5})) & 245 & \chi_{q,t}(L(Y_{2,-2})) \\
134 & t^{-1/2}\chi_{q,t}(L(Y_{2,0})) & 246 & t^{1/2}\chi_{q,t}(L(Y_{2,-2}Y_{1,-5})) \\
135 & t^{-1/2}\chi_{q,t}(L(Y_{2,0}Y_{1,-3})) & 256 & t^{1/2}\chi_{q,t}(L(Y_{2,-2}Y_{2,-4})) \\
136 & \chi_{q,t}(L(Y_{2,0}Y_{1,-3}Y_{1,-5})) & 345 & 1 \\
145 & t^{-1/2}\chi_{q,t}(L(Y_{2,0}Y_{2,-2})) & 346 & t^{1/2}\chi_{q,t}(L(Y_{1,-5})) \\
146 & \chi_{q,t}(L(Y_{2,0}Y_{2,-2}Y_{1,-5})) & 356 & t^{1/2}\chi_{q,t}(L(Y_{2,-4})) \\
156 & \chi_{q,t}(L(Y_{2,0}Y_{2,-2}Y_{2,-4})) & 456 & t^{1/2} \\
\end{array}
\]
In particular, the four non-wrapping consecutive
quantum Pl\"ucker coordinates have images
$t^{-1/2},1,1,t^{1/2}$, whereas the two wrapping
cyclic consecutive quantum Pl\"ucker coordinates
$P_{126}$ and $P_{156}$ correspond to $(q,t)$-characters of nontrivial
simple modules.

For the reduced tableau $T$ with columns $124$ and $356$,
we have $M_T=Y_{1,-1}Y_{2,-4}$ and $\chi(T)=0$. Thus
$\Phi_\chi(\pi_\chi(b_T))=\chi_{q,t}(L(Y_{1,-1}Y_{2,-4}))$.

\section{An explicit formula for
\texorpdfstring{$(q,t)$}{(q,t)}-characters of simple
modules in type \texorpdfstring{$A$}{A}}
\label{sec:kl-characters}

We use Theorem~\ref{thm:intro-A} and the formula in~\cite{Brundan2005} for the
transition from standard monomials to the dual canonical
basis to express the $(q,t)$-character of an arbitrary simple module
in $\mathcal C_\ell$ in terms of fundamental characters
and ordinary Kazhdan--Lusztig polynomials of symmetric
groups. We use $v$ in the Grassmannian formulas and $t$
in the character formulas, identifying $v=t$ under the
isomorphism. The variable of the Kazhdan--Lusztig
polynomials will be denoted by $z$.

\subsection{Kazhdan--Lusztig coefficients and the transition matrix}
\label{subsec:kl-coefficients}

For $M=\prod_{b,i}y_{b,i}^{u_{b,i}}\in\mathcal P_\ell^+$,
write $U(M)=U(\mathbf u)$ for the tableau defined in
Section~\ref{subsec:source-standard-basis}.

Write
$U(M)=(U_1,\ldots,U_w)$, $w=\sum_{b,i}u_{b,i}$, $d=kw$.
Let $\mu=(k,\ldots,k)$, with $w$ parts, and let
$\nu=(\nu_1,\ldots,\nu_n)$ be the content of $U(M)$,
where $\nu_j$ is the number of occurrences of $j$
in $U(M)$.

We denote by $\SSYT(\mu,\nu)$ the set
of semistandard tableaux of column shape $\mu$
and content $\nu$.

Let $\mathfrak S_d$ be the symmetric group on
$\{1,\ldots,d\}$, with length function $\ell(-)$,
and let $w_0$ be its longest element.
We write permutations in one-line notation, so that
$x=x(1)\cdots x(d)$ for $x\in\mathfrak S_d$.
Consider the parabolic subgroups
$\mathfrak S_\mu=\mathfrak S_k\times\cdots\times\mathfrak S_k$ and $\mathfrak S_\nu=\mathfrak S_{\nu_1}\times\cdots\times\mathfrak S_{\nu_n}$
of $\mathfrak S_d$. Let $\mathcal D_\nu$ be the set
of minimal length representatives of the cosets
$\mathfrak S_\nu\backslash\mathfrak S_d$.

For a tableau $A\in\SSYT(\mu,\nu)$,
let $\gamma(A)$ be the word obtained by reading its
columns from left to right, and reading the entries
of each column in decreasing order. Thus, if
$A_j=(a_{1,j}<\cdots<a_{k,j})$, $1\leq j\leq w$,
then
\[
\gamma(A)=
(a_{k,1},\ldots,a_{1,1},
 \ldots,
 a_{k,w},\ldots,a_{1,w}).
\]
This is the column reading used in~\cite[Section~2]{Brundan2005},
with the bottom-aligned diagrams translated to our top-aligned convention.

Let
\[
\mathbf i_\nu=
(\underbrace{1,\ldots,1}_{\nu_1},
 \ldots,
 \underbrace{n,\ldots,n}_{\nu_n}).
\]
We use the right action of $\mathfrak S_d$ on words
given by
$(\mathbf a\cdot x)_j=a_{x(j)}$.
Denote by $x_A$ the unique element of $\mathcal D_\nu$
such that
$\gamma(A)=\mathbf i_\nu\cdot x_A$.
Equivalently, $x_A$ is obtained by replacing the
occurrences of each letter in $\gamma(A)$, from
left to right, by the corresponding consecutive
integers in $\{1,\ldots,d\}$.

For $x,y\in\mathfrak S_d$, let $P_{x,y}(z)$ be the
ordinary Kazhdan--Lusztig polynomial, normalized by $P_{x,x}=1$
and $\deg P_{x,y}\leq(\ell(y)-\ell(x)-1)/2$ for $x<y$,
with $P_{x,y}=0$ unless $x\leq y$ in the Bruhat order.
Equivalently, $z^{-\ell(y)/2}\sum_{x\leq y}P_{x,y}(z)T_x$
is the bar-invariant Kazhdan--Lusztig element in the Hecke algebra
with $(T_s-z)(T_s+1)=0$ for simple reflections $s$; the bar sends
$z^{1/2}$ to $z^{-1/2}$ and $T_x$ to $T_{x^{-1}}^{-1}$.
Here $z$ is independent of $q,v,t$.

For $A,B\in\SSYT(\mu,\nu)$, set
\begin{equation}
\label{eq:kl-tableau-coefficient}
h_{A,B}(v)
=
(-v)^{\ell(x_B)-\ell(x_A)}
\sum_{\sigma\in
\mathcal D_\nu\cap
\mathfrak S_\nu x_A\mathfrak S_\mu}
(-1)^{\ell(\sigma)+\ell(x_B)}
P_{\sigma w_0,x_Bw_0}(v^2).
\end{equation}

By~\cite[Remark~14 and Theorem~26]{Brundan2005},
$h_{A,B}(v)=k^*_{A,B}(v^{-1})$ in the notation of~\cite{Brundan2005}.
Since the rescaling factor $v^{-\beta(A)/2}$ in
\eqref{eq:main-Brundan-rescaling} depends only
on the shape and content of $A$, we obtain
\begin{equation}\label{eq:kl-standard-monomial-expansion}
 \widetilde P_A
 =
 \sum_{B\in\SSYT(\mu,\nu)}
 h_{A,B}(v)b_B.
\end{equation}

For $M,M'\in\mathcal P_\ell^+$, define
\begin{equation}
\label{eq:kl-monomial-coefficient}
c(M,M';t)
=
\sum_{\substack{
B\in\SSYT(\mu,\nu)\\
M_B=M'}}
h_{U(M),B}(t),
\end{equation}
where $\mu$ and $\nu$ are the column shape and
content of $U(M)$, respectively. The sum is
understood to be zero if there is no such tableau.

As observed in Section~\ref{sec:main-support},
tableaux of fixed content are distinguished by their
associated monomials. Thus the sum has at most one
term. Writing it in this form makes the passage to
the central quotient explicit.

By the triangularity established above, we have
\begin{equation}
\label{eq:kl-coefficient-triangularity}
c(M,M;t)=1,
\qquad
c(M,M';t)=0
\quad\text{unless }M'\leq M,
\end{equation}
where $\leq$ is the partial order on monomials.
Moreover,
$c(M,M';t)\in t^{-1}\mathbb Z[t^{-1}]$ if  $M'<M$.

For $M=1$, we use the empty tableau and set
$c(1,1;t)=1$.

The coefficients $c(M,M';t)$ give the transition from the standard
basis $\mathcal S_K$ to the canonical basis $\mathcal B_K$. Indeed, applying
$\Phi_\chi\circ\mathcal D$ to~\eqref{eq:kl-standard-monomial-expansion}
with $A=U(M)$ and using Theorem~\ref{thm:general-based-isomorphism} and~\eqref{eq:general-spectral-standard}, we obtain
\begin{equation}\label{eq:kl-standard-character-expansion}
 \mathcal E(M)
 =\sum_{B\in\SSYT(\mu,\nu)}h_{U(M),B}(t)\chi_{q,t}(L(M_B))
 =\sum_{M'\leq M}c(M,M';t)\chi_{q,t}(L(M')).
\end{equation}

\subsection{The formula for \texorpdfstring{$(q,t)$}{(q,t)}-characters of simple modules}
\label{subsec:kl-character-formulas}

We invert the transition matrix in~\eqref{eq:kl-standard-character-expansion}
to express the $(q,t)$-characters of simple modules in terms of normalized products
of fundamental characters.

\begin{theorem}
\label{thm:kl-character-formula}
For every $M\in\mathcal P_\ell^+$, we have
\begin{equation}
\label{eq:kl-character-formula}
\chi_{q,t}(L(M))
=
\sum_{s\geq 0}(-1)^s
\sum_{M=M_0>M_1>\cdots>M_s}
\left(
\prod_{a=0}^{s-1}
c(M_a,M_{a+1};t)
\right)
\mathcal E(M_s).
\end{equation}
All sums are finite. For $s=0$, the corresponding
term is $\mathcal E(M)$.
\end{theorem}

\begin{proof}
To prove \eqref{eq:kl-character-formula}, fix
$M\in\mathcal P_\ell^+$ and consider the finite
set of dominant monomials $M'\leq M$ in
$\mathcal P_\ell^+$. With respect to a linear
extension of the partial order on this set, the
matrix
$C(t)=\bigl(c(M',M'';t)\bigr)_{M',M''}$
is unitriangular. Write
$C(t)=I+N(t)$,
where $N(t)$ is strictly triangular and hence
nilpotent. Therefore,
$C(t)^{-1} = \sum_{s\geq 0}(-1)^sN(t)^s$.
For $s\geq 1$, the $(M,M')$-entry of $N(t)^s$ is
\[
\sum_{M=M_0>M_1>\cdots>M_s=M'}
\prod_{a=0}^{s-1}c(M_a,M_{a+1};t).
\]
Applying the inverse matrix to
\eqref{eq:kl-standard-character-expansion}
gives \eqref{eq:kl-character-formula}.
\end{proof}

For a fundamental monomial $M$, the tableau $U(M)$
has a single column, so
\eqref{eq:kl-character-formula} reduces to
$\chi_{q,t}(L(M))=\mathcal E(M)$.

\begin{example}\label{ex:kl-Gr36}
Consider $\Gr(3,6)$, corresponding to type $A_2$
with $\ell=2$. Let $y_1=Y_{1,-1}$, $y_2=Y_{2,-4}$
and $M=y_1y_2$. The corresponding fundamental columns
are $(1,2,4)$ and $(3,5,6)$. Thus
\[
 A=U(M)=
 \begin{array}{|c|c|}
 \hline
 1&3\\ \hline
 2&5\\ \hline
 4&6\\ \hline
 \end{array}.
\]
Here $\mu=(3,3)$ and $\nu=(1,1,1,1,1,1)$, so
$\mathfrak S_\nu$ is trivial and
$\mathcal D_\nu=\mathfrak S_6$.
The five semistandard tableaux of this shape and
content, written as pairs of columns, and their
coefficients in~\eqref{eq:kl-tableau-coefficient} are
\[
\begin{array}{c|c|c|c}
 B & x_B & \ell(x_B) & h_{A,B}(t)\\ \hline
 (123,456)&321654&6&t^{-1}\\
 (124,356)&421653&7&1\\
 (125,346)&521643&8&0\\
 (134,256)&431652&8&0\\
 (135,246)&531642&9&0
\end{array}
\]
Indeed, $x_A$ is the longest element of
$x_A\mathfrak S_\mu$. A nonzero summand in
\eqref{eq:kl-tableau-coefficient} requires
$\sigma w_0\leq x_Bw_0$, and hence
$\ell(x_B)\leq\ell(\sigma)\leq\ell(x_A)=7$.
This gives the last four coefficients in the table.
For $B=(123,456)$, the only contributing permutation
is $\sigma=x_A$, since $x_B\notin x_A\mathfrak S_\mu$.
Moreover, $x_Aw_0=356124$ and $x_Bw_0=456123$
form a covering relation in the Bruhat order.
Thus $P_{356124,456123}(z)=1$ and
$h_{A,B}(t)=(-t)^{-1}(-1)^{13} P_{356124,456123}(t^2)=t^{-1}$.
Since both columns of $B$ are trivial, $M_B=1$.
Hence $c(M,M;t)=1$, $c(M,1;t)=t^{-1}$, and all
other coefficients vanish.

By~\eqref{eq:kl-standard-character-expansion},
$\mathcal E(M)=\chi_{q,t}(L(M))+t^{-1}$.
Since $J(y_1,y_2)=2$, we have
$\mathcal E(M)=t^{-1}\chi_{q,t}(L(y_1))*\chi_{q,t}(L(y_2))$.
Consequently,
\begin{equation}
\label{eq:kl-Gr36-example}
 \chi_{q,t}(L(Y_{1,-1}Y_{2,-4}))
 =t^{-1}\chi_{q,t}(L(Y_{1,-1}))*\chi_{q,t}(L(Y_{2,-4}))-t^{-1}.
\end{equation}
Put $Z_I=\pi_\chi(P_I)$. Theorem~\ref{thm:intro-A} gives the corresponding
identity $\pi_\chi(b_A)=v^{-1}Z_{124}Z_{356}-v^{-1}$.

At $t=1$, this gives the $q$-character identity
$\chi_q(L(M))=\chi_q(L(y_1))\chi_q(L(y_2))-1$.
\end{example}

\section{Quantum cyclic shifts and cluster quasi-automorphisms}
\label{sec:cyclic-shifts}

In this section, we prove Theorem~\ref{thm:intro-C}.
We first obtain the normalized quantum cyclic shift
from the cocycle isomorphism of Launois and
Lenagan~\cite{LaunoisLenagan2011}.
We then establish its action on the dual canonical basis
of $A_v$ and specialize the shifted basis elements to
the central quotient $B_\chi$.
We explain why the shift does not descend to this quotient
and establish the quantum cluster structure after
localization at the remaining frozen variables.
A correction of the normalized quantum cyclic shift by
frozen Laurent monomials then gives a quantum cluster
quasi-automorphism of the localized quotient.
Finally, we determine its order and transport it to
the localized quantum Grothendieck ring.

\subsection{The normalized quantum cyclic shift}
\label{subsec:normalized-cyclic-shift}

Use the cyclic permutation $\sigma$, the grading lattice
$\mathcal L$, and the integers $d_\alpha$ introduced above
Theorem~\ref{thm:intro-C}. For every $k$-element subset $I$,
$\sigma^j I=\{\sigma^j(i):i\in I\}$ is a set, and the quantum
Pl\"ucker coordinate $P_{\sigma^j I}$ is defined using its elements
in increasing order, as in~\eqref{eq:quantum-plucker-coordinate},
without any quantum antisymmetry factor.
For $\alpha,\beta\in\mathcal L$, put
$H(\alpha,\beta)=d_\alpha\beta_n-\alpha_n d_\beta$.
This is an integral skew-symmetric bilinear form.

To use the isomorphism of Launois and Lenagan, extend the field temporarily to $\mathbb Q(v^{1/(2k)})$.
In the notation of~\cite[Lemma~5.3]{LaunoisLenagan2011},
$p^k=v^2$; we take $p=v^{2/k}$. Their cocycle is
$v^{D(\alpha,\beta)}$, where
\begin{equation}\label{eq:cyclic-LL-cocycle}
 D(\alpha,\beta)=\frac{2}{k}\alpha_n(|\beta|-\beta_n)
 =2\alpha_n d_\beta-\frac{2}{k}\alpha_n\beta_n.
\end{equation}
Unlike $\Omega$, the form $D$ is not skew-symmetric and
in general is not integral on $\mathcal L$.
The skew part of $D$ satisfies
\begin{equation}\label{eq:cyclic-cocycle-comparison}
 D(\beta,\alpha)-D(\alpha,\beta)=2H(\alpha,\beta).
\end{equation}

Let $\rho_v^{\rm LL}$ denote the linear map on the underlying
vector space of $A_v$ induced by the isomorphism of
\cite[Theorem~5.9]{LaunoisLenagan2011}.
By that theorem and~\eqref{eq:cyclic-LL-cocycle},
\begin{equation}\label{eq:cyclic-LL-product}
 \rho_v^{\rm LL}(fg)=v^{-D(\alpha,\beta)}
 \rho_v^{\rm LL}(f)\rho_v^{\rm LL}(g),
 \qquad \rho_v^{\rm LL}(P_I)=v^{-2\mathbf1_{n\in I}}P_{\sigma I}.
\end{equation}
Set
\begin{equation}\label{eq:cyclic-normalization}
 s(\alpha)=\alpha_n(d_\alpha+1)-\frac{\alpha_n(\alpha_n-1)}{k}.
\end{equation}
A direct calculation gives
\[
 s(\alpha+\beta)-s(\alpha)-s(\beta)
 =\alpha_n d_\beta+\beta_n d_\alpha-\frac{2\alpha_n\beta_n}{k}
 =H(\alpha,\beta)+D(\alpha,\beta).
\]
For a quantum Pl\"ucker coordinate, $s(\deg P_I)=2\mathbf1_{n\in I}$.
Consequently, on the extended coefficient field, the map
\begin{equation}\label{eq:cyclic-normalization-comparison}
 \widehat\rho_v(f)=v^{s(\alpha)}\rho_v^{\rm LL}(f),
 \qquad f\in (A_v)_\alpha,
\end{equation}
has the following properties.

\begin{proposition}\label{prop:cyclic-normalized-shift}
There is a unique $\bbk$-linear bijection
$\widehat\rho_v:A_v\to A_v$ sending the homogeneous component
of content $\alpha$ to the component of content $\sigma\alpha$
and satisfying
\begin{equation}\label{eq:cyclic-normalized-product}
 \widehat\rho_v(fg)=v^{H(\alpha,\beta)}
 \widehat\rho_v(f)\widehat\rho_v(g),
 \qquad \widehat\rho_v(P_I)=P_{\sigma I},
\end{equation}
for homogeneous $f,g$ of contents $\alpha,\beta$.
Moreover, $\widehat\rho_v$ commutes with the bar involution,
preserves the integral form spanned by the standard monomials
over $\mathbb Z[v^{1/2},v^{-1/2}]$, and satisfies $\widehat\rho_v^n=1$.
\end{proposition}

\begin{proof}
The identities for $s$ and~\eqref{eq:cyclic-LL-product} give
\eqref{eq:cyclic-normalized-product} directly from
\eqref{eq:cyclic-normalization-comparison}. Although the construction
uses the extended coefficient field, the exponents $H(\alpha,\beta)$
are integral, and the quantum Pl\"ucker coordinates are mapped to
quantum Pl\"ucker coordinates. Since they generate $A_v$, the map
restricts to a $\bbk$-linear map on $A_v$, uniquely determined by
\eqref{eq:cyclic-normalized-product}.
For a product of quantum Pl\"ucker coordinates this gives
\begin{equation}\label{eq:cyclic-word}
 \widehat\rho_v(P_{I_1}\cdots P_{I_d})
 =v^{\sum_{a=1}^d(2a-d-1)\mathbf1_{n\in I_a}}
 P_{\sigma I_1}\cdots P_{\sigma I_d}.
\end{equation}
In particular the map preserves the integral span of all
products of quantum Pl\"ucker coordinates, which equals the
integral span of standard monomials by integral quantum
straightening~\cite[Theorem~26]{Brundan2005}.
Since $H$ is skew-symmetric and each quantum Pl\"ucker coordinate
is bar invariant, the map commutes with the bar by induction
on the number of factors in a product.

Iterating~\eqref{eq:cyclic-normalized-product} $n$ times, the
sum of exponents is $d_\alpha|\beta|-|\alpha|d_\beta=0$.
Hence $\widehat\rho_v^n$ is an algebra automorphism of $A_v$
fixing every quantum Pl\"ucker coordinate, and thus is the
identity. Its inverse also preserves the integral form.
\end{proof}

At $v^{1/2}=1$, Proposition~\ref{prop:cyclic-normalized-shift}
gives the classical signed cyclic shift $\rho$ on
$A_{\rm cl}=\mathbb C[\Gr(k,n)]$, with $\rho^n=1$;
see~\cite[Section~2.5]{Lam2019}.
The substitution $X_{r,a}\mapsto X_{r,a+1}$ for $a<n$ and
$X_{r,n}\mapsto(-1)^{k-1}X_{r,1}$ gives $\rho(P_I)=P_{\sigma I}$,
as the two signs cancel when $n\in I$.
It is an orientation-preserving cluster automorphism, since
rotation of the boundary labels takes a Grassmannian seed to its
rotated seed, preserving its quiver~\cite[Section~5]{FraserBraid2020}.
In the quantum case, $\widehat\rho_v$ is an invertible linear map
of $A_v$, but it is not an algebra automorphism. For example, put
$I=[1,k]$ and $J=[1,k-1]\cup\{n\}$. By~\eqref{eq:frozen-commutation},
$P_IP_J=vP_JP_I$ and $P_{\sigma I}P_{\sigma J}=v^{-1}P_{\sigma J}P_{\sigma I}$,
since $\sigma J=I$ and $\sigma I=[2,k+1]$.
Since $A_v$ is a domain and $v$ is an indeterminate, a map sending
$P_D$ to $P_{\sigma D}$ for every $k$-element subset $D$ cannot preserve the ordinary
quantum multiplication. The correction
factor in~\eqref{eq:cyclic-normalized-product} accounts for this
change in the commutation relations. In particular, the uncorrected
quantum cyclic shift is not a quantum cluster quasi-automorphism,
which by definition must be an algebra homomorphism.

\subsection{Comparison with the Launois--Lenagan cocycle}

Our twist is defined by the skew-symmetric form $\Omega$
in~\eqref{eq:Omega} and makes the chosen consecutive minors
central. The twist of Launois and Lenagan serves a different
purpose: it makes the cyclic permutation of the quantum
Pl\"ucker coordinates into an isomorphism from a twisted
algebra to the original quantum Grassmannian, up to scalar
factors on the generators. The relation between their cyclic shift and our normalization is
\eqref{eq:cyclic-normalization-comparison}.
More explicitly, $D+H=(D+D^{\mathsf T})/2$ is the symmetric
part of $D$. The identity for $s$ shows that this symmetric
part is a coboundary: the rescaling $f\mapsto v^{-s(\alpha)}f$
is an isomorphism from the algebra with multiplication
$v^{D(\alpha,\beta)}fg$ to the algebra with multiplication
$v^{-H(\alpha,\beta)}fg$. By~\eqref{eq:cyclic-normalized-product},
$\widehat\rho_v$ maps the latter algebra isomorphically to $A_v$.
Their composition is $\rho_v^{\rm LL}$. Thus the symmetric part of the Launois--Lenagan cocycle
accounts for the scalar normalization, and its skew-symmetric
part gives the exponent $-H$ of the twisted multiplication.
The normalized shift itself is defined over our original
coefficient field by Proposition~\ref{prop:cyclic-normalized-shift}.

\subsection{The action on the dual canonical basis}

Let $\operatorname{pr}$ be Sch\"utzenberger promotion on rectangular
semistandard tableaux, with the convention that the content is
transformed by $\sigma$. Thus $\operatorname{pr}$ removes the
entries $n$, moves the empty boxes to the upper left by reverse
jeu de taquin, increases the remaining entries by one, and fills
the empty boxes with $1$. This is the inverse of promotion in
the convention that removes the entries $1$, decreases the
remaining entries by one, and fills the empty boxes with $n$.
For a one-column tableau $I$, $\operatorname{pr}(I)$ is the column
obtained by listing the elements of $\sigma I$ in increasing order.
For a general tableau, $\operatorname{pr}$ is the usual tableau
promotion and is not defined by applying $\sigma$ separately to
its columns. If $T$ has content $\gamma$, then
$\operatorname{pr}(T)$ has content $\sigma\gamma$.

The term promotion has a different meaning
in~\cite{EPSBTW2026,LOT2026}, where it refers to
geometric maps associated with plabic tangles.
Here we use promotion only for the Sch\"utzenberger
operation on tableaux.

The classical cyclic shift of the
dual canonical basis is described in
\cite[Proposition~5.5]{Rhoades2010} and
\cite[Theorem~2]{Lam2019}.

\begin{theorem}\label{thm:cyclic-canonical-basis}
For every $T\in\SSYT(k,[n])$, we have
\begin{equation}\label{eq:cyclic-canonical-basis}
 \widehat\rho_v(b_T)=b_{\operatorname{pr}(T)}.
\end{equation}
Equivalently, if $T$ has $d$ columns and $r$ entries equal to $n$,
then
\begin{equation}\label{eq:cyclic-LL-canonical-basis}
 \rho_v^{\rm LL}(b_T)=
 v^{-r(d+1)+r(r-1)/k}b_{\operatorname{pr}(T)}.
\end{equation}
\end{theorem}

\begin{proof}
We use the cyclic invariance of the canonical bases of
coinvariant spaces~\cite[Proposition~28.2.4 and Section~28.2.9]{Lusztig1993}.

\smallskip
\noindent\emph{The tensor product realization.}
Let $V$ be the natural module of $U=U_v(\mathfrak{sl}_k)$
with basis $u_1,\ldots,u_k$, where
$E_i u_j=\delta_{j,i+1}u_i$, $F_i u_j=\delta_{j,i}u_{i+1}$ and
$K_i u_j=v^{\delta_{i,j}-\delta_{i+1,j}}u_j$.
The coproduct in~\cite[Section~3]{Brundan2005} is
$\Delta(E_i)=1\otimes E_i+E_i\otimes K_i^{-1}$,
$\Delta(F_i)=K_i\otimes F_i+F_i\otimes1$ and
$\Delta(K_i)=K_i\otimes K_i$.
In the convention of~\cite[Section~5]{Brundan2005}, the quantum
symmetric power $S^a(V)$ is the degree-$a$ part of the tensor algebra
quotient by $u_j u_i-v^{-1}u_i u_j$ for $i<j$. Its monomial basis
consists of $u_{i_1}\cdots u_{i_a}$ for $i_1\leq\cdots\leq i_a$.
The divided power $\widetilde S^a(V)\subset V^{\otimes a}$ has basis
\[
 x^*_{\mathbf i}=\sum_{\mathbf j\sim\mathbf i}
 v^{\operatorname{inv}(\mathbf j)}
 u_{j_a}\otimes\cdots\otimes u_{j_1},
 \qquad i_1\leq\cdots\leq i_a,
\]
where the sum is over distinct rearrangements of $\mathbf i$ and
$\operatorname{inv}(\mathbf j)=|\{(p,q):p<q,\ j_p>j_q\}|$.
These are the bar-invariant vectors of
\cite[Equation~(5.3)]{Brundan2005}. For $a=0$, both powers are
the trivial module. Set
\begin{equation}\label{eq:cyclic-tensor-spaces}
 C_\gamma=S^{\gamma_n}(V)\otimes\cdots\otimes S^{\gamma_1}(V),
 \qquad
 D_\gamma=\widetilde S^{\gamma_1}(V)\otimes\cdots\otimes
 \widetilde S^{\gamma_n}(V).
\end{equation}
Here $\gamma=(\gamma_1,\ldots,\gamma_n)$ is a nonnegative content.
The two monomial bases are indexed by the same arrays of row
multiplicities $h=(h_{i,a})$, with $\sum_i h_{i,a}=\gamma_a$,
using the symmetric monomials in $C_\gamma$ and the divided power
vectors $x^*_{\mathbf i}$ in $D_\gamma$ in the indicated factor orders.
By transposing the quantum matrix and applying
\cite[Theorem~15]{Brundan2005}, we identify $C_\gamma$ with the
content-$\gamma$ component of the quantum matrix algebra.
Under this identification the pure symmetric monomial basis consists of
$\prod_{a=n}^{1}\prod_{i=1}^k X_{i,a}^{h_{i,a}}$,
with the columns in decreasing order and the rows within each column
in increasing order. The maximal quantum minor corresponding to
$I=\{a_1<\cdots<a_k\}$ has the expansion
\[
 P_I=\sum_{\pi\in\mathfrak S_k}(-v^{-1})^{\ell(\pi)}
 X_{\pi(k),a_k}\cdots X_{\pi(1),a_1}.
\]
This follows from~\eqref{eq:quantum-plucker-coordinate} by
reordering the commuting entries in each term and applying the matrix bar.
Each such minor is a quantum determinant vector for the row action of $U$.
Indeed, the two terms in the action of $E_i$ or $F_i$ arising from
permutations that differ by exchanging $i$ and $i+1$ cancel by the
coproduct formulas above; the weight of the determinant vector is zero
for $U$. Hence $(A_v)_\gamma\subseteq C_\gamma^U$, where
$C_\gamma^U$ denotes the invariant subspace for the counit $\varepsilon$
of $U$.
If $|\gamma|=kd$, the multiplicity of the trivial $U$-module in
$C_\gamma$ is the number of semistandard tableaux of shape $(d^k)$
and content $\gamma$, by the Pieri rule and the generic quantum
character formula. This is also the dimension of $(A_v)_\gamma$ by
standard monomial theory. If $k\nmid|\gamma|$, both spaces are zero.
Thus $(A_v)_\gamma=C_\gamma^U$.

The tensor product bar is obtained from the factor bars using the
quasi-$R$-matrix $\Theta$: on two factors it is
$\psi(x\otimes y)=\Theta(\overline x\otimes\overline y)$,
iterated for more factors. Let $\psi_D$ denote this bar on $D_\gamma$;
see~\cite[Sections~3--5]{Brundan2005}.
The bilinear pairing $(\cdot,\cdot)_B$ of
\cite[Section~5]{Brundan2005} gives a sesquilinear pairing
\[
 \langle x,z\rangle_\gamma=(x,\psi_D(z))_B.
\]
For this pairing the two pure monomial bases, indexed by $h,h'$
in the opposite factor orders, pair to $\delta_{h,h'}$ by
\cite[Equation~(5.4)]{Brundan2005} and the pairing formula following it. The dual canonical basis of $C_\gamma$
and the canonical basis of $D_\gamma$ are also dual, by
\cite[Theorem~7]{Brundan2005}.
Write $\tau_U$ for the algebra anti-automorphism exchanging $E_i$
and $F_i$ and fixing $K_i$, and $\overline{\phantom{x}}_U$ for the
quantum group bar, the semilinear algebra involution fixing $E_i,F_i$
and sending $K_i$ to $K_i^{-1}$. By~\cite[Lemma~3]{Brundan2005},
$\langle ux,z\rangle_\gamma=\langle x,\overline{\tau_U(u)}_U z\rangle_\gamma$.
Since $\overline{\tau_U(\ker\varepsilon)}_U=\ker\varepsilon$,
$C_\gamma^U$ is the annihilator of $(\ker\varepsilon)D_\gamma$,
and is thus dual, for this pairing, to the coinvariant space
$(D_\gamma)_U=D_\gamma/(\ker\varepsilon)D_\gamma$.
Since $D_\gamma$ is finite-dimensional and semisimple,
$(\ker\varepsilon)D_\gamma$ is the sum of its nontrivial isotypic
components. By~\cite[Section~27.2.5]{Lusztig1993}, the canonical basis
of $D_\gamma$ induces a basis of $(D_\gamma)_U$, with the other
basis elements spanning the kernel. Its dual basis is therefore
the subset of the dual canonical basis of $C_\gamma$ contained in
$C_\gamma^U$.

By Section~\ref{subsec:homogeneous-standard-monomials}, the elements
$b_T$ are the matrix dual canonical basis elements in $(A_v)_\gamma$,
with the normalization given there. To compare the two normalizations,
note that the ordinary matrix transpose
is the composition of the matrix bar and the antilinear
anti-isomorphism of~\cite[Equations~(6.3)--(6.4)]{Brundan2005}.
Consequently it preserves the dual canonical basis normalized by
the matrix bar, by~\cite[Remark~17]{Brundan2005}.
In the transposed matrix the row and column degrees are $\gamma$
and $(d,\ldots,d)$, respectively. By~\cite[Remark~17]{Brundan2005}
and~\eqref{eq:main-Brundan-rescaling}, the conversion from the
basis in $C_\gamma^U$ to the basis $b_T$ is multiplication by
$v^{\beta(T)/2}$, where $\beta(T)=(kd^2-\sum_a\gamma_a^2)/2$.
This factor is unchanged under $\gamma\mapsto\sigma\gamma$.

\smallskip
\noindent\emph{Rotation of invariant tensors.}
Let $2\varrho$ be the sum of the positive coroots of type $A_{k-1}$,
put $g_0=K_{-2\varrho}$, and write
$a(h)=\langle2\varrho,h\rangle=\sum_{i=1}^k(k+1-2i)h_i$
for a row weight $h=(h_1,\ldots,h_k)$.
Thus $g_0$ acts by $v^{-a(h)}$ on the weight space of weight $h$.
Set $\vartheta_r=(-1)^{(k-1)r}$. Define
$R_\gamma:C_\gamma\to C_{\sigma\gamma}$ by
\begin{equation}\label{eq:cyclic-tensor-rotation}
 R_\gamma(x_n\otimes\cdots\otimes x_1)
 =\vartheta_{\gamma_n}x_{n-1}\otimes\cdots\otimes x_1
 \otimes g_0^{-1}x_n.
\end{equation}

To apply~\cite[Section~28.2.9]{Lusztig1993}, denote its parameter
by $u$. The change $u\mapsto v^{-1}$, $E_i^{\rm L}\mapsto E_i$,
$F_i^{\rm L}\mapsto F_i$, $K_i^{\rm L}\mapsto K_i^{-1}$
intertwines the quantum group bars and preserves the weights.
It takes the coproduct in~\cite[Section~3.1.4]{Lusztig1993}
to the opposite of the coproduct above. Reversing tensor order
therefore identifies the tensor products. Each divided power
$\widetilde S^a(V)$ is the irreducible based module of highest
weight $a\omega_1$; its canonical monomial basis is the one in
\cite[Equation~(5.3)]{Brundan2005}. The normalization is unchanged
under $u\mapsto v^{-1}$, since these monomial basis vectors are
obtained from the highest weight vector by successive divided
powers along simple root strings and the quantum factorials
are invariant under parameter inversion.

Under the parameter change and the flip, the quasi-$R$-matrix of
\cite[Section~4.1.2]{Lusztig1993} satisfies the intertwining identity
for the quasi-$R$-matrix of~\cite[Section~3]{Brundan2005} and has
constant term $1$ in the completion of
$\sum_{\nu\geq0}U^+_\nu\otimes U^-_{-\nu}$. By uniqueness they agree.
Thus the tensor product bars correspond, and the congruence modulo
$u^{-1}\mathbb Z[u^{-1}]$ becomes the congruence modulo
$v\mathbb Z[v]$. By the characterization of the tensor product bases
in~\cite[Theorem~27.3.2]{Lusztig1993} and
\cite[Section~5]{Brundan2005}, the based module with factors
$\gamma_n,\ldots,\gamma_1$ in the former convention becomes
$D_\gamma$ in the latter convention.

In type $A$ the weight function $\mathbf n$ of
\cite[Section~2.3.1]{Lusztig1993} is $\langle2\varrho,-\rangle$.
Since $a(h)\equiv(k-1)\sum_i h_i\pmod2$, the scalar
$(-1)^{a(h)}u^{-a(h)}$ in the rotation formula of
\cite[Section~28.2.9]{Lusztig1993} becomes
$\vartheta_{\gamma_n}v^{a(h)}$ on the moved factor.
Reversing the tensor order, the first-to-last rotation of that
formula therefore gives the based isomorphism
\begin{equation}\label{eq:cyclic-coinvariant-rotation}
 [z_1\otimes\cdots\otimes z_n]\longmapsto
 [\vartheta_{\gamma_n}g_0^{-1}z_n\otimes z_1\otimes\cdots
 \otimes z_{n-1}]
\end{equation}
from $(D_\gamma)_U$ to $(D_{\sigma\gamma})_U$.
Denote the corresponding map on the full tensor spaces by $J_\gamma$.
On the dual pure monomial bases, $R_\gamma$ and $J_\gamma$ rotate
matching factors and multiply them by the same scalar
$\vartheta_{\gamma_n}v^{a(h)}$. Its product with its bar is $1$,
so
\[
 \langle R_\gamma x,J_\gamma z\rangle_{\sigma\gamma}
 =\langle x,z\rangle_\gamma.
\]
It follows that $R_\gamma$ maps the annihilators of the coinvariant
kernels onto one another and, as the inverse adjoint of the based
isomorphism~\eqref{eq:cyclic-coinvariant-rotation}, maps their dual
canonical bases onto one another. The normalization factor
$v^{\beta(T)/2}$ is unchanged by rotation. Hence the direct sum $R$
of the $R_\gamma$ preserves the dual canonical basis of $A_v$.

\smallskip
\noindent\emph{Identification with the normalized cyclic shift.}
In the braided tensor product realization of
\cite[Section~3 and Theorem~15]{Brundan2005}, replacing the
$U_v(\mathfrak{gl}_k)$ braiding by the $U_v(\mathfrak{sl}_k)$
braiding on symmetric powers of degrees $a,b$ multiplies it by
$v^{-ab/k}$, since the scalar product of the projected weights is
the original scalar product minus $ab/k$. The tensor product realization uses decreasing column order. In a product of tensors
of contents $\alpha,\beta$, the crossings involve exactly the pairs
of columns $i<j$ from the first and second tensor, respectively.
Thus the multiplication defined by the $\mathfrak{sl}_k$ braiding,
over the extended coefficient field, is
\begin{equation}\label{eq:cyclic-sl-multiplication}
 f\cdot_{\mathfrak{sl}}g=v^{-Q(\alpha,\beta)/k}fg,
 \qquad Q(\alpha,\beta)=\sum_{i<j}\alpha_i\beta_j.
\end{equation}

We show that $R$ preserves this multiplication. For two $U$-module
algebras $A,B$ with braided tensor product multiplication, the map
$\sum_i a_i\otimes b_i\mapsto\sum_i g_0b_i\otimes a_i$ gives an
algebra isomorphism $(A\otimes B)^U\to(B\otimes A)^U$.
Indeed, for the antipode $S$ of $U$, the above coproduct gives
$S(E_i)=-E_iK_i$, $S(F_i)=-K_i^{-1}F_i$, and hence
$S^2(u)=g_0ug_0^{-1}$. Invariance of
$f=\sum_i a_i\otimes b_i$ and $g=\sum_j c_j\otimes d_j$ gives
\[
 \sum_i a_i\otimes u b_i=\sum_i S(u)a_i\otimes b_i,
 \qquad
 \sum_j u c_j\otimes d_j=\sum_j c_j\otimes S^{-1}(u)d_j.
\]
These are the antipode identities applied to invariant tensors.
The first identity and $S^2(u)=g_0ug_0^{-1}$ give
\[
 \Delta(u)\sum_i g_0b_i\otimes a_i
 =\sum_{i,(u)}g_0b_i\otimes u_{(2)}S^{-1}(u_{(1)})a_i
 =\varepsilon(u)\sum_i g_0b_i\otimes a_i.
\]
Writing $\mathcal R=\sum_\nu r_\nu\otimes s_\nu$ for the
universal $R$-matrix (acting on the finite-dimensional homogeneous
components), the rotation of the product $fg$ is
\[
 \begin{aligned}
 \sum_{i,j,\nu}g_0(r_\nu b_i)g_0(d_j)\otimes a_i(s_\nu c_j)&=\sum_{i,j,\nu}g_0(b_i)(S(s_\nu)g_0(d_j))
 \otimes (S(r_\nu)a_i)c_j  \\
 & =\sum_{i,j,\nu}g_0(b_i)(s_\nu g_0(d_j))\otimes(r_\nu a_i)c_j.
 \end{aligned}
\]
The first equality follows from the two invariance identities,
and the second from $(S\otimes S)(\mathcal R)=\mathcal R$.
The last expression is the product of the rotated tensors.
The inverse map sends $b\otimes a$ to $a\otimes g_0^{-1}b$;
it preserves invariants by the same antipode identities.

Apply the inverse of this isomorphism to the first factor of the
decreasing-column tensor product, corresponding to column $n$.
Since the sign $\vartheta_{\gamma_n}$ is multiplicative in the degree,
$R$ preserves the $\mathfrak{sl}_k$ multiplication after cyclic
relabelling. By~\eqref{eq:cyclic-sl-multiplication},
\[
 R(fg)=v^{(Q(\alpha,\beta)-Q(\sigma\alpha,\sigma\beta))/k}R(f)R(g).
\]
Now $Q(\sigma\alpha,\sigma\beta)-Q(\alpha,\beta)
=\alpha_n|\beta|-|\alpha|\beta_n$, so the exponent is
$H(\alpha,\beta)$.

It remains to compute $R(P_I)$ using the decreasing-column expansion above.
If $n\notin I$, we have $R(P_I)=P_{\sigma I}$ directly.
If $a_k=n$, moving the last entry $i=\pi(k)$ of the permutation
to the first position changes its length by $2i-k-1$.
The ratio of the two minor coefficients is
$(-v^{-1})^{2i-k-1}=(-1)^{k-1}v^{k+1-2i}$,
which is precisely the factor in~\eqref{eq:cyclic-tensor-rotation}.
Thus $R(P_I)=P_{\sigma I}$ in both cases. Since the quantum
Pl\"ucker coordinates generate $A_v$,
Proposition~\ref{prop:cyclic-normalized-shift} gives $R=\widehat\rho_v$.

Finally, the specialization of $\widehat\rho_v$ at $v^{1/2}=1$
is the classical cyclic shift sending $P_I$ to $P_{\sigma I}$.
By~\cite[Theorem~2]{Lam2019}, based on
\cite[Proposition~5.5]{Rhoades2010}, it sends the classical dual
canonical basis element indexed by $T$ to that indexed by
$\operatorname{pr}(T)$. The classical specializations of the $b_T$
are linearly independent by Section~\ref{subsec:homogeneous-standard-monomials}.
Since $\widehat\rho_v$ permutes the quantum basis, its permutation
is therefore exactly $T\mapsto\operatorname{pr}(T)$. This proves
\eqref{eq:cyclic-canonical-basis}. Formula~\eqref{eq:cyclic-LL-canonical-basis}
follows from~\eqref{eq:cyclic-normalization-comparison} and
\eqref{eq:cyclic-normalization}.
\end{proof}

\subsection{Specialization to the central quotient}

We now apply the quantum cyclic shift to the dual canonical basis of $A_v$
before taking the central quotient. Throughout
this subsection the non-wrapping consecutive columns are
$C_a=[a,a+k-1]$, $1\leq a\leq s$. Neither the cocycle $\Omega$
nor the quotient map $\pi_\chi$ is changed.
Recall from Section~\ref{subsec:central-specialization} that
$\mathcal D(f)=v^{-\chi(\deg f)/2}\pi_\chi(f)$ for homogeneous
$f$. In particular, $\mathcal D(P_{C_a})=1$ and
$\mathcal D(b_T)=\widehat b_{T_{\mathrm{red}}}$.

For $0\leq j<n$, Theorem~\ref{thm:cyclic-canonical-basis} gives
$\widehat\rho_v^j(b_T)=b_{\operatorname{pr}^j(T)}$.
Combining~\eqref{eq:normalized-reduction} with
Theorem~\ref{thm:general-based-isomorphism} therefore gives the
specialization assertion of Theorem~\ref{thm:intro-C}(i).
For every $T\in\SSYT(k,[n])$ and $0\leq j<n$, we have
\begin{equation}\label{eq:cyclic-tableau-images}
\begin{aligned}
 \mathcal D(\widehat\rho_v^j(b_T))
 =\widehat b_{(\operatorname{pr}^j(T))_{\mathrm{red}}},\quad 
 \Phi_\chi\bigl(\mathcal D(\widehat\rho_v^j(b_T))\bigr) =\chi_{q,t}(L(M_{\operatorname{pr}^j(T)})).
\end{aligned}
\end{equation}
In particular, specializing the unnormalized shifted element gives
\[
 \Phi_\chi\bigl(\pi_\chi(\widehat\rho_v^j(b_T))\bigr)
 =t^{\chi(\operatorname{pr}^j(T))/2}
 \chi_{q,t}(L(M_{\operatorname{pr}^j(T)})).
\]

For the classical case, denote by $b_T^{\rm cl}$ the specialization
of $b_T$ at $v^{1/2}=1$, and let $\pi_{\rm cl}$ be the quotient map
in~\eqref{eq:classical-isomorphism}, from the homogeneous coordinate
ring to $\mathbb C[\Gr(k,n,\sim)]$. Write
$\operatorname{ch}(T)=\pi_{\rm cl}(b_T^{\rm cl})$ for the basis element
in this quotient, using the notation of~\cite{CDFL20}.
For the classical cyclic shift $\rho$ of Section~\ref{subsec:normalized-cyclic-shift}, we have
\[
 \pi_{\rm cl}(\rho^j(b_T^{\rm cl}))
 =\operatorname{ch}(\operatorname{pr}^j(T))
 =\operatorname{ch}((\operatorname{pr}^j(T))_{\mathrm{red}}).
\]
Thus we may first apply the quantum cyclic shift to $b_T$, then
specialize $v^{1/2}=1$ and impose $P_{C_a}=1$ for the original
non-wrapping consecutive columns. The result is the same as first
specializing $v^{1/2}=1$, applying the classical cyclic shift on
the homogeneous coordinate ring, and then imposing these same
relations.

The cyclic shift must be applied before taking the quotient.
The obstruction to descent is already present in the classical case.
Let $\mathcal I_{\rm cl}=\langle P_{C_a}-1:1\leq a\leq s\rangle$
and $B_{\rm cl}=A_{\rm cl}/\mathcal I_{\rm cl}$.
The map $\rho$ induces an automorphism of $B_{\rm cl}$ if and only if
$\rho(\mathcal I_{\rm cl})=\mathcal I_{\rm cl}$.
However, $\rho(P_{C_s}-1)=P_{\sigma C_s}-1$, where
$C_s=[s,n]$ and $\sigma C_s=\{1,s+1,\ldots,n\}$ is a wrapping
consecutive column. It is not a trivial column and is itself reduced.
Consequently $\pi_{\rm cl}(P_{\sigma C_s})\ne 1$ by the classical
basis of the quotient, so $\rho$ does not preserve $\mathcal I_{\rm cl}$.
In the quantum case, the same obstruction is expressed by
\[
 \mathcal D(P_{C_s})=1,\qquad
 \mathcal D(\widehat\rho_v(P_{C_s}))=\mathcal D(P_{\sigma C_s})\ne 1,
\]
where the inequality follows from Proposition~\ref{thm:free}.
Equivalently, $P_{C_s}-v^{\chi(C_s)/2}$ belongs to
$\ker\pi_\chi$, but its image under $\widehat\rho_v$ does not.
In particular, there is no linear map on $B_\chi$ making
$\pi_\chi\circ\widehat\rho_v$ factor through $\pi_\chi$.
Thus the cyclic shift does not descend to the central quotient, as stated
in Theorem~\ref{thm:intro-C}(i).

One can select the unique reduced representative of each tableau
class and assign $\widehat b_T$ to
$\widehat b_{(\operatorname{pr}(T))_{\rm red}}$ for reduced $T$.
This assignment, extended linearly, is not induced by the quantum
cyclic shift on arbitrary representatives and is not injective. Indeed, the wrapping
consecutive column $\sigma^{-1}C_1$ is reduced, but its shift is
$C_1$, which represents the unit, just as the empty tableau does.
The same argument applies to the classical basis. Hence even
choosing reduced representatives does not give a cyclic action
by automorphisms of the quotient.

\begin{example}
Consider the tableau $T=(124\mid356)$ for $\Gr(3,6)$, as in
Example~\ref{ex:kl-Gr36}, where the vertical bar separates the columns.
Removing $6$ and applying reverse jeu de
taquin moves the entries $5$, $3$, and $1$ successively into the
empty box. Increasing the remaining entries by one and filling
the last empty box with $1$ gives
\[
 T=\begin{ytableau}
 1 & 3\\
 2 & 5\\
 4 & 6
 \end{ytableau},\qquad
 \operatorname{pr}(T)=\begin{ytableau}
 1 & 2\\
 3 & 4\\
 5 & 6
 \end{ytableau}.
\]
Shifting the columns separately gives $(235\mid146)$, which is not
semistandard in that order. The KL calculation in
Example~\ref{ex:kl-Gr36} and~\eqref{eq:kl-standard-monomial-expansion}
give $\widetilde P_T=b_T+v^{-1}b_{(123\mid456)}$.
Since $b_{(123\mid456)}=\widetilde P_{(123\mid456)}$ and both
tableaux have normalization exponent $\beta=3$ by~\eqref{eq:main-beta},
we obtain
\[
 b_T=v^{-3/2}\bigl(P_{124}P_{356}-v^{-1}P_{123}P_{456}\bigr)
\]
in $A_v$. By~\eqref{eq:cyclic-normalized-product},
\[
 \widehat\rho_v(P_{124}P_{356})=vP_{235}P_{146},\qquad
 \widehat\rho_v(P_{123}P_{456})=vP_{234}P_{156}.
\]
Consequently,
\[
 \widehat\rho_v(b_T)
 =v^{-1/2}\bigl(P_{235}P_{146}-v^{-1}P_{234}P_{156}\bigr)
 =b_{(135\mid246)},
\]
where the last equality follows from~\eqref{eq:cyclic-canonical-basis}.
Here $d=2$ and $r=1$, so
$\rho_v^{\rm LL}(b_T)=v^{-3}b_{(135\mid246)}$ by
\eqref{eq:cyclic-LL-canonical-basis}. We have
$M_T=Y_{1,-1}Y_{2,-4}$ and
$M_{\operatorname{pr}(T)}=Y_{2,0}Y_{1,-3}Y_{2,-2}Y_{1,-5}$.
By~\eqref{eq:cyclic-tableau-images}, the shifted element has image
\[
 \Phi_\chi(\mathcal D(\widehat\rho_v(b_T)))
 =\chi_{q,t}(L(Y_{2,0}Y_{1,-3}Y_{2,-2}Y_{1,-5})).
\]
\end{example}

\subsection{The quantum cluster structure of the localized quotient}
\label{subsec:cyclic-quasi-automorphism}

We keep the cocycle $\Omega$ and the specialization of the
non-wrapping consecutive minors fixed. We first establish the quantum
cluster structure of the quotient localized at the remaining cyclic
consecutive minors.
For $1\leq a\leq n$, write $F_a=P_{C_a}$, where $C_a$ is now
interpreted cyclically, and put $z_a=\mathcal D(F_a)$.
Thus $z_a=1$ for $1\leq a\leq s$. The remaining $z_a$ are
nonzero by Proposition~\ref{thm:free}. Cyclic consecutive minors
quasi-commute with every quantum Pl\"ucker coordinate, and this
property is preserved by the cocycle twist and central specialization.
Since $B_\chi$ is a domain by Theorem~\ref{thm:general-based-isomorphism},
these elements are regular normal elements and their powers generate
an Ore set. We write
$B_\chi^{\mathrm{loc}}=B_\chi[z_{s+1}^{-1},\ldots,z_n^{-1}]$ and $z=z_{s+1}$.
Since the $z_a$ are bar invariant, the modified bar extends uniquely to
this localization with $\overline{z_a^{-1}}^{\,\chi}=z_a^{-1}$ for $s<a\leq n$.
For pairwise quasi-commuting bar-invariant elements $x_1,\ldots,x_r$ with
$x_ix_j=v^{\lambda_{ij}}x_jx_i$, the normalized product is
\[
 [x_1^{a_1}\cdots x_r^{a_r}]
 =v^{-\frac12\sum_{i<j}a_ia_j\lambda_{ij}}x_1^{a_1}\cdots x_r^{a_r},
\]
where negative exponents are allowed for invertible factors.
For a quantum seed with $N$ variables and $r$ mutable variables
(listed first), its extended exchange matrix $\widetilde B$ is an
$N\times r$ integer matrix, with principal part $B$ consisting of
the first $r$ rows. The skew-symmetric commutation matrix $\Lambda$
satisfies $\widetilde B^{\mathrm T}\Lambda=(D\ \ 0)$ for a positive
diagonal matrix $D$. Writing $X^a=[x_1^{a_1}\cdots x_N^{a_N}]$
in its seed torus, the exchange relation in direction $i$ replaces
$x_i$ by $X^{-e_i+[b_i]_+}+X^{-e_i+[-b_i]_+}$, where $b_i$ is the
$i$th column of $\widetilde B$ and $[\cdot]_+$ is the entrywise
positive part; see~\cite{BZ05}.
\begin{proposition}\label{prop:localized-cluster-structure}
The algebra $B_\chi^{\mathrm{loc}}$ is a quantum cluster algebra,
with seeds obtained from the Grassmannian quantum seeds by the
cocycle twist and the specialization $F_1=\cdots=F_s=1$ via
$\mathcal D$. The extended exchange matrices are obtained by deleting
the rows of the specialized frozen variables.
\end{proposition}

\begin{proof}
Use the quantum Grassmannian cluster structure of
\cite[Theorem~1.3]{JKS2022}, over the generic coefficient field
(the quantum minor seeds and mutations are defined over $\bbk$).
Choose a seed of quantum Pl\"ucker coordinates $x_i$, with extended
exchange matrix $\widetilde B$ and commutation matrix $\Lambda$,
so that $x_ix_j=v^{\Lambda_{ij}}x_jx_i$. Let $G$ have columns
$\alpha_i=\deg x_i$. The content grading is a cluster grading:
$G\widetilde B=0$, as follows from homogeneity of the Grassmannian
exchange relations; see~\cite[Section~2.5]{CDFL20}. Under mutation
at $i$, the new degree is $G(-e_i+[b_i]_+)$, where $b_i$ is the
$i$th column of $\widetilde B$. The usual matrix mutation rule
then gives the grading identity at every mutated seed.

The twisted commutation matrix is
$\Lambda^\Omega=\Lambda+G^{\mathrm T}\Omega G$.
In particular,
$\widetilde B^{\mathrm T}\Lambda^\Omega=\widetilde B^{\mathrm T}\Lambda$.
A normalized Laurent monomial in the twisted seed has the same
underlying vector as the corresponding normalized Laurent monomial
in the original seed: the scalar factors in the twisted products
cancel the change in the commutation matrix. This follows by
bilinearity, also for negative powers in the seed torus. Thus the
two terms in each quantum mutation, which have the same content,
are unchanged as vectors by twisting. Induction on the mutation
sequence shows that the twisted algebra with all cyclic consecutive
minors inverted is the quantum cluster algebra of the twisted seed.

Let $p$ delete the $s$ rows corresponding to $F_1,\ldots,F_s$.
Their rows and columns in $\Lambda^\Omega$ vanish by
Proposition~\ref{prop:center}. Consequently the reduced matrices
$\widetilde B_0=p\widetilde B$ and $\Lambda_0=p\Lambda^\Omega p^{\mathrm T}$
are compatible, with the same positive diagonal block as before.

We justify that the specialization produces the claimed seed,
including the independence of its variables. Let $A^*$ be the twisted
algebra with all cyclic consecutive minors inverted, and $\mathcal T$
the full twisted seed torus. By the quantum Laurent phenomenon,
$A^*\subseteq\mathcal T$. The central ideal
$J=\langle F_a-1:1\leq a\leq s\rangle$ in $A^*$ extends to the
ideal $J\mathcal T$ of $\mathcal T$. The quotient of $\mathcal T$
by this ideal is the quantum torus of $\Lambda_0$, since the
specialized variables are central independent torus generators.
The quotient $A^*/J$ is $B_\chi^{\mathrm{loc}}$ via $\mathcal D$.
If $f\in A^*\cap J\mathcal T$, express it as a finite sum
$\sum_a(F_a-1)h_a$ with $h_a\in\mathcal T$. Multiply on the right
by a seed monomial $u$ with sufficiently large nonnegative exponents
so that every $h_a u$ is a polynomial in the seed generators.
Then $fu\in J$. The image of $u$ is nonzero in $A^*/J$, because
every seed Pl\"ucker coordinate has nonzero image by
Theorem~\ref{thm:general-based-isomorphism}, and the quotient is
a domain. Hence $f\in J$. This proves that $A^*/J$ embeds in the
reduced seed torus.

Normalized Laurent monomials specialize to the normalized
monomials of the reduced torus, with exponent vector projected by
$p$. Quantum mutation therefore commutes with the specialization.
More explicitly, deleting frozen entries commutes with taking the
positive and negative parts of an exchange column on the retained
indices, and the deleted variables have value $1$. Inductively the
specialized cluster variables are exactly the variables obtained by
mutation from the reduced seed (in particular they are nonzero).
Since $A^*$ is generated by its cluster variables and invertible
frozen variables, the specialized cluster algebra equals
$B_\chi^{\mathrm{loc}}$.
\end{proof}

\subsection{The corrected quantum cyclic shift}

We obtain a map on the quotient by correcting the cyclic shift
with normalized frozen Laurent factors, and establish its action on
the localized dual canonical basis.
For quantum cluster quasi-homomorphisms, we use the
definition in~\cite[Definition~2.1]{CHLQuasi2024}.
See also~\cite[Section~3]{KQW2023} for the related
theory of variation maps and quantum twist endomorphisms.
A quantum cluster quasi-homomorphism between seeds with the same number $r$
of mutable variables and principal parts $B=B'$ is an algebra
homomorphism $f$ satisfying
\[
 f(x_i)=X'^{e_i+h_i}\ (1\leq i\leq r),\qquad
 f(x_j)=X'^{l_j}\ (j>r),\qquad f(X^{b_i})=X'^{b'_i}\ (1\leq i\leq r),
\]
where $h_i,l_j$ are integer exponent vectors supported on the
frozen indices of the target seed, and $b_i,b'_i$ are the respective
exchange columns. The monomials $X^{b_i}$ and $X'^{b'_i}$ are the
quantum cluster $X$-variables. A quasi-automorphism is a
quasi-homomorphism from a cluster algebra to itself admitting another
quasi-homomorphism whose two compositions with it agree with the
identity on mutable cluster variables up to frozen Laurent monomials
and powers of $v^{1/2}$.
In our construction this quasi-inverse is an actual inverse.

For homogeneous elements of contents $\alpha,\beta\in\mathcal L$,
set $h(\alpha)=\alpha_n$ and define the skew-symmetric integral form
\begin{equation}\label{eq:cyclic-multiplicative-defect}
 E(\alpha,\beta)=\Omega(\alpha,\beta)-\Omega(\sigma\alpha,\sigma\beta)
             +2(d_\alpha\beta_n-\alpha_n d_\beta).
\end{equation}
Then~\eqref{eq:star} and~\eqref{eq:cyclic-normalized-product} give
$\widehat\rho_v(f\star g)=v^{E(\alpha,\beta)/2}
\widehat\rho_v(f)\star\widehat\rho_v(g)$.
Put $c=\deg F_s$ and $\lambda(\alpha)=E(c,\alpha)$.
Since $c_n=1$ and $E(\alpha,\beta)=0$ when $\alpha_n=\beta_n=0$
(the restriction of $\Omega$ to the first $n-1$ coordinates is
invariant under simultaneous translation by one), decomposing
$\alpha$ as $(\alpha-\alpha_nc)+\alpha_nc$ gives
\begin{equation}\label{eq:cyclic-defect-factorization}
 E(\alpha,\beta)=h(\alpha)\lambda(\beta)-h(\beta)\lambda(\alpha).
\end{equation}

\begin{theorem}\label{thm:cyclic-quasi-automorphism}
There is a unique quantum cluster quasi-automorphism $\widetilde\rho_v$
of $B_\chi^{\mathrm{loc}}$ satisfying, for every $k$-element subset
$I\subseteq[n]$ and every rectangular semistandard tableau $T$,
\begin{equation}\label{eq:cyclic-quasi-formula}
\begin{aligned}
 \widetilde\rho_v(\mathcal D(P_I))
 = [z^{-\mathbf1_{n\in I}}\mathcal D(P_{\sigma I})],\quad 
 \widetilde\rho_v(\mathcal D(b_T))
 = [z^{-r_n(T)}\mathcal D(b_{\operatorname{pr}(T)})],
\end{aligned}
\end{equation}
where $r_n(T)$ is the number of entries equal to $n$ in $T$.
In particular, the second formula is independent of the choice
of homogeneous tableau representative of a quotient basis element.
Uniqueness follows because the images of the Pl\"ucker coordinates
and the inverses of the remaining frozen variables generate the
localized algebra.
\end{theorem}

\begin{proof}
For homogeneous $f\in A_v^\Omega$ of content $\alpha$, set
$A_f=\mathcal D(\widehat\rho_v(f))$. Since $F_s$ is central in
$A_v^\Omega$, the product formula with $E$ gives
$zA_f=v^{-\lambda(\alpha)}A_fz$. Define
$\Psi(f)=v^{-h(\alpha)\lambda(\alpha)/2}z^{-h(\alpha)}A_f$
and extend linearly. This is the normalized product
$[z^{-h(\alpha)}A_f]$ if $f$ is bar invariant.
To check multiplicativity, let $a=h(\alpha)$, $b=h(\beta)$,
$u=\lambda(\alpha)$ and $w=\lambda(\beta)$ for homogeneous $f,g$.
Using $A_f z^{-b}=v^{-bu}z^{-b}A_f$, we get
$\Psi(f)\Psi(g)=v^{-au/2-bw/2-bu}z^{-a-b}A_f A_g$.
By~\eqref{eq:cyclic-defect-factorization}, the same expression is
$\Psi(f\star g)$, since
$-(a+b)(u+w)/2+(aw-bu)/2=-au/2-bw/2-bu$.
Thus $\Psi$ is a unital algebra homomorphism.
For $a<s$, $F_a$ has no entry $n$, hence $\Psi(F_a)=1$.
For $F_s$, use $h(c)=1$ and $\lambda(c)=0$ to obtain
$\Psi(F_s)=z^{-1}z=1$. The map factors through $\mathcal D$,
whose kernel is generated by $F_a-1$, $1\leq a\leq s$.
For $s+1\leq a<n$ it sends $z_a$ to $[z^{-1}z_{a+1}]$,
and sends $z_n$ to $z^{-1}$. These images are invertible, so
it extends to an endomorphism of $B_\chi^{\mathrm{loc}}$.
The formulas~\eqref{eq:cyclic-quasi-formula} follow from the action
of $\widehat\rho_v$ on quantum Pl\"ucker coordinates and from
Theorem~\ref{thm:cyclic-canonical-basis}.

For the quantum cluster quasi-automorphism property
(\cite[Definition~2.1]{CHLQuasi2024}), choose a Pl\"ucker seed
from a maximal weakly separated collection and its cyclic shift.
Here two $k$-subsets $I,J$ are weakly separated if there are no
four cyclically ordered distinct indices $a,b,c,d$ with
$a,c\in I\setminus J$ and $b,d\in J\setminus I$.
Both are seeds of the Grassmannian cluster algebra by~\cite{Scott2006,JKS2022}.
Rotation preserves their extended quivers, when the vertices
are indexed by the original labels; see also~\cite[Section~5]{FraserBraid2020}.
Let $\widetilde B$ be this common full exchange matrix, and
$\Lambda^\Omega,\Lambda'^\Omega$ the full twisted commutation
matrices. Proposition~\ref{prop:localized-cluster-structure}
provides the specialized seeds. Index frozen vertices by their
labels $F_a$ before rotation. Let $p$ delete $F_1,\ldots,F_s$,
and let $p'$ delete $F_n,F_1,\ldots,F_{s-1}$. These are exactly the
vertices specialized to $1$ in the source and target respectively.
In particular the principal exchange matrices coincide.

Let $h=(h(\alpha_i))_i$ be the row of last content coordinates
of the full source seed and $e_s$ the unit vector at its vertex
$F_s$. The exponent map before deleting rows is
$Q=I-e_s h$. For a seed variable with exponent $e_i$, its image
has normalized exponent $p'Qe_i$. For each deleted source
vertex $F_a$, $p'Qe_a=0$: if $a<s$ then $h(\deg F_a)=0$,
and if $a=s$ then $Qe_s=0$. Hence $R$ defined by $Rp=p'Q$
is a well-defined integer matrix on the reduced seed lattices.
The grading identity gives $h\widetilde B=0$, so
\begin{equation}\label{eq:cyclic-quasi-exchange-matrix}
 R\widetilde B_0=Rp\widetilde B=p'Q\widetilde B
 =p'\widetilde B=\widetilde B'_0.
\end{equation}

The commutation formula with $E$ implies
$\Lambda'^\Omega_{ij}=\Lambda^\Omega_{ij}-E(\alpha_i,\alpha_j)$.
Since the row of $F_s$ in $\Lambda^\Omega$ is zero, the
corresponding row in $\Lambda'^\Omega$ has entries
$-\lambda(\alpha_j)$. Writing $h_i=h(\alpha_i)$, we obtain
\[
 (Q^{\mathrm T}\Lambda'^\Omega Q)_{ij}
 =\Lambda'^\Omega_{ij}+h_i\lambda(\alpha_j)-h_j\lambda(\alpha_i)
 =\Lambda^\Omega_{ij}.
\]
The deleted rows and columns are zero in the respective
commutation matrices, so this gives
\begin{equation}\label{eq:cyclic-quasi-commutation-matrix}
 R^{\mathrm T}\Lambda'_0 R=\Lambda_0.
\end{equation}
The images of the mutable seed variables are normalized products
of the corresponding rotated variables with frozen Laurent
monomials; the images of the frozen variables are frozen Laurent
monomials. Equations~\eqref{eq:cyclic-quasi-exchange-matrix}
and~\eqref{eq:cyclic-quasi-commutation-matrix} are exactly the
matrix conditions of~\cite[Proposition~2.8]{CHLQuasi2024}.

Finally $R$ has block form $\left(\begin{smallmatrix}I&0\\H_0&L_0\end{smallmatrix}\right)$,
where, with the remaining frozen variables ordered cyclically as
$w_i=z_{s+i}$ ($1\leq i\leq k-1$), the columns of $L_0$ are
$-e_1+e_{i+1}$ for $i<k-1$ and $-e_1$ for $i=k-1$.
Subtracting the last column from all earlier columns shows
$\det L_0=(-1)^{k-1}$. Thus $R^{-1}$ is integral, has the same block form and satisfies
the reversed matrix equations. By~\cite[Proposition~2.8]{CHLQuasi2024},
it induces a quantum cluster quasi-homomorphism inverse to
$\widetilde\rho_v$ on the seed tori and hence on the cluster algebras.
If there are no mutable variables, the same conclusion follows directly
from the isomorphism of the two quantum tori.
Thus $\widetilde\rho_v$ is a quantum cluster quasi-automorphism.
\end{proof}

In the classical case the same construction, without the
commutation correction, gives an automorphism of
$\mathbb C[\Gr(k,n,\sim)][z_{s+1}^{-1},\ldots,z_n^{-1}]$.
On the images of the homogeneous Pl\"ucker coordinates it is
$\widetilde\rho(\pi_{\rm cl}(P_I))=
 z^{-\mathbf1_{n\in I}}\pi_{\rm cl}(P_{\sigma I})$, where now
$z=\pi_{\rm cl}(F_{s+1})$. Indeed the content grading makes this
formula multiplicative before taking the quotient, the first $s$
frozen variables have image $1$, and the matrix $R$ above gives
an invertible coefficient transformation preserving the exchange
columns. It is a cluster quasi-automorphism in the sense of
\cite{FraserQuasi2016}. The corresponding formula on the quotient
basis is $\widetilde\rho(\operatorname{ch}(T))=
 z^{-r_n(T)}\operatorname{ch}(\operatorname{pr}(T))$.
The localization is essential for this construction. In fact,
the last remaining frozen variable is sent to $z^{-1}$, which
is not in the unlocalized algebra. To see the latter assertion
in the quantum case, use the filtration by the positive weights
$i(k-i)$ from the proof of Proposition~\ref{thm:general-defining}.
The correction terms in the defining relations have strictly
smaller weight, so the associated graded algebra is a quantum
affine space by the ordered monomial basis. It is a domain and
has degree-zero component $\bbk$. Thus every unit of $B_\chi$
is a scalar. Since $z$ is a nontrivial basis element, it is not
a unit. In the classical case the same assertion follows from
the polynomial presentation of the Grothendieck ring.

For brevity, write $\tau=\widetilde\rho_v$,
$X_T=\mathcal D(b_T)=\widehat b_{T_{\mathrm{red}}}$ and
$w_b=z_{s+b}$, $1\leq b<k$. In particular,
$\mathcal B_\chi=\{X_T:T\text{ is reduced}\}$. For $\eta\in\mathbb Z^{k-1}$,
let $w^\eta=[w_1^{\eta_1}\cdots w_{k-1}^{\eta_{k-1}}]$ denote the
normalized frozen Laurent monomial. The notation $[w^\eta X_T]$ is well defined since
the frozen variables quasi-commute with homogeneous basis elements.
Indeed, their commutation with every quantum Pl\"ucker coordinate
is determined by its content, and this property extends by
homogeneity and survives the twist and specialization.

\begin{proposition}\label{prop:cyclic-localized-basis}
The set
\[
 \mathcal B_{\mathrm{loc}}=\{[w^\eta X_T]:
 \eta\in\mathbb Z^{k-1},\ T\text{ reduced}\}
\]
(with repetitions omitted) is a bar-invariant basis of
$B_\chi^{\mathrm{loc}}$, called the localized dual canonical basis.
The map $\widetilde\rho_v$ permutes $\mathcal B_{\mathrm{loc}}$.
\end{proposition}

\begin{proof}
First, normalized multiplication by any $w_b$ takes a quotient
basis element to a quotient basis element. Indeed, every cyclic
consecutive column is $\sigma^j C_1$ for some $j$. Applying
$\widehat\rho_v^j$ to the product $F_1 b_{\operatorname{pr}^{-j}(T)}$
in~\eqref{eq:untwisted-frozen}, and using the multiplication formula
\eqref{eq:cyclic-normalized-product} and
Theorem~\ref{thm:cyclic-canonical-basis}, shows that
$P_{\sigma^j C_1}b_T$ is a power of $v^{1/2}$ times a dual canonical
basis element of $A_v$. Twisting and applying $\mathcal D$
preserve this assertion. The normalized product is bar invariant,
so the scalar power is removed by normalization. Consequently,
$[w^\eta X_T]$ is a quotient basis element whenever
$\eta\in\mathbb Z_{\geq0}^{k-1}$.

For $a\in\mathbb Z_{\geq0}^{k-1}$, the set of elements
$[w^{-a}X_T]$ with $T$ reduced is a basis of the vector space
$w^{-a}B_\chi$: left multiplication by the invertible $w^{-a}$
is injective, and normalization rescales each element by a nonzero
scalar. These sets of basis elements are nested with respect to the
coordinatewise order on $a$. Indeed, if $a\leq b$, then
$[w^{b-a}X_T]$ is a quotient basis element and
$[w^{-b}[w^{b-a}X_T]]=[w^{-a}X_T]$, since the two normalized
expressions differ at most by a power of $v^{1/2}$ and are both
bar invariant. The union of these bases is precisely
$\mathcal B_{\mathrm{loc}}$, and it is a basis because any finite
collection of elements of the Ore localization admits a common
frozen monomial denominator.

By~\eqref{eq:cyclic-quasi-formula}, $\tau$ sends the normalized
Pl\"ucker generators and their frozen inverses to bar-invariant
elements. It therefore commutes with the bar on the whole
localization and sends normalized frozen Laurent monomials
to normalized frozen Laurent monomials. They also show
that $\tau(X_T)$ belongs to $\mathcal B_{\mathrm{loc}}$.
As $\tau$ respects normalized products, it maps every element of
$\mathcal B_{\mathrm{loc}}$ to another element of this basis.
Since $\tau$ is an algebra automorphism, it permutes the basis.
\end{proof}

In particular, the corrected quantum cyclic shift preserves the
dual canonical basis after localization, not in general the
unlocalized dual canonical basis $\mathcal B_\chi$. For example,
$\tau(z_n)=z^{-1}\notin B_\chi$.

\subsection{The order of the corrected cyclic shift}

We determine the $n$th iterate of the corrected shift and its order
as an algebra automorphism.

\begin{proposition}\label{prop:cyclic-quasi-order}
Let $r_i(T)$ be the multiplicity of $i$ in $T$, and set
$c_b(T)=\sum_{i\equiv b\ ({\rm mod}\ k)}r_i(T)$, for $1\leq b\leq k$.
The corrected shift satisfies
\begin{equation}\label{eq:cyclic-quasi-nth}
 \tau^n(X_T)=[w^{a(T)}X_T],\qquad
 a(T)_b=c_{b+1}(T)-c_b(T)\quad(1\leq b<k).
\end{equation}
Its order as an algebra automorphism is
\begin{equation}\label{eq:cyclic-quasi-order}
 \operatorname{ord}(\tau)=
 \begin{cases}
 k,&n=k+1,\\
 kn,&n\geq k+2,\ \gcd(k,n)=1,\\
 \infty,&\gcd(k,n)>1.
 \end{cases}
\end{equation}
The same assertions hold for the classical corrected shift, with
ordinary products instead of normalized products.
\end{proposition}

\begin{proof}
The frozen-variable formulas in the proof of
Theorem~\ref{thm:cyclic-quasi-automorphism} give
$\tau(w^\eta)=w^{L\eta}$, where $L$ is the integer matrix determined
by $Le_i=e_{i+1}-e_1$ for $1\leq i<k-1$ and $Le_{k-1}=-e_1$.
Here $e_i$ denotes the standard basis of the frozen exponent lattice.
As in the proof of Proposition~\ref{prop:cyclic-localized-basis},
$\tau$ preserves normalized products, which can be grouped in any
order. Repeated application of
\eqref{eq:cyclic-quasi-formula} gives, for $j\geq1$,
\begin{equation}\label{eq:cyclic-quasi-iterate}
 \tau^j(X_T)=[w^{a_j(T)}X_{\operatorname{pr}^j(T)}],\qquad
 a_j(T)=-\sum_{h=0}^{j-1}L^{j-1-h}e_1\,r_n(\operatorname{pr}^h(T)).
\end{equation}
Indeed $a_{j+1}(T)=La_j(T)-e_1r_n(\operatorname{pr}^j(T))$.
No additional scalar occurs because the products are normalized.

For $0\leq h<n$ we have $r_n(\operatorname{pr}^h(T))=r_{n-h}(T)$,
and $\operatorname{pr}^n(T)=T$. Hence
$a_n(T)=-\sum_{i=1}^n r_i(T)L^{i-1}e_1$.
The vectors $f_i=L^{i-1}e_1$, $1\leq i\leq k$, are
$f_1=e_1$, $f_i=e_i-e_{i-1}$ for $2\leq i<k$, and $f_k=-e_{k-1}$.
They sum to zero, span the lattice, and have no other independent
linear relation. The matrix $L$ permutes them cyclically. Thus
$L$ has order $k$ and $I+L+\cdots+L^{k-1}=0$. This gives
\eqref{eq:cyclic-quasi-nth}. Adding any non-wrapping consecutive
column increases all the $c_b(T)$ by one, so $a(T)$ is unchanged
by replacing $T$ with its reduced representative.

Put $g=\gcd(k,n)$ and $p=k/g$. By
\eqref{eq:cyclic-quasi-nth},
\begin{equation}\label{eq:cyclic-quasi-block-iterate}
 \tau^{np}(X_T)=[w^{(I+L^n+\cdots+L^{(p-1)n})a(T)}X_T].
\end{equation}
Moreover, $\tau^{np}$ fixes the frozen Laurent monomials because
$k$ divides $np$.
If $g=1$, the sum of matrices is zero, so $\tau^{kn}=1$ on
all the basis elements and frozen inverses that generate the algebra.
For $n=k+1$, every Pl\"ucker coordinate is frozen; the localized
algebra is generated by its remaining frozen Laurent monomials,
whose exponent transformation has order $k$. This gives the first
case of~\eqref{eq:cyclic-quasi-order}.

Suppose $n\geq k+2$ and $g=1$. The order of $\tau$ is a multiple
of $k$ by its action on the frozen exponent lattice, which embeds
in the seed torus by Proposition~\ref{prop:localized-cluster-structure}.
It is also a multiple of $n$. To see this, specialize first to
$v^{1/2}=1$ and then set the remaining frozen variables equal to $1$.
The corrected shift preserves these latter relations in the classical
algebra and induces the cyclic permutation of the Pl\"ucker labels
on the quotient by all cyclic consecutive minors minus one.
This comparison of orders is valid because the normalized products
in~\eqref{eq:cyclic-quasi-formula} and all their iterates on the
Pl\"ucker generators are defined over $\mathbb Z[v^{1/2},v^{-1/2}]$
in the localized integral form. Equivalently, one can specialize
these identities inside the integral Laurent seed torus, where
normalized Laurent monomials become ordinary Laurent monomials.
The specialization of an identity of iterates is therefore an
identity for the classical corrected shift, and then for the
all-frozen specialization.

Distinct non-frozen Pl\"ucker coordinates have distinct images in
the all-frozen specialization by~\cite[Corollary~3.3]{FeiginLi2026}:
they occur among the noncrossing monomials without frozen factors
that form its basis. Take $I_0=\{1,\ldots,k-1,k+1\}$, which is not
cyclic consecutive when $n\geq k+2$. Its cyclic orbit has length
$n$ since $k$ and $n$ are coprime: a nonidentity cyclic symmetry
would force its cardinality to be a multiple of the length of an
orbit of that symmetry, a nontrivial divisor of $n$.
It follows that any order of $\tau$ is divisible by $n$, hence by
$kn$. Together with the upper bound, this proves the second case.

Finally suppose $g>1$, so $n\geq k+2$. For the one-column tableau
$I_0$ just used, $a(I_0)=f_k-f_1$. With the subscripts on $f_i$
read modulo $k$, the exponent in~\eqref{eq:cyclic-quasi-block-iterate}
is
$b=\sum_{j=0}^{p-1}(f_{k+jn}-f_{1+jn})$.
It is nonzero: the positive and negative terms belong to distinct
residue classes modulo $g$, and the only linear relation among the
$f_i$ is generated by their sum. Since $\tau^{np}$ fixes the frozen
monomials, $\tau^{rnp}(X_{I_0})=[w^{rb}X_{I_0}]$ for every $r\geq1$.
These elements are distinct, since $X_{I_0}\ne0$, the algebra embeds
in the seed torus's division ring, and the frozen Laurent monomials
of distinct exponents are linearly independent. Thus $\tau$ has
infinite order. The proof for the classical map is identical,
using its commutative Laurent seed torus.
\end{proof}

In particular, $\tau^n$ sends every mutable variable of a Pl\"ucker
seed to itself up to a normalized frozen Laurent monomial. Thus it
is equivalent to the identity as a cluster quasi-automorphism in
the sense of~\cite{CHLQuasi2024}, even when $\tau$ itself has
infinite order.

\subsection{The corrected shift on the quantum Grothendieck ring}

We identify the localized quantum cluster structures under the main
isomorphism and transport the corrected cyclic shift, its action on
the localized canonical basis, and its order to the quantum
Grothendieck ring. We put
\[
 Z_i=\chi_{q,t}\left(L\left(\prod_{b=1}^{m}Y_{i,i-2b}\right)\right),
 \quad 1\leq i<k,\qquad
 K^{\mathrm{loc}}=K[Z_1^{-1},\ldots,Z_{k-1}^{-1}].
\]
For $1\leq a<k$, the wrapping consecutive column $C_{s+a}$
consists of $[1,a]$ and $[s+a,n]$. Its highest monomial is
$\prod_{b=1}^m Y_{k-a,k-a-2b}$ by~\eqref{eq:scaled-plucker-images}.
Thus $\Phi_\chi(z_{s+a})=Z_{k-a}$, and Theorem~\ref{thm:intro-A}
extends to an isomorphism $\Phi_\chi^{\mathrm{loc}}$ of the localizations.
The bar of $K$ extends to $K^{\mathrm{loc}}$ since the $Z_i$ are bar
invariant, with $\overline{Z_i^{-1}}=Z_i^{-1}$.

For clarity, this is also a quantum cluster isomorphism for the
usual Kirillov--Reshetikhin initial seed on $K^{\mathrm{loc}}$.
Recall that Kirillov--Reshetikhin modules have highest monomials of the
form $\prod_{a=0}^{r-1}Y_{i,p+2a}$. Here the rectangular Pl\"ucker seed after specialization consists
of the normalized minors with labels
$J_{i,r}=[1,k-i]\cup[k-i+r+1,k+r]$, with $1\leq i<k$ and $1\leq r\leq m$,
where $r=m$ labels the remaining frozen variables. By the column
formula for $M_J$, Theorem~\ref{thm:intro-A} gives
\[
 \Phi_\chi(\mathcal D(P_{J_{i,r}}))
 =\chi_{q,t}\left(L\left(\prod_{b=1}^{r}Y_{i,i-2b}\right)\right).
\]
The classical extended exchange matrices agree by
\cite[Theorem~2.1]{CDFL20}; the quantum seed on the Grothendieck
ring uses the same exchange matrix and the bar-normalized
Kirillov--Reshetikhin characters, by
\cite[Remark~8.4.2 and Theorems~8.4.3, 9.4.1]{Qin2017},
with the level category and height function as above.
With the commutation convention
$x_ix_j=t^{\Lambda_{ij}}x_jx_i$, the commutation matrices agree
because these are the images under an algebra isomorphism of
nonzero quasi-commuting elements. Consequently normalized seed
monomials and quantum mutations agree. This identifies the
localized cluster structures, not merely the underlying algebras.

The seed identification above, Theorem~\ref{thm:cyclic-quasi-automorphism},
and Propositions~\ref{prop:cyclic-localized-basis} and~\ref{prop:cyclic-quasi-order}
give the following corollary by conjugation. Since $\Phi_\chi$ intertwines
the bar involutions and sends $z$ to $Z_{k-1}$, it also preserves the
normalized products appearing in the formulas.

\begin{corollary}\label{cor:cyclic-quasi-grothendieck}
The map $\mathcal R_t=\Phi_\chi^{\mathrm{loc}}\circ\widetilde\rho_v
\circ(\Phi_\chi^{\mathrm{loc}})^{-1}$ is a quantum cluster
quasi-automorphism of $K^{\mathrm{loc}}$. For every tableau $T$,
\begin{equation}\label{eq:cyclic-quasi-characters}
 \mathcal R_t(\chi_{q,t}(L(M_T)))=
 [Z_{k-1}^{-r_n(T)}\chi_{q,t}(L(M_{\operatorname{pr}(T)}))].
\end{equation}
Here the normalized product uses the bar and the multiplication
in the localized quantum Grothendieck ring. Moreover, $\mathcal R_t$
permutes the localized canonical basis of $K^{\mathrm{loc}}$,
namely
\[
 \Phi_\chi^{\mathrm{loc}}(\mathcal B_{\mathrm{loc}})
 =\{[Z^\eta\chi_{q,t}(L(M))]:
 \eta\in\mathbb Z^{k-1},\ M\in\mathcal P_\ell^+\},
\]
where $Z^\eta=[Z_1^{\eta_1}\cdots Z_{k-1}^{\eta_{k-1}}]$ and
repetitions are omitted. The order of $\mathcal R_t$
is given by~\eqref{eq:cyclic-quasi-order}, and its $n$th iterate is
obtained from~\eqref{eq:cyclic-quasi-nth} by replacing $w_b$ with
$Z_{k-b}$.
\end{corollary}

This completes the proof of Theorem~\ref{thm:intro-C}(ii).

\begin{example}
For $\Gr(3,6)$, put $w_1=\mathcal D(P_{156})$ and
$w_2=\mathcal D(P_{126})$. The two remaining frozen variables
commute (by the frozen commutation formula and the twist).
The corrected shift acts by
$\widetilde\rho_v(w_1)=w_1^{-1}w_2$ and
$\widetilde\rho_v(w_2)=w_1^{-1}$. Thus the image of $P_{126}$ under
$\widetilde\rho_v\circ\mathcal D$ is no longer $1$, in contrast to
uncorrected specialization. For example,
$\widetilde\rho_v(\mathcal D(P_{124}))=\mathcal D(P_{235})$ and
$\widetilde\rho_v(\mathcal D(P_{346}))=v^{-1}w_1^{-1}\mathcal D(P_{145})$,
where $\lambda(\deg P_{346})=2$ by~\eqref{eq:cyclic-multiplicative-defect}.
For $T=(124\mid356)$ in the preceding example,
$r_6(T)=1$ and $\operatorname{pr}(T)=(135\mid246)$. Therefore
\[
\begin{aligned}
 \widetilde\rho_v(\mathcal D(b_T))&=[w_1^{-1}\mathcal D(b_{(135\mid246)})],\\
 \mathcal R_t(\chi_{q,t}(L(Y_{1,-1}Y_{2,-4})))
 &= [Z_2^{-1}\chi_{q,t}(L(Y_{2,0}Y_{1,-3}Y_{2,-2}Y_{1,-5}))].
\end{aligned}
\]
For $I=124$ the residue counts modulo $3$ are $(2,1,0)$.
Thus~\eqref{eq:cyclic-quasi-nth} gives
$\widetilde\rho_v^{6r}(\mathcal D(P_{124}))=[w_1^{-r}w_2^{-r}\mathcal D(P_{124})]$
for $r\geq1$, whereas the sixth power fixes both frozen variables. On the
representation-theoretic side this reads
\[
 \mathcal R_t^{6r}(\chi_{q,t}(L(Y_{1,-1})))
 =[Z_2^{-r}Z_1^{-r}\chi_{q,t}(L(Y_{1,-1}))].
\]
Consequently the corrected shift has infinite order, although the
homogeneous quantum cyclic shift has order six.
\end{example}

\end{document}